\documentclass[11pt,a4paper,reqno]{amsart}

\usepackage[colorlinks,bookmarksopen,bookmarksnumbered,citecolor=red,urlcolor=teal]{hyperref}
\usepackage{amsaddr}
\usepackage{epsfig,amsmath,amssymb,mathrsfs}
\usepackage{graphicx}
\usepackage{float} 
\usepackage{textcomp, eufrak}
\usepackage{verbatim}
\usepackage{fancyhdr}
\usepackage{BOONDOX-uprscr}

\usepackage[numbers]{natbib}
\usepackage[dvipsnames]{xcolor}

\usepackage{enumitem}
\usepackage{soul}
\usepackage{centernot}
\usepackage{textgreek}
\usepackage{caption}
\usepackage[mathscr]{euscript}
\usepackage{appendix}
\usepackage{footmisc}
\usepackage{scalerel}

\let\oldtocsection=\tocsection
 
\let\oldtocsubsection=\tocsubsection
 
\let\oldtocsubsubsection=\tocsubsubsection
 
\renewcommand{\tocsection}[2]{\hspace{0em}\oldtocsection{#1}{#2}}
\renewcommand{\tocsubsection}[2]{\hspace{1em}\oldtocsubsection{#1}{#2}}
\renewcommand{\tocsubsubsection}[2]{\hspace{2em}\oldtocsubsubsection{#1}{#2}}

\numberwithin{equation}{section}

\newtheorem{theorem}{Theorem}[section]
\newtheorem{definition}[theorem]{Definition}
\newtheorem{lem}[theorem]{Lemma}
\newtheorem{prop}[theorem]{Proposition}
\newtheorem{rem}[theorem]{Remark}
\newtheorem{cor}[theorem]{Corollary}

\numberwithin{equation}{section}
\numberwithin{table}{section}

\newcommand{\R}{\mathbb{R}}

\newcommand{\N}{\mathbb{N}}

\newcommand{\Rp}{\mathbb{R}^+}
\newcommand{\F}{\mathcal{F}}

\newcommand{\PS}{(\Omega, \mathcal{F}, \mathbb{P})}

\newcommand{\I}{\mathbb{I}}

\newcommand{\Bb}{\mathcal{B}}

\newcommand{\D}{\mathcal{D}}

\newcommand{\T}{\mathbb{T}}

\newcommand{\p}{\mathbb{P}}

\newcommand{\bn}{\begin{definition}}
\newcommand{\en}{\end{definition}} 
\newcommand{\bt}{\begin{theorem}}                
\newcommand{\et}{\end{theorem}}
 \newcommand{\bnm}{\begin{enumerate}}              
\newcommand{\enm}{\end{enumerate}}
\newcommand{\br}{\begin{rem}} 
\newcommand{\er}{\end{rem}}

\newcommand{\om}{\omega}

\newcommand{\Om}{\Omega}

\newcommand{\btm}{\begin{itemize}}
 \newcommand{\etm }{\end{itemize}}

\newcommand{\E}{\mathbb{E}}

\newcommand{\Rd}{\R^d}

\newcommand{\Ic}{\mathcal{I}}

\newcommand{\Scd}{\mathbb{S}^{d-1}}
\newcommand{\DK}{\mathcal{D}_{\textsc{kl}}}

\newcommand{\Solm}{\phi_{t_0,t}}
\newcommand{\LG}{\mathcal{L}}

\newcommand{\Xt}{\mathcal{M}}
\newcommand{\XXt}{\mathfrak{M}}
\newcommand{\Df}{\D_{\varphi}}
\newcommand{\ccdot}{{\,\cdot\,}}
\newcommand{\cccdot}{\scaleobj{1.25}{\ccdot}}
\newcommand{\Ms}{{\mathfrak{P}\hspace{.01cm}}}
\newcommand{\Ns}{\mathfrak{Y}}

\newcommand{\Ps}{{\mathfrak{P}\hspace{.01cm}}}

\newcommand{\Mt}{{\mathscr{M}\hspace{.02cm}}}
\newcommand{\St}{\hspace{.02cm}{\mathscr{S}\hspace{.02cm}}}

\newcommand{\XX}{\mathcal{X}}
\newcommand{\xxt}{\scaleobj{0.7}{\mathcal{X}}}

\newcommand{\brr}{\mathring{b}}
\newcommand{\Wd}{\mathcal{W}_d}
\newcommand{\Yt}{\tilde{X}_t}

\newcommand{\Sden}{\rho^{\hspace{0.03cm}\mu}_t}
\newcommand{\betaLk}{\varTheta_{\mu\nu}}

\newcommand{\tx}{{\hspace{.005em}\textsc{x}}}
\newcommand{\ty}{\textsc{y}}
\newcommand{\ve}{\varepsilon}
\newcommand\simbar[1]{\tilde{\bar{#1}}}
\newcommand{\fr}{\mathfrak{r}}

\newcommand{\fI}{\mathfrak{I}}
\newcommand{\fF}{\mathfrak{F}}

\newcommand{\PP}{\hspace*{.000cm}\mathord{\raisebox{-0.095em}{\scaleobj{.86}{\includegraphics[width=1em]{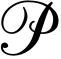}}}}\hspace*{.02cm}}

\newcommand{\mm}{{\hspace{0.02cm}\scaleobj{1.05}{\textrm{\textmugreek}}\hspace{0.01cm}}}
\newcommand{\mmo}{ {\hspace{0.02cm}\scaleobj{1.05}{\textrm{\textmugreek}}_{\textrm{\raisebox{-.15em}{\hspace{0.08em}$t_0$}}} }}
\newcommand{\mmt}{ {\hspace{0.02cm}\scaleobj{1.05}{\textrm{\textmugreek}}_{\textrm{\raisebox{-.15em}{\hspace{0.08em}$t$}}} \hspace{.1em}}}

\newcommand{\rrr}{\scaleobj{1.07}{\varrho} }
\newcommand{\rrt}{\scaleobj{1.07}{\varrho}^{\hspace{-0.02em}\mm}_t }
\newcommand{\rrto}{\scaleobj{1.07}{\varrho}^\mm_{t_0} }

\newcommand{\ff}{f}

\newcommand{\sff}{f}

\newcommand{\frat}{{\textstyle\frac{1}{2}}}
 \newcommand{\betaL}{\beta_{\mu\nu}}

\renewcommand\leftmark{{Lagrangian uncertainty quantification and information inequalities for stochastic flows}}

\renewcommand\rightmark{{Lagrangian uncertainty quantification and information inequalities for stochastic flows}}
\begin{document}
\title{\vspace*{-1.5cm} Lagrangian uncertainty quantification and information inequalities for stochastic flows}

\author{\vspace*{-.6cm} Micha\l \;Branicki$^{\,\dagger*}$ and Kenneth Uda$^{\dagger}$}

\address{\vspace*{-0.4cm} $^\dagger$ Department of Mathematics, University of Edinburgh, Scotland, UK \\
\hspace*{-1.9cm}$^*$ The Alan Turing Institute for Data Science, London, UK}

\email{M.Branicki@ed.ac.uk, K.Uda@ed.ac.uk}
\thanks{This work was supported by the Office of Naval Research grant ONR N00014-15-1-2351.}

\begin{abstract} 
 We develop a systematic information-theoretic framework for quantification and mitigation of error in probabilistic Lagrangian (i.e., path-based) predictions which are obtained from dynamical systems generated by uncertain (Eulerian) vector fields.  This work is motivated by the desire to improve Lagrangian  predictions in complex dynamical systems based either on analytically simplified or data-driven models. We derive a hierarchy of general information bounds on uncertainty in estimates of statistical observables  $\E^{\nu}[\ff]$,  evaluated on trajectories of  the approximating dynamical system, relative to the ``true'' observables  $\E^{\mu}[\ff]$  in terms of certain $\varphi$-divergences, $\Df(\mu\|\nu)$, which quantify discrepancies between probability measures $\mu$ associated with the original dynamics and their approximations $\nu$. We then derive two distinct bounds on $\Df(\mu\|\nu)$ itself in terms of the Eulerian fields.  This new framework provides a rigorous way for quantifying and mitigating uncertainty in Lagrangian predictions due to  Eulerian model error.

\medskip
\noindent\textbf{Keywords}: Lagrangian uncertainty quantification (LUQ), information inequalities, stochastic flows, $\varphi$-divergence, information theory, information geometry,  expansion rates.
 
\end{abstract}

\maketitle

\setcounter{tocdepth}{2}
\vspace*{.5cm}\tableofcontents

\vspace*{-0.3cm}\section{Introduction}\label{intro}
Given a probability space  $(\Om_\mm,\mathcal{H}_\mm, \Ms^\mm_{t_0})$,  consider a  dynamical system on a smooth finite-dimensional manifold $\XXt$ generated by a continuous map, $\phi^\mm_{t_0,t}\,{:}\;\XXt\,{\times}\, \Omega_\mm \rightarrow \XXt$, $t\in \Ic:=[t_0,t_0+T)$, such that $\phi^\mm_{t_0,t_0}(\ccdot,\om)=\textrm{id}_\XXt$, and $\Ms^\mm_{t_0}$ is the law of $\phi^\mm_{t_0,\scaleobj{1.1}{\ccdot}}$. Given the paths, $ t\mapsto \phi^\mm_{t_0,t}(\xxt,\om)$, $\om\in \Om_\mm$, labelled by the (potentially uncertain) initial conditions $\xxt\in \XXt$, and a measurable functional $f$ on these paths,  we refer to the estimation of {\it observables} $\xxt\mapsto \E\big[\ff\big(\phi^\mm_{t_0,t}(\xxt)\big)\big]$ as a {\it Lagrangian prediction}.  
This terminology follows from trajectory-based studies of transport in finite-dimensional dynamical systems (e.g., \cite{wiggins92,Shadden05,lekien07,wiggins06,macKay87,ottino89,rmked90,rypina20} for $f={\rm id}_\XXt$ or $f = \mathbb{I}_\mathscr{A}$, $\mathscr{A}\subseteq \Xt$, amongst a plethora of other publications) 
Many dynamical systems encountered in applications are  defined over high-dimensional manifolds, and they are   often generated by solutions of ordinary/stochastic differential equations (SDE/ODE) which, in turn, are often generated by (Eulerian) fields solving   partial differential equations;  
 examples range from fluid dynamics, to neural networks, to systems biology and molecular dynamics. 
 Approximations of the original dynamics 
 result in a loss of information, making the subsequent Lagrangian predictions  uncertain and often unreliable in ways which are difficult~to~assess directly from pointwise errors in the  Eulerian~fields generating the original and the approximate dynamics.  

 In this work we focus on developing a framework for {\it Lagrangian Uncertainty Quantification} (LUQ) which is concerned with a probabilistic quantification and mitigation of the error in  estimates of path-based observables (including the fate of the trajectories for $f={\rm id}_\XXt$), and which arises from approximations of the dynamics and uncertainties in the initial conditions. This task requires the study  of (path space) probability measures $\Ms^\mm_{t_0}$ and  the evolution of their time-marginal measures $\mmo\mapsto \mmt$, $\Ms^\mm_{t_0}\circ (\phi^\mm_{t_0,t})^{-1}=\mm_t$, where $\mmo\in \PP(\XXt)$ is a probability measure on the initial conditions $\xxt\in \XXt$. It is worth stressing from the outset that a number  of  results derived in what follows does not rely on Markovianity of the  underlying dynamics, or on the SDE/ODE formulation,  although this  is perhaps the most tractable and practically useful~setup.

 In order to focus attention and outline the LUQ framework, consider the original/reference dynamics induced by an  SDE in the Stratonovich form\footnote{We start from the Stratonovich form of the SDE rather than the It\^o form, since the former is consistent with the physical limit which leads to idealised stochastic perturbations in the deterministic dynamics (e.g., \cite{Gardiner10}).} on a smooth manifold $\XXt$ given by 
\begin{align}
 dX_t^\mm &= \textstyle \mathring{b}^{\mm}(t,X_t^\mm)dt+ \sigma^{\mm}(t,X_t^\mm)\circ d  W^\mm_{t-t_0},\qquad X_{t_0}^\mm \sim   \mmo, \quad t\in \Ic, \label{SDE1}
\end{align} 
where  the (Eulerian) fields $\xxt\mapsto \mathring{b}^\mm(\,\cdot\,,\xxt)$ and $\xxt\mapsto\sigma^\mm(\,\cdot\,,\xxt)$ generate the dynamics of (\ref{SDE1}) with  solutions  $X_t^\mm(\om) \,{=}\, \phi^\mm_{t_0,t}(\xxt,\om)$ in $\XXt$, and the uncertain initial condition $\xxt\in \XXt$ is distributed according  $\mmo\in \PP(\XXt)$;  $W_t^\mm$ is a
the Wiener processes of an appropriate dimension. Similarly, the approximate dynamics of (\ref{SDE1}) on a linear subspace $\Xt \subseteq \mathfrak{M} $ is induced~by  
\begin{align}
 dX_t^\nu &= \textstyle \mathring{b}^{\nu}(t,X_t^\nu)dt+ \sigma^{\nu}(t,X_t^\nu)\circ d W^\nu_{t-t_0},\qquad X_{t_0}^\nu \sim \nu_{t_0}, \quad t\in \Ic, \label{SDE2}
\end{align}
 with solutions $X_t^\nu(\om) \,{=}\, \phi^\nu_{t_0,t}(x,\om)$  in $\Xt$,  $\nu_{t_0}\in \PP(\Xt)$, and $W^\nu_t$ independent  of $W^\mm_t$. 

\begin{figure}
\captionsetup{width=1\linewidth}

\vspace{-.4cm}\includegraphics[height = 7.9cm]{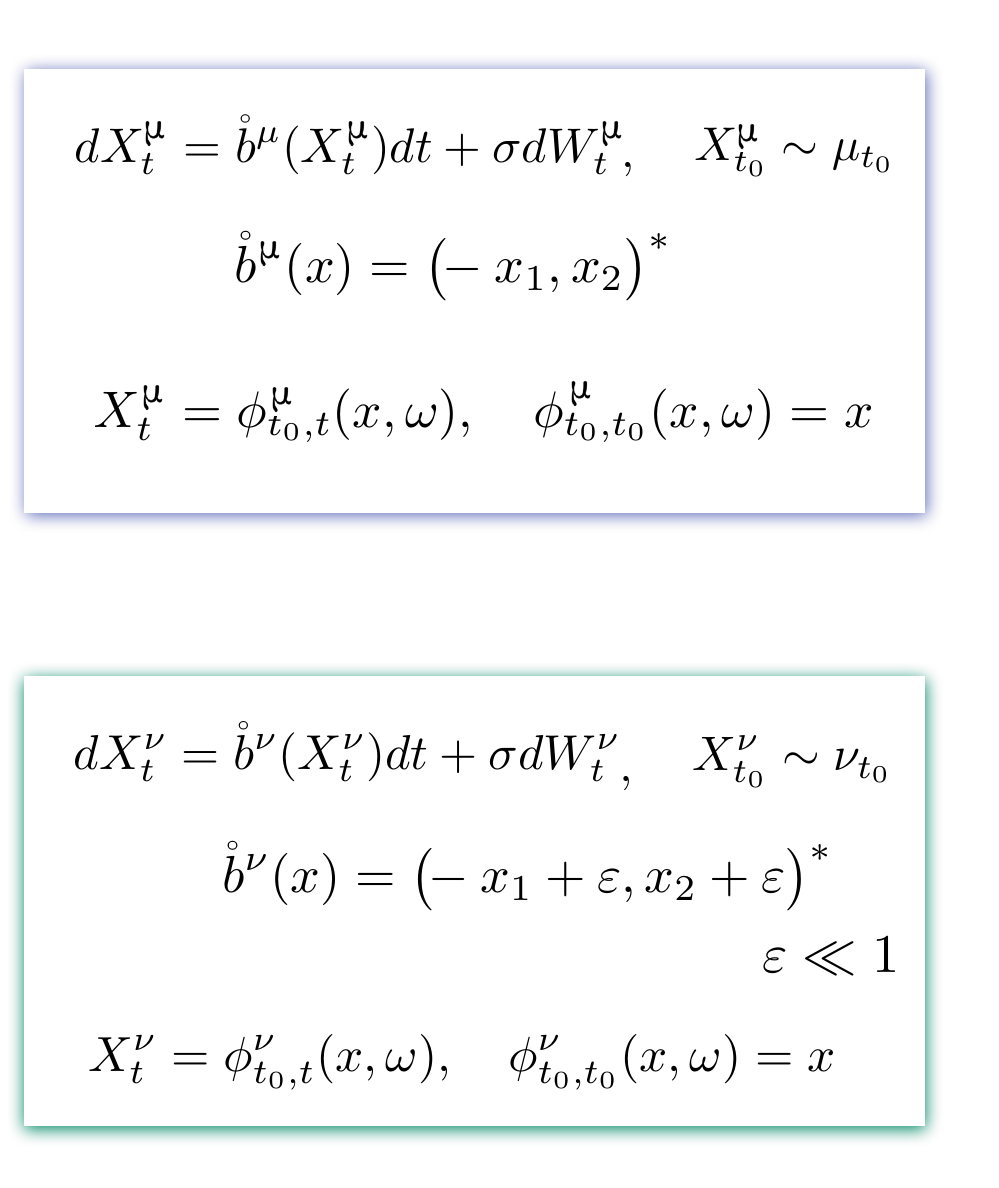}\hspace{.8cm}\includegraphics[height = 8.3cm]{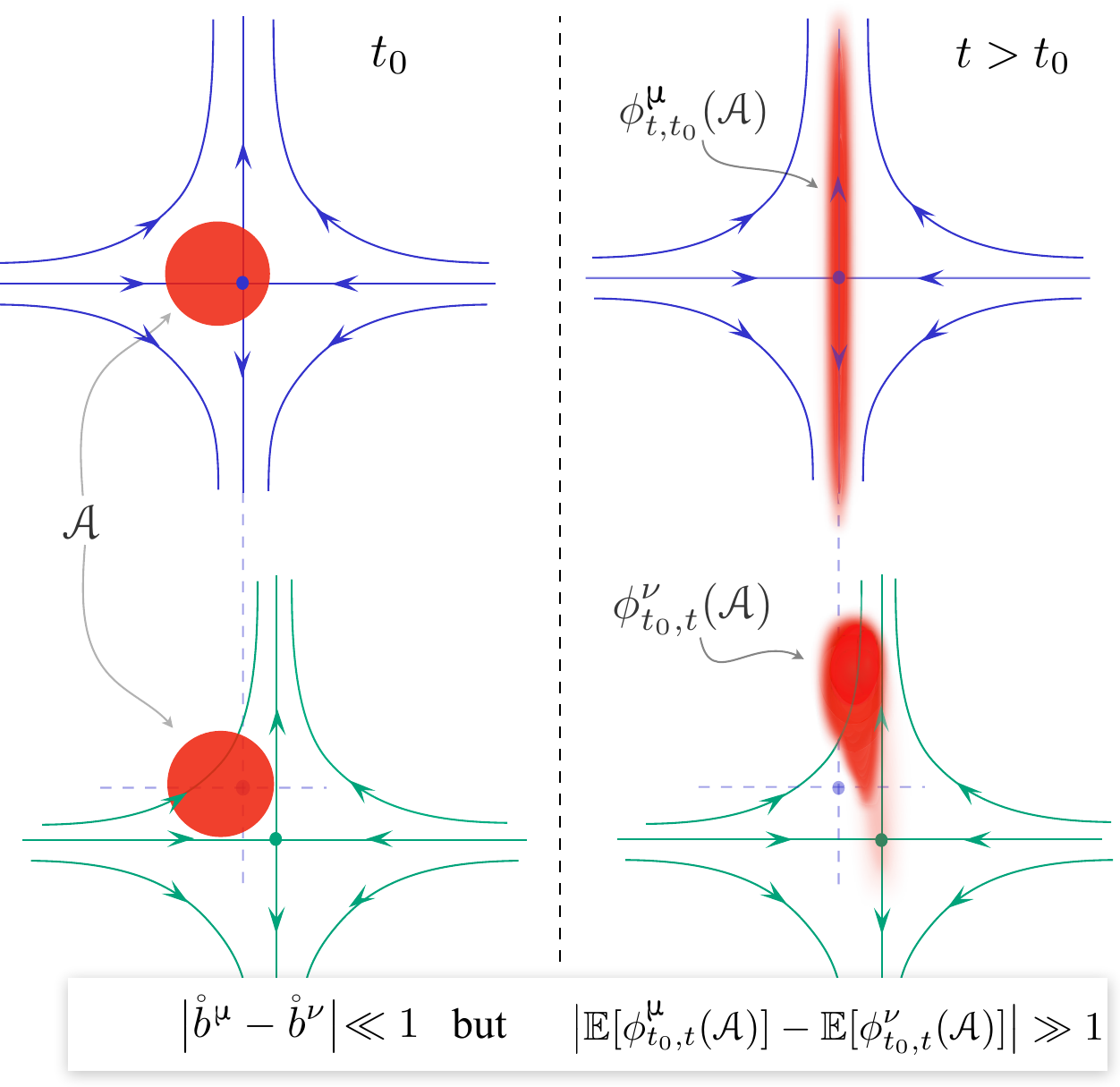}

\vspace{-0.45cm}\caption{\footnotesize
A simple  illustration of  differences between Eulerian and Lagrangian predictions.  
The discrepancy in the structure of trajectories starting from the same set  of initial conditions~$\mathscr{A}$, and  evolved under the dynamics $\phi^\text{\textmugreek}_{t_0,t}(\mathscr{A},\cdot\,)$, $\phi^\nu_{t_0,t}(\mathscr{A},\cdot\,)$,  can become large even if $|\brr^\text{\textmugreek}-\brr^\nu| \ll 1$, $\sigma^\text{\textmugreek} = \sigma^\nu$, $\forall\, t\geqslant t_0$. This is a vastly simplified sketch (for small diffusion $\sigma$) of a situation in which  one compares/predicts the evolution of a patch of a  passive tracer under the velocity field generated by, e.g., two oceanographic models.  Suitable bounds for quantification and mitigation of uncertainty in such predictions are outlined in \S\ref{out_main}, and a detailed worked example associated with a different setting of dimensionally-reduced approximations of a stochastic slow-fast dynamics is discussed~in~\S\ref{toysec}.}\label{motiv_1}\vspace*{-.3cm}
\end{figure}

\begin{figure}[t]
\centering
\captionsetup{width=1\linewidth}
\vspace{-0.7cm}\includegraphics[height = 9.5cm, width = 15cm]{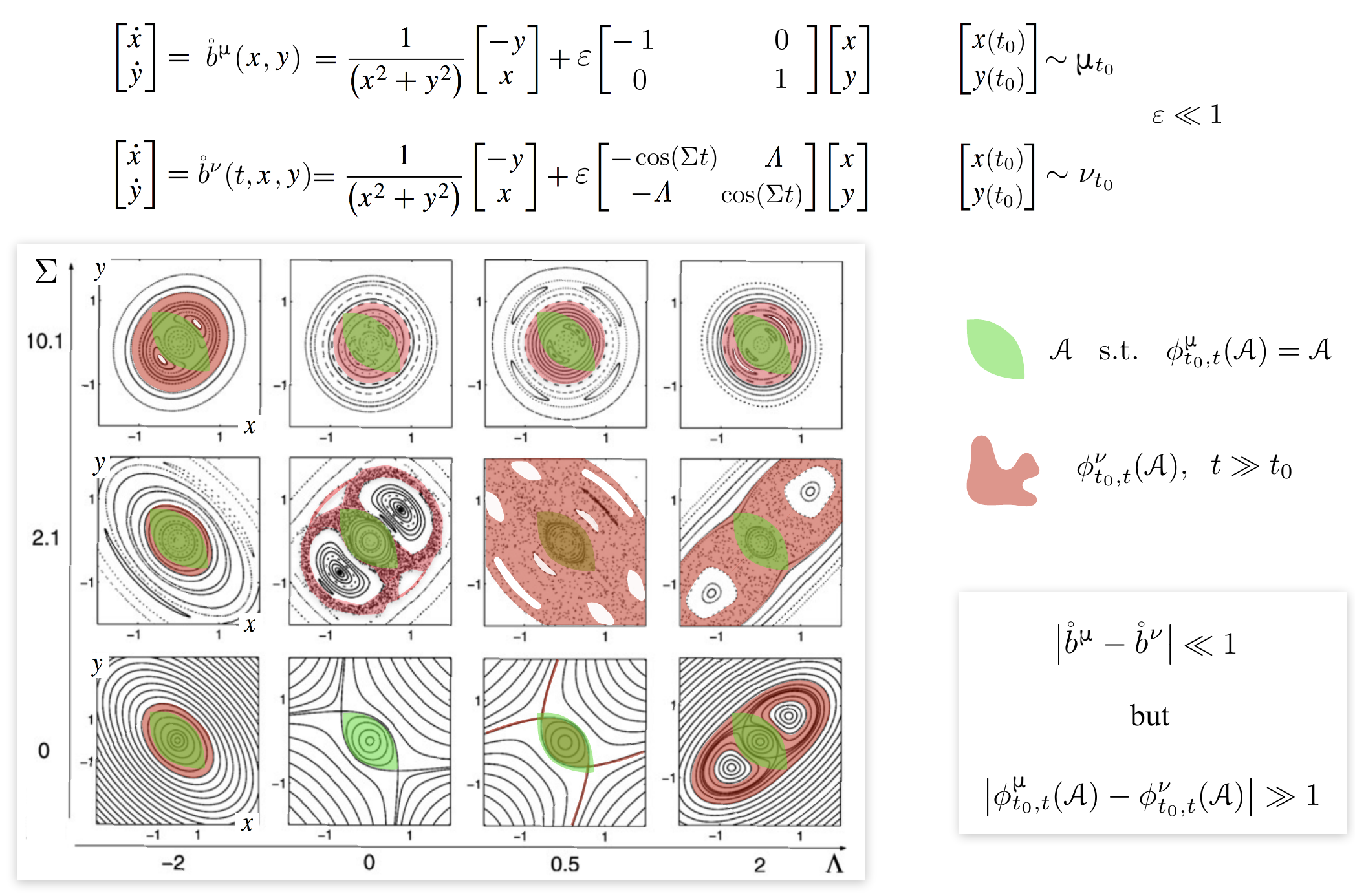}
\vspace{-0.45cm}\caption{\footnotesize An illustration of sensitivity of Lagrangian predictions to Eulerian biases/uncertainties (Hamiltonian, time-periodic,  deterministic setting on a fixed Poincar\'e section). The discrepancy in the structure of trajectories starting from the same set $\mathscr{A}$ of initial conditions, and  evolved either via  $\phi^\text{\textmugreek}_{t_0,t}(\mathscr{A})$ induced by $\brr^\text{\textmugreek}$ or via $\phi^\nu_{t_0,t}(\mathscr{A})$ induced by $\brr^\nu$ can become large even if $|\brr^\text{\textmugreek}-\brr^\nu|\ll 1$, $\forall\, t\geqslant t_0$ (e.g.,  KAM theorem \& Arnold diffusion \cite{kolm64,arnld64}). Here, $\phi^\text{\textmugreek}_{t_0,t}(\mathscr{A})=\mathscr{A}$, $\forall t\geqslant t_0$, and the approximating vector field is such that $\brr^\nu = \brr^\text{\textmugreek}+\varepsilon \hspace{0.04cm}\tilde b(t,\Lambda,\Sigma)$, \,$\tilde b(t,0,0) = 0$. Different values of (uncertain) parameters $\Lambda,\Sigma$  in the $\varepsilon$-small term of $\brr^\nu$ lead to $\big|\phi^\text{\textmugreek}_{t_0,t}(\mathscr{A})-\phi^\nu_{t_0,t}(\mathscr{A})\big|\gg1$ for $t\gg t_0$; the red patches correspond to $\phi^\nu_{t_0,t}(\mathscr{A})$ while the green patch indicates $\phi^\text{\textmugreek}_{t_0,t}(\mathscr{A})$. See \S\ref{out_main} for an outline of suitable bounds.}\label{motiv_2}\vspace*{-.35cm}
\end{figure}

Despite superficial similarities,  LUQ  is distinctly different from uncertainty quantification in the Eulerian case which is concerned with estimating the lack of information between the fields $(\mathring{b}^\mm,\sigma^\mm)$ and their approximations $(\mathring{b}^\nu,\sigma^\nu)$; e.g.,  \cite{MajdaGershgorin11a,MajdaGershgorin11b,branmaj12,branmaj12c,branmaj12d,branmaj14,brani15}.  Importantly, Eulerian accuracy does not, in general, imply Lagrangian accuracy; simple examples are sketched in Figures~\ref{motiv_1} and ~\ref{motiv_2}. Even if
$\big|\mathring{b}^\mm-\mathring{b}^\nu\big|\,{\ll}\, 1$\,\footnote{\,Throughout, $|\ccdot|$ denotes the $L^2$ norm and $\|\ccdot\|_\textsc{hs}$ denotes the Hilbert-Schmidt norm; see the Glossary.} , $\|\sigma^\mm-\sigma^\nu\|_\textsc{hs}\,{\ll}\, 1$, for all $x\in \Xt=\XXt$, $t\in \Ic$, this does not  imply
 that  $\int_\mathscr{A}\big| \E[\ff(\phi^\mm_{t_0,t}(x))]- \E[\ff(\phi^\nu_{t_0,t}(x))]\big|dx\,{\ll}\,1$, $\mathscr{A}\subset\Xt$, $t\in \Ic$, since the path space laws  of $\phi_{t_0,\scaleobj{1.25}{\ccdot}}^\mm$ and $\phi_{t_0,\scaleobj{1.25}{\ccdot}}^\nu$ can be very different\footnote{\,In this work we quantify discrepancies between probability measures in terms of premetrics given by certain divergences (see \S\ref{phi_def}),  since they are the most suitable for our subsequent work on information-geometric inference and estimation in Lagrangian prediction problems. For $\big|\mathring{b}^\mm-\mathring{b}^\nu\big|, \|\sigma^\mm-\sigma^\nu\|_\textsc{hs}\sim\mathcal{O}(\varepsilon)$, one can obtain a Gronwall-type bound on $|\phi_{t_0,t}^\mm(x,\ccdot)-\phi_{t_0,t}^\nu(x,\ccdot)|$, $t\in [t_0, \, t_0+T)$, in terms of $\varepsilon \hspace{.04cm}T$. However, such a bound becomes uninformative for any fixed $\varepsilon$ (even if $\varepsilon \ll 1$) and  a sufficiently large $T$; it also does not imply closeness of the laws of $\phi_{t_0,\scaleobj{1.25}{\ccdot}}^\mm$, $\phi_{t_0,\scaleobj{1.25}{\ccdot}}^\nu$, in terms of divergences or  $L^1$ type norms. More versatile and  tighter bounds are derived in~\S\ref{out_main}.}; this fact is well known in dynamical systems,  bifurcation theory, nearly-integrable chaotic dynamics (e.g., KAM theorem (\cite{kolm64,wiggins92}),  Arnold Diffusion~\cite{arnld64}, and  noise induced phenomena (e.g., \cite{berg06,Doan_2018}).  To compound matters, one is often interested in minimising the error in Lagrangian predictions based on a parametric family of approximate (Eulerian) models where $\inf_{\varepsilon}\big|\mathring{b}^\mm-\mathring{b}^\nu_\varepsilon\big|$, $\inf_\varepsilon\|\sigma^\mm-\sigma^\nu_\varepsilon\|_\textsc{hs}$ are not small, $\sup_{\varepsilon}\big|\mathring{b}^\mm-\mathring{b}^\nu_\varepsilon\big|$, $\sup_\varepsilon\|\sigma^\mm-\sigma^\nu_\varepsilon\|_\textsc{hs}$ might not exist, yet an optimal model for Lagrangian predictions is needed. 

Quantification and mitigation of error and uncertainty in Lagrangian predictions based on uncertain Eulerian fields 
requires the following three major steps:
\begin{itemize}[leftmargin=0.8cm] 

\vspace{-0.0cm}\item[\bf (i)] Determination of an appropriate probabilistic measure of discrepancy between two Lagrangian (path-based) predictions, which is sensitive to discrepancies in features relevant in applications,  is  computationally tractable, and is such that it can be  utilised in information-geometric  analysis of statistical estimation and inference on families of models.

\item[\bf (ii)] Identification of the most important Lagrangian structures, i.e., subsets of paths generated by $t\mapsto\phi^\mm_{t_0,t}(\xxt,\om)$ or their local averages around $\xxt\in \XXt$, $t\mapsto\int_\XXt\E\big[\phi^\mm_{t_0,t}(\xi)\big]\hspace{.04cm}\text{\textmugreek}^{\xxt}_{t_0}(d\xi)$, and a framework which allows for a systematic `tuning' of such structures in order to minimise the loss of information in Lagrangian predictions from simplified Eulerian models.

\item[\bf (iii)] Derivation of bounds on the error in  Lagrangian observables estimated from approximate models of the original dynamics, and bounding errors in the underlying probability measures. 
\end{itemize}

\noindent  Here, we combine (i)--(iii) in a systematic fashion in order to  develop a unified and  systematic  framework for the analysis of uncertainty in Lagrangian predictions, and we derive a hierarchy of bounds which allow for sensitivity analysis and mitigation of error in such~predictions. 

\begin{samepage}
The contents of this article are as follows. Section \ref{setup_sec} outlines the general setup, the hierarchy of main results and notation. In section \ref{Phi_D}, we introduce $\varphi$-divergences and  discuss  their main properties, following \cite{MBUda18}. Then, we utilise these $\varphi$-divergences to derive general information inequalities which  provide a unified framework for uncertainty quantification when estimating observables; in particular, these inequalities apply to path-based observables which we focus~on.  In section~\ref{SFlow}, we recall some relevant definitions and results concerned with  stochastic flows, which are then used in section~\ref{Main_Bounds} to characterise specific information bounds that are crucial in path-based (Lagrangian) considerations. Abstract bounds on the uncertainty in Lagrangian predictions in terms of the Eulerian error are obtained  in section~\ref{Inf_IP}; these bounds are particularly useful in theoretical considerations. Section~\ref{ftdr_sec} discusses a class of bounds which are amenable to numerical considerations, and  which involve scalar fields of path-based {\it divergence rates} in stochastic flows derived in \cite{MBUda18}. Extension of the results obtained for  time-marginal probability measures sections \ref{Inf_IP} and \ref{ftdr_sec}  to measures on  path space are discussed in section~\ref{Path-space}. Section~\ref{tests}  deals with  an application of our results to a toy example of a slow-fast SDE. We close with some remarks on future work in section~\ref{conclusions}. A number of technical proofs are  discussed in the Appendix. A glossary of  frequently used notation is included at the end of the paper for readers'~reference.
\end{samepage}

\section{Main results} \label{setup_sec}

 \vspace{-0.2cm} Here, we outline  the framework for uncertainty quantification in Lagrangian predictions, and we highlight the main results. Frequently used notation is summarised in the Glossary at the end of the paper. A simple example illustrating some of the developed tools is discussed in~\S\ref{tests}.

\addtocontents{toc}{\protect\setcounter{tocdepth}{1}}

\vspace*{-0.1cm}\subsection{Problem setup and notation}\label{setp}

We  denote the Wiener probability space associated with an $m$-dimensional Wiener process by  $\PS$; i.e., 
$\Om\simeq\mathcal{C}_{0}(\R;\R^m)$ is identified with a subspace of continuous functions $\mathcal{C}(\R;\R^m)$ which vanish at zero (e.g., \cite{Arnold1}). $\F$ is the Borel  $\mathfrak{S}$-algebra on~$\Om$ generated by the Wiener process, and $\p$ is the Wiener measure on $\mathcal{F}$. 

The original dynamics is generated by a continuous map, $\phi^\mm_{t_0,t}\,{:}\;\XXt\,{\times}\, \Omega_\mm \rightarrow \XXt$, 
which satisfies 
$$\phi^\mm_{t_0,t_0}(\ccdot,\om)=\textrm{id}_\XXt, \quad \phi^\mm_{t_0,t}(\ccdot,\om) = \phi^\mm_{u,t}(\phi^\mm_{t_0,u}(\ccdot,\om),\om),  \quad \forall\, t,u\in \Ic:=[t_0,\,t_0+T)\subset \R,\;\; \textrm{a.a.} \;\,\om, $$
where $\om\in \Om_\mm$ is a sample space associated with the probability space $(\Om_\mm,\mathcal{H}_\mm,\Ms^\mm_{t_0})$; $\Om_\mm$ and $\mathcal{H}_\mm$ can be identified with, respectively, $\Om$ and $\mathcal{F}$ in a standard way (e.g. \cite{Arnold1,oksendal,Stroock79}).

 Throughout,   $\mathfrak{M} \,{=}\,\Xt\times \mathcal{Y}$, $ \dim\XXt = \ell$,  where $\Xt = \R^d$ or a flat torus $\Xt\,{=}\,\bar{\T}^d$; similarly, $\mathcal{Y} = \R^{\ell-d}$ or $\mathcal{Y} = \bar{\T}^{\ell-d}$. The sets  of  all probability measures on $\XXt$ and $\Xt$ are denoted, respectively, by $\PP(\XXt)$ and  $\PP(\Xt)$.  The approximate dynamics is generated by  $\phi^\nu_{t_0,t}\,{:}\;\Xt\,{\times}\, \Omega_\nu \rightarrow \Xt$, $\Om_\nu\simeq \Om$, for $t\in \Ic$, with the same properties on~$\Xt\subseteq \XXt$ as those of $\phi^\mm_{t_0,t}$ on~$\XXt$. Under some weak assumptions discussed  in~\S\ref{SFlow} the maps $\phi^\mm_{t_0,t}$, $\phi^\nu_{t_0,t}$ represent  {\it stochastic flows} which we will rely on throughout. 

  Path spaces associated with $t\mapsto \phi^\mm_{t_0,t}$ and $t\mapsto \phi^\nu_{t_0,t}$ are given by $\mathcal{W}^\mm_\ell = \mathcal{C}(\Ic,\XXt)\simeq \Om_\mm$ and $\mathcal{W}^\nu_d = \mathcal{C}(\Ic,\Xt)\simeq \Om_\nu$. Probability measures induced by the laws of $\phi^\mm_{t_0,\scaleobj{1.25}{\ccdot}}(\xxt,\ccdot)$ and $\phi^\nu_{t_0,\scaleobj{1.25}{\ccdot}}(x,\ccdot)$ on $\mathcal{W}^\mm_\ell$ and $\mathcal{W}^\nu_d$ are denoted by $\Ms^\mm_{t_0,x}\in \PP(\mathcal{W}_\ell)$, $\Ns^\nu_{t_0,x}\in \PP(\mathcal{W}_d)$, and   we set 
 \begin{equation}\label{avg_mu_om}
  \E\big[g\big(\phi^{\mm}_{t_0,\scaleobj{1.25}{\ccdot}}(\xxt)\big)\big] :=\!\scaleobj{.92}{\int_{\Om_\mm}} g\big(\phi_{t_0,\cccdot}^\mm(\xxt,\om)\big)\Ps_{t_0,\xxt}^\mm(d\om), \qquad g\in \mathbb{M}(\mathcal{W}_\ell),
  \end{equation}
\begin{equation}\label{avg_nu_om}
 \E\big[f\big(\phi^{\nu}_{t_0,\scaleobj{1.25}{\ccdot}}(x)\big)\big] :=\!\scaleobj{.92}{\int_{\Om_\nu}} f\big(\phi_{t_0,\cccdot}^\nu(x,\om)\big)\Ns_{t_0,x}^\nu(d\om), \qquad f\in \mathbb{M}(\mathcal{W}_d). 
 \end{equation}
 Given the (random) paths  $t\mapsto \phi^\mm_{t_0,t}(\xxt,\om)$,  and $t\mapsto \phi^\nu_{t_0,t}(x,\om)$,  labelled by their  initial conditions $\xxt\in \XXt$, $x\in \Xt$, we will focus on locally averaged path-based observables  
\begin{equation}\label{avg_obss}
x\mapsto\int_\Xt\E\big[\ff\big(\phi^\nu_{t_0,\scaleobj{1.25}{\ccdot}}(\xi)\big)\big]\nu^x_{t_0}(d\xi), \qquad \xxt\mapsto\int_\XXt\E\big[\ff\big(\pi^\nu_\mm\circ\phi^\mm_{t_0,\scaleobj{1.25}{\ccdot}}(\zeta)\big)\big]\hspace{.04cm}\text{\textmugreek}_{t_0}^{\xxt}(d\zeta),
\end{equation}
where $\mm_{t_0}^{\xxt}\in \PP(\XXt)$ and $\nu^x_{t_0}\in \PP(\Xt)$   are concentrated on some initial conditions $\xxt\in \XXt$, $x\in \Xt$,  and $\pi^\nu_\mm\,{:} \;\XXt\rightarrow \Xt$  is a projection onto $\Xt\subseteq \XXt$; i.e., for  $\xxt=(x,y)\in \XXt$,   $\pi^\nu_\mm(x,y) = x\in \Xt$.

 For the most part, either for convenience or by necessity (see~\S\ref{Inf_IP}), the original and approximate dynamics will be assumed to be generated  by the SDE's (\ref{SDE1}) and (\ref{SDE2}) with standard conditions on the coefficients $(\brr^\mm, \sigma^\mm)$, $(\brr^\nu, \sigma^\nu)$
  for existence and uniqueness of global solutions on $\Ic$ which  can be represented by stochastic flows (e.g., \cite[Thm.~4.5.1]{Kunitabook} or \S\ref{SFlow}); namely 
 $$X^\mm_t(\om) = \phi^\mm_{t_0,t}(\xxt,\om),\; \;\xxt\in \XXt, \qquad\;\; X^\nu_t(\om) = \phi^\nu_{t_0,t}(x,\om),\; \;x\in \Xt, \qquad\quad  \p\,\textrm{-\,a.s.},$$
where the uncertainty in the initial conditions, or a distribution of initial conditions of interest, is prescribed by $\mmo\in \PP(\XXt)$ and $\nu_{t_0}\in \PP(\Xt)$. Many results derived in what follows apply to a broader class of flows than those induced by strong solutions of SDEs but the Markovian setting serves as a useful reference setup. However, in the framework we develop below, it will be important that the stochastic (semi-martingale) flows are diffeomorphisms so that the associated probability measures are sufficiently regular; in the case of SDE dynamics, sufficient conditions for the solutions to generate flows of diffeomorphisms are outlined later in~\S\ref{SFlow}.

\vspace{.1cm}
Let $(\mmt)_{t\in \Ic}$, $\mmt\,{\in}\, \PP(\XXt)$, and $(\nu_t)_{t\in \Ic}$, $\nu_t\,{\in}\, \PP(\Xt)$, denote a measure-valued processes induced by the underlying dynamics on, respectively, $\XXt$ and $\Xt$. If the dynamics is induced by the SDEs (\ref{SDE1}) and  (\ref{SDE2}), the time-marginals of the laws of $\phi^{\mm}_{t_0,\cccdot}$ and  $\phi^{\nu}_{t_0,\cccdot}$  satisfy (in the weak sense) the {\it forward Kolmogorov equations} (see, e.g., \cite{Stroock79,Roc-Kry} and Definition \ref{weak_kol}) for $t\in \Ic$
\begin{align}\label{F_Kol}
(a)\quad  \partial_t\mmt = \LG^{\mm*}_t\hspace*{0.01cm}\mmt, \quad \mmo\in \PP(\XXt),
\hspace{1.3cm} (b) \quad \partial_t\nu_t = \LG^{\nu*}_t\hspace*{0.01cm}\nu_t, \quad \nu_{t_0}\in \PP(\Xt),
\end{align}
where $\LG^{{ \mu }*}_t$ is the $L^2(\XXt,\mmt)$ dual of the  second-order differential operator $\LG^{ \mm}_t$  given by 
\begin{equation}\label{gen}
 \LG_t^\mm f(\xxt) = \sum_{i=1}^\ell b_i^\mm(t,\xxt)\partial_{\xxt_i}f(\xxt)+\frat\sum_{i,j=1}^\ell a_{ij}^\mm(t,\xxt)\partial^2_{\xxt_i\xxt_j}f(\xxt), \qquad f\in \mathcal{C}^{2}(\XXt),
\end{equation}

\vspace{-.1cm}
\noindent with $ b_i^\mm(t,\xxt) \,{:=}\, \brr_i^\mm(t,\xxt)\,{+}\,c_i^\mu(t,\xxt)$, $c_i^\mm(t,\xxt) := \frac{1}{2}\sum_{k,j=1}^{m,d} \sigma_{jk}^\mm(t,\xxt)\partial_{\xxt_j} \sigma_{ik}^\mm(t,\xxt)$ the Stratonovich correction, and $a_{ij}^\mm :=  \sum_{k=1}^m \sigma_{ik}^\mm\sigma_{jk}^\mm$; analogous notation applies to $\LG_t^\nu$ associated~with~(\ref{SDE2}) on $\Xt$. For deterministic dynamics, i.e., $\sigma^\mm=\sigma^\nu=0$, the evolution of the time-marginal measures $(\mmt)_{t\in \Ic}$, $(\nu_t)_{t\in \Ic}$  is given simply by the push forward of the probability measures $\mmo$, $\nu_{t_0}$ on the initial conditions  under, respectively, $\phi^\mm_{t_0,t}(\xxt,\om) = \psi_{t_0,t}^\mm(\xxt)$ and $\phi_{t_0,t}^\nu(x,\om) = \psi_{t_0,t}^\nu(x)$ for $t\in \Ic$.
  
 If $\mm_t$,  $\nu_t$ solve (\ref{F_Kol}),  there exist  path space probability measures, $\Ms_{t_0}^\mm\in \PP(\mathcal{W}_\ell)$, $\Ns_{t_0}^\nu\in \PP(\mathcal{W}_d)$, such that (see~\S\ref{SFlow})
 \begin{align*}
\E^{\mm_t}[g] &=  \!\scaleobj{.92}{\int_{\XXt}}g(\xxt)\mm_t(d\xxt) =\!\scaleobj{.92}{\int_{\XXt}}\scaleobj{.92}{\int_{\Om_\mm}} g\big(\Solm^\mm(\xxt,\om)\big)\Ps_{t_0,\xxt}^\mm(d\om)\mmo(d\xxt)=: \E[g\big(\Solm^\mm\big)], \quad\, \,g\in \mathcal{C}^{2}(\XXt),  \\[.13cm]
\E^{\nu_t}[f]  &=  \!\scaleobj{.92}{\int_{\Xt}}f(x)\nu_t(dx) = \!\scaleobj{.92}{\int_{\Xt}}\scaleobj{.92}{\int_{\Om_\nu}} f\big(\Solm^\nu(x,\om)\big)\Ns_{t_0,x}^\nu(d\om)\nu_{t_0}(dx)=:\E[f\big(\Solm^\nu\big)], \quad f\in \mathcal{C}^2(\Xt).
 \end{align*}   

\vspace{.15cm}
\noindent In particular, for $g(\xxt) = (f\,{\circ}\,\pi^\nu_\mm)(\xxt)$  so that $g(\xxt) = f(\pi^\nu_\mm(x,y)) = f(x)$, we have  
\begin{align*}
\E^{\mu_t}[f] =  \!\!\scaleobj{.92}{\int_{\XXt}}\!( f\hspace{.04cm}{\circ}\,\pi^\nu_\mm)(\xxt)\mm_t(d\xxt) \!=\! \!\scaleobj{.92}{\int_{\XXt}}\scaleobj{.92}{\int_{\Om_\mm}} \!f\big(\pi^\nu_\mm\hspace{.02cm}{\circ}\,\Solm^\mm(\xxt,\om)\big)\Ps_{t_0,\xxt}^\mm(d\om)\mmo(d\xxt) = \E\big[f\big(\pi^\nu_\mm\hspace{.02cm}{\circ}\,\Solm^\mm\big)\big]. 
\end{align*}
Note that one can consider locally averaged Lagrangian observables  $\E^{\mu^x_t}[f]$, $\E^{\nu^x_t}[f]$,  by choosing  $\text{\textmugreek}_{t_0}^{\xxt}$ and $\nu^x_{t_0}$ to be concentrated on some initial conditions $\xxt\in \XXt$, $x\in \Xt$ (see (\ref{avg_obss}) and \S\ref{ftdr_sec}).

\vspace{.1cm}
 In order to carry out meaningful analysis\footnote{\label{fokt2} \,If $\Xt\ne\XXt$, probability measures $\mm\in \PP(\XXt)$ and $\nu\in \PP(\Xt)$ are generally singular. Thus, when comparing  the dynamics on $\XXt$ with its approximation on $\Xt$, 
the loss of information in the approximation can be considered as infinite.  We adopt a more `constructive approach' and 
compare probabilistic aspects of the original dynamics on $\XXt$ with its approximation on $\Xt\subset\XXt$. This approach is in line with assessing reduced-order models in applications.} when  $\Xt\subset \XXt$, one has to consider  projections of the time-marginal probability measures $\mmt\in \PP(\XXt)$  associated with the original dynamics onto appropriate probability measures  $\mu_t\in \PP(\Xt)$ which are obtained via marginalisation of their Lebesgue densities\footnote{\,Assumptions enforced in the sequel on the underlying dynamics guarantee the existence of densities w.r.t.~the Lebesgue measures on $\XXt$ and $\Xt$ by the properties of the forward Kolmogorov equation; see \S\ref{SFlow}.}; i.e., for $\mathfrak{M} = \Xt\times \mathcal{Y}$, $\xxt=(x,y)\in \XXt$,  $x\in \Xt$, $y\in \mathcal{Y}$, we set  
 \begin{equation}\label{margs0}
 \mmt(dxdy) = \rrt(x,y)dxdy, \qquad  \mu_t(dx) = \rho_t^\mu(x)dx, \quad \rho^\mu_t(x) = \scaleobj{.9}{\int_\mathcal{Y}}\,\rrt(x,y)dy.
 \end{equation}

The discrepancy in the information content between the original and approximate dynamics on $\Xt$  is considered  via $\varphi$-divergences, which for time-marginal measures are defined by 
\begin{align}\label{phidiv}
\Df(\mu_t\|\nu_t) = \scaleobj{.9}{\int_\Xt}\,\varphi\big(d\mu_t/d\nu_t\big)d\nu_t, \qquad \mu_t,\nu_t\in \PP(\Xt),
\end{align}
where $d\mu_t/d\nu_t$ is the Radon-Nikodym derivative, and $\varphi\in \mathcal{C}^2(\R^+)$ is a strictly convex scalar function (see \S\ref{phi_def}); loss of information $\Df\big(\Ms_{t_0}^\mm\|\Ns_{t_0}^\nu\big)$ in probability measures  $\Ms_{t_0}^\mm\in \PP(\mathcal{W}_\ell)$, $\Ns_{t_0}^\nu\in \PP(\mathcal{W}_d)$ on path spaces is considered within this framework in \S\ref{Path-space}.
 
\vspace{.1cm}
The family of $\varphi$-divergences contains some well-known premetrics and metrics (e.g., Kullback-Leibler divergence, 
the Hellinger distance, or the Total Variation distance; see \S\ref{phi_def}).  
Apart from their utility for deriving a hierarchy of bounds outlined below,  the choice of $\varphi$-divergences over other metrics is driven by their unique properties suitable for information-geometric considerations in the context of  statistical estimation and inference, which is a focus of a forthcoming~work.

\subsection{Outline of main results}\label{out_main} In the following sections we develop a framework for the analysis of uncertainty in Lagrangian predictions of path-based observables, and we derive a hierarchy of bounds which allows for the mitigation of error in this setup. For brevity, we outline the results for time-marginal probability measures; see  \S\ref{Path-space} for extensions to path space measures. We obtain
\begin{itemize}[leftmargin=.7cm]
\item[\bf (a)] $\varphi$-information inequalities for observables (see \S\ref{Info_ineq} for details)\\[.2cm]
\hspace*{.0cm}$\hat{\mathcal{K}}_{\varphi,\sff}^{\nu}\big(\!-\Df(\mu_t\|\nu_t)\big)\leqslant  \E^{\mu_t}[\ff] - \E^{\nu_t}[\ff] \leqslant \mathcal{K}_{\varphi,\sff}^{\nu}\big(\Df(\mu_t\|\nu_t)\big), \hspace{.6cm}  \,t\in \Ic,\;\;\mu_t,\nu_t\in \PP(\Xt)$,  
\\[.2cm]
where $\Ic=[t_0,\,t_0+T)$, $\E^{\mu}[\ff]:= \int \!\! \ff d \mu$, and   $\mathcal{K}_{\varphi,\sff}^{\nu}(s)\rightarrow 0$, $\hat{\mathcal{K}}_{\varphi,\sff}^{\nu}(-s)\rightarrow 0$ as $s\,{\downarrow}\, 0$. These bounds   are tight and general; i.e., they are not restricted to Markovian processes or solutions of SDE's/ODE's, although in the latter setting we have for path-based observables 
$$\hspace{.8cm}\E^{\mu_t}[f] = \E\big[f\big(\pi^\nu_\mm\hspace{0.04cm}{\circ}\,\Solm^\mm\big)\big],\qquad\quad   \E^{\nu_t}[f] = \E\big[f\big(\Solm^\nu\big)\big], $$
as highlighted in \S\ref{setp} and~\S\ref{SFlow}. For autonomous dynamics the above inequality  can be interpreted as a bound on the difference between solutions of two backward Kolmogorov equations in terms of a $\varphi$-divergence between the solutions of the corresponding forward Kolmogorov equations, but we  emphasise the dynamical systems' interpretation. 
\end{itemize}

\noindent Moreover, we derive two distinct bounds on $\Df(\mu_t\|\nu_t)$  which can be combined with (a), namely: 

\vspace{-0.05cm}\begin{itemize}[leftmargin = .8cm]
\item[\bf (b)] Bounds on $\Df(\mu_t\|\nu_t)$ in terms a functional on certain ``reconstructed'' fields  involving the coefficients $(\brr^\mm, \sigma^\mm)$ of the original SDE dynamics (\ref{SDE1}), and the coefficients $(\brr^\nu, \sigma^\nu)$ in the approximation (\ref{SDE2}). The overall  structure of these bounds can be expressed as (see~\S\ref{Inf_IP})
\begin{equation*}
\hspace{1.5cm}\Df(\mu_t\|\nu_t)\leqslant \Upsilon_{\varphi,t}^{\mm,\nu}\big(\,\brr^\mm,\brr^\nu, \sigma^\mm,\sigma^\nu\,\big), \hspace{1.2cm}  \,t\in \Ic, \;\mu_t,\nu_t\in \PP(\Xt),\;\mu_{t_0} = \nu_{t_0}.
\end{equation*}
These bounds hold for dynamics induced by SDE's with sufficiently non-degenerate diffusion coefficients (we assume uniform ellipticity to simplify matters), and they cannot yet be directly extended to the deterministic case, unless viscosity-type limits are employed. 
\item[\bf (c)] Bounds on $\Df( \mu_t\|\nu_t)$ in terms of probabilistic measures of expansion rates 
$$\hspace{1.3cm}\Df(\mu_t\|\nu_t)\leqslant \big\vert \Df( \mu_t\|\mu_{t_0})-\Df(\nu_t\|\nu_{t_0})\big\vert, \hspace{.7cm}   \,t\in \Ic, \;\mu_t,\nu_t\in \PP(\Xt), \;\mu_{t_0} = \nu_{t_0}.\hspace{.4cm}$$
Similar to (a), these bounds are not restricted to time-marginal measures induced by solutions of SDEs or ODEs  (see \S\ref{ftdr_sec}, and see \S\ref{Path-space} for a generalisation to path space measures).
 \end{itemize} 
 The above bounds are non-uniform in $T$, unless the underlying dynamics have stationary or cyclo-stationary measures. In general, none of the path-based criteria in (a)-(c) are trivially linked to the  minimisation of the Eulerian discrepancies $|\brr^\mm-\brr^\nu|$ and $\|\sigma^\mm-\sigma^\nu\|_\textsc{hs}$. A combination of bounds (a) and (b)  provides an analytically tractable connection between the Eulerian error  and the uncertainty in Lagrangian predictions from the  estimates of $\E\big[\ff\big(\phi^\nu_{t_0,t}(x)\big)\big]$.
 The bound in~(c) can be cast in terms of  the difference between {\it Finite-Time $\varphi$-Divergence Rate} ($\varphi$-FTDR) fields~\cite{MBUda18} which utilise a recently developed probabilistic framework for quantifying local expansion rates in stochastic flows;  roughly, the scalar $\varphi$-FTDR fields are defined via
 $$\hspace{3cm} x\mapsto \Df(\mu^x_t\|\mu^x_{t_0}), \qquad x\mapsto \Df(\nu^x_t\|\nu^x_{t_0}),\hspace{1.cm}   \,t\in \Ic,\;\mu_t,\nu_t\in \PP(\Xt),$$
 where the time-marginal probability measures, $\mu^x_t$, $\nu^x_t$, evolve from  $\mu^x_{t_0} = \nu^x_{t_0}$ concentrated on a neighbourhood of $x\in \Xt$. Importantly, the bound (c) can be utilised within a computational framework, and the combination of (a) and (c) implies that error in Lagrangian predictions can be mitigated by tuning  $\varphi$-FTDR fields in the approximate dynamics to optimally reproduce the original expansion rate~fields.   The fact that different conditions are required for accuracy of Eulerian and Lagrangian predictions should not be surprising to experts working on transport in dynamical systems, and geophysical/oceanographic Lagrangian predictions.  However, to the best of our knowledge, this is the first time where systematic and  rigorous bounds have been derived for improving Lagrangian estimates. An introductory  illustration of the developed tools is discussed~in~\S\ref{tests}.

\vspace{0.cm}\section{Information measures and information inequalities}\label{Phi_D}
Here, we first introduce a class of generalised information divergences (\S\ref{phi_def}) defined via a family of specifically normalised strictly convex functions.  Then, in \S\ref{Info_ineq}  we derive a class of new versatile information inequalities which provide bounds on the error in estimation of observables in terms of $\varphi$-divergences; these inequalities can be related to  a number of well-known bounds but they are tighter and apply to a larger class of observables. The resulting information inequalities are subsequently combined with additional bounds on the $\varphi$-divergences themselves (see~\S\ref{Main_Bounds}) to form a general probabilistic/information-theoretic framework for Lagrangian uncertainty quantification (see~\S\ref{out_main} for an outline and \S\ref{tests} for a simple illustration of forthcoming applications). 

\addtocontents{toc}{\protect\setcounter{tocdepth}{2}}
\subsection{$\varphi$-divergences}\label{phi_def}
Consider a class of strictly convex functions $\varphi: [0,\infty) \rightarrow (-\infty, \infty)$ satisfying the following  {\it normality conditions}
\begin{align}\label{Normality}
\varphi(1) = 0, \quad \nabla\varphi(1) = 0, \quad \inf_{x>0}\varphi(x)>-\infty.
\end{align}
Let $\mu$ and $\nu$ be probability measures on a Polish space $(\mathcal{\XX}, \Bb(\XX))$,
i.e.,  $\mu,\nu\in\PP(\XX)$ are not necessarily induced  by  the SDEs (\ref{SDE1}) and  (\ref{SDE2}), and they can represent  probability measures on the path space or their time-marginals. Then, the $\varphi$-divergence between $\mu$ and $\nu$ is defined by\,\footnote{\,Definition of $\Df$ in (\ref{phi_dA}) is related to that of $f$-divergence due to Csisz\'ar \cite{Csiszar72,csiszar91,Ciszar08}.  However, the conditions~(\ref{Normality}) are often not imposed and $f$-divergences might not even be premetrics. Here, the constraint imposed on $\varphi$ which generates $\Df$ removes the symmetries $\mathcal{D}_{f+c(u-1)} \,{=}\, \mathcal{D}_f$, $\mathcal{D}_{cf} \,{=}\, \mathcal{D}_f$, $c\ne0$, present in $f$-divergences. }
\begin{align}\label{phi_dA}
\Df(\mu\|\nu) = \begin{cases} \int_{\Xt}\varphi\left( \frac{d\mu}{d\gamma}\big/\frac{d\nu}{d\gamma}\right)d\nu, &\;\;\ \textrm{if }\; \mu,\nu\ll\gamma,\;\varphi\left( \frac{d\mu}{d\gamma}\big/\frac{d\nu}{d\gamma}\right)\in L^1(\XX, \nu),\\
+\infty, &\;\;\; \textrm{otherwise},
\end{cases}
\end{align}
where $\gamma$ is any finite, positive reference measure on $(\XX, \Bb(\XX))$. Note that the definition in (\ref{phi_dA}) is independent of the reference measure due to the uniqueness of the Radon-Nikodym derivative;  this property implies invariance of the above definition w.r.t.\ a diffeomorphic change of variables.  

  In general, $\Df$ is not symmetric and it does not satisfy the triangle inequality (see \cite{MBUda18} and also, e.g., \cite{Csiszar72,csiszar91,Ciszar08}). However,   $\Df$ is {\it information monotone} in the sense that 
  $$\Df(\mu\|\nu)\geqslant \Df(\mu_\mathcal{A}\|\nu_\mathcal{A})$$
for $\mu_{\mathcal{A}}(B)= \mu(A_n\cap B)$, $\nu_{\mathcal{A}}(B) =\nu(A_n\cap B)$ for all $B\in \Bb(\XX)$ and for any measurable partition $\mathcal{A} = \{A_n: n\in \N\}$ of $\XX$. Information monotonicity also implies that $\Df$ is a {\it premetric}; i.e.,  $\Df(\mu\|\nu)\geqslant 0$ and $\Df(\mu\|\nu) =0$ iff $\mu = \nu$ almost everywhere. Importantly, $\varphi$-divergences belong to a class of convex integrals which admit the following duality representation (e.g.,~\cite{Gmeas}) 
\begin{align}\label{Dua_l}
\Df(\mu\|\nu) &= \sup_{f\in \mathcal{C}_\infty(\XX)}\Big\{\int_{\XX} f(x)\mu(dx) - \int_{\XX}\varphi^*(f(x))\nu(dx)\Big\}
\\[.2cm]
&\hspace{0cm} =\sup_{f\in \mathcal{C}_\infty(\XX)}\Big\{\langle f, \mu\rangle - \langle \varphi^*(f), \nu\rangle\Big\}, \qquad\qquad  \nu, \mu \in \PP(\XX), \notag
\end{align}
where $\varphi^*$ is the Legendre-Fenchel convex conjugate of $\varphi$; i.e., 
\begin{align}\label{Phistar}
\varphi^*(\xi) = \sup_{u\geqslant 0}\big\{u \,\xi - \varphi(u)\big\} \quad \forall\;\xi\in \mathbb{R}.
\end{align}
It follows immediately from the above representation that $\Df: \PP(\XX)\times\PP(\XX) \rightarrow\Rp$ is lower semicontinuous, and that it is jointly strictly convex in its arguments. Information-monotonicity of  $\varphi$-divergences and the above properties allow to uniquely determine a special Riemannian geometry on the manifold of probability measures in which a Pythagorean-like decomposition and (non-metric) geodesic projections are crucial in applications of {\it information-geometric} framework to statistical estimation (e.g., \cite{chentsov72,amari00,amari09,Amari10, Amari16}). Moreover, information monotonicity is naturally imposed by physical constraints when simplifying/coarse-graining the original dynamics.  The suitability of  a given $\varphi$-divergence for uncertainty quantification depends on the application, and on the considered submanifold of probability measures (e.g., \cite{amari00,amari09,Amari16,Csiszar72,csiszar91,Ciszar08}).  Given our future aim of exploiting the information-geometric framework for  uncertainty quantification in families of reduced-order models, we consider the whole family of $\varphi$-divergences. 

\vspace{.1cm}
Various well-known divergences used in information theory, probability theory and statistics are derived from (\ref{phi_dA}) with an appropriate choice the convex function~$\varphi$, namely:

 \vspace*{-0.1cm}
 \renewcommand{\arraystretch}{1.2}
 \begin{table}[H]
\begin{tabular}{ | c | c|}
  \cline{1-2}
  $\Df$ - notation &  \hspace{3.5cm}$\varphi(u)$ \qquad $(u=d\mu/d\nu)$\ \\
  \hline
  KL-divergence - $\D_\textsc{kl}(\mu\|\nu)$ & $u\log u-u+1$\\
  \hline
  Hellinger distance - $\D_\textsc{h}(\mu\|\nu)$ or $d_H(\mu,\nu)$& $(\sqrt{u} -1)^2$ \\
  \hline
  Total variation - $\D_{\textsc{tv}}(\mu\|\nu)$ or $\|\mu-\nu\|_{\textsc{tv}}$& $\frac{1}{2}\vert u-1\vert$\\
  \hline
 $\chi^\alpha$-divergence - $\D_{\chi^\alpha}(\mu\|\nu)$ or $\chi^\alpha(\mu\|\nu)$& \hspace{3.1cm} $|u-1|^\alpha, \qquad 1\leqslant \alpha<\infty$ \\ 
 \hline
 $\alpha$-divergence - $\D_\alpha(\mu\|\nu)$& \hspace{1.1cm}$\begin{cases} \;\;\;\frac{4}{1-\alpha^2}(1-u^{(1+\alpha)/2}),  \hspace{.4cm}\alpha\neq \pm 1,\\ \;\;\;u\log u, \hspace{2.43cm} \alpha = 1,\\ -\log u, \hspace{2.35cm}\;\;\;\alpha=-1.\end{cases}$\\
  \hline
\end{tabular} 
\label{phi_table}
\captionsetup{width=1\linewidth}

\vspace{.1cm}\caption{\footnotesize Examples of well-known $\varphi$-divergences in a form satisfying (\ref{Normality}) and common notations.
}\label{phidlist}
\end{table}

\vspace*{-.5cm}
\noindent Finally, we define an Orlicz subspace\footnote{\,The set of measurable functions $L_{\varphi*}(\XX;\nu)$ is a subspace of a larger Orlicz space $\mathfrak{L}_{\varphi^*}(\XX;\nu)$  

\vspace{-0.5cm}$$\mathfrak{L}_{\varphi^*}(\XX;\nu) := \big\{ f\in\mathbb{M}(\XX): \quad  \exists\;\gamma>0,\quad \textstyle\int_{\XX}\varphi^*(\gamma f)d\nu<\infty\big\}, \quad \nu\in \PP(\XX).$$} associated with the strictly convex $\varphi$ in (\ref{Normality}) and defined~by 
\begin{align}\label{Orlicz}
L_{\varphi*}(\XX;\nu):=\Big\{f\in \mathbb{M}(\XX): \;\forall \gamma>0,\; \int_{\XX}\varphi^*(\gamma f)d\nu+ \int_{\XX}\varphi^*(-\gamma f)d\nu<\infty\Big\},\quad \nu\in \PP(\XX),
\end{align}
with the Orlicz norm $\Vert f\Vert_{\varphi^*} = \inf\big\{ a\,{>}\,0\,\,{:}\; \int_{\XX}\varphi^*\left(f/a\right)d\nu\leqslant 1\big\}$.  The convex conjugate $\varphi^*$  in~(\ref{Phistar}) is locally bounded and the normality conditions (\ref{Normality})  ensure that  $\varphi^*$ is a Young function, i.e.,
\begin{itemize}[leftmargin = .8cm]
\item[(i)]
$\varphi^*$ is lower-semicontinuous, $\varphi^*(0) = 0$  and $\varphi^*$ is not identically zero, and  

\item[(ii)] $(-\alpha, \alpha)\subset \text{dom}\,\varphi^*$ for some $\alpha>0$. 
\end{itemize}
The above properties imply that the Orlicz subspace $L_{\varphi^*}(\XX; \nu)$ is well-defined and nontrivial in the sense that $ L_{\varphi^*}(\XX; \nu) \neq\{0\}$. The Orlicz subspace $L_{\varphi^*}(\XX;\nu)$ will contain the class of observables for  which the information inequality in subsection \S\ref{Info_ineq} will be formulated.

It is worth stressing that a number of other divergences, including Chernoff \cite{cher52}, Renyi \cite{ren61}, Bregman \cite{breg67} divergences, or proper metrics like the Wasserstein distance (e.g., \cite{Gmeas}), have been extensively used  in various contexts  including information theory, statistical inference, optimisation, image processing, and neural networks; e.g.,~\cite{burnh02,niels88,Amari16, Vajda06, Amari10, Csiszar72, Ciszar08}). However, these divergences and metrics are not, in general, information monotone and/or do not have the right decomposition properties,   and are thus not suitable for our purposes.

\subsection{Information inequalities via $\varphi$-divergences}\label{Info_ineq}
Here, we derive a sequence of general bounds on the error in estimates of observables termed {\it information inequalities}, which are expressed in terms of $\varphi$-divergences. These inequalities provide tight error bounds tailored to a given observable by utilising the variational formulation of $\Df$ in (\ref{Dua_l}), and  they provide an extension of analogous bounds  developed for the KL-divergence in \cite{dupuis16, dupuis13, Jing12} to a much larger class of admissible observables. The results discussed below rely solely on convex-analytic tools, and are not confined to the Lagrangian framework, which we focus on in the subsequent sections. 
If the measures considered below are generated by the SDEs/ODEs~in~(\ref{SDE1})-(\ref{SDE2}),   one sets $\XX=\Xt\subseteq \XXt$ in the derivations below with $\Xt = \R^d$ or  $\Xt\,{=}\,\bar{\T}^\ell$ (see \S\ref{setp})), and $\mu_t,\nu_t\in \PP(\Xt)$. For extensions to  path space probability measures see \S\ref{Path-space}. 

\smallskip
The main bound derived in this section has  the form (see Theorem~\ref{iNF_IN})
\begin{align}\label{AmInfo}
\mathfrak{B}_{\varphi,-}\big(\mu\|\nu;\ff\big)\leqslant \E^{\mu}[\ff] - \E^{\nu}[\ff]\leqslant \mathfrak{B}_{\varphi,+}\big(\mu\|\nu;\ff\big),\qquad \mu,\nu\in \PP(\XX), \;\ff\in L_{\varphi*}(\XX,\nu),
\end{align}
where $\E^{\mu}[\ff]:= \int_{\Xt}\!\ff d\mu$, $\mathfrak{B}_{\varphi,+}(\mu\|\nu;\ff) \,{=}\, \mathfrak{B}_{\varphi,-}(\mu\|\nu;\ff) \,{=}\, 0$ iff $\mu \,{=}\, \nu$ or if $\ff$ is constant $\nu\,\text{-\,a.s.}$ 
 Then, in Proposition \ref{InformaRe}, we derive a representation  $\mathfrak{B}_{\varphi,\pm}(\mu\|\nu;\ff)$   which allows us to re-write  (\ref{AmInfo}) as 
\begin{align}\label{Csiszar}
\hat{\mathcal{K}}_{\varphi,\sff}^{\nu}\big(-\Df(\mu\|\nu)\big)\leqslant  \E^\mu[\ff] - \E^\nu[\ff] \leqslant \mathcal{K}_{\varphi,\sff}^{\nu}\big(\Df(\mu\|\nu)\big),  
\end{align}
where 
 $\mathcal{K}_{\varphi,\sff}^{\nu}(s)\rightarrow 0$, $\hat{\mathcal{K}}_{\varphi,\sff}^{\nu}(-s)\rightarrow 0$ as $s\,{\downarrow}\, 0$.  In \S\ref{Inf_IP} we develop  the bound (\ref{Csiszar}) further in the context of SDE's in order to bound $\Df(\mu_t\|\nu_t)$ in terms of differences between the  vector fields generating the dynamics in (\ref{SDE1}) and (\ref{SDE2}); this step is important for bounding the error in estimates of (Lagrangian) observables explicitly in terms of the (Eulerian) vector fields generating the true and approximate underlying dynamics. Specific bounds on path-based observables follow from the fact that $\E^{\mu_t}[f] = \E\big[f\big(\pi^\nu_\mm\hspace{.04cm}{\circ}\hspace{0.04cm}\Solm^\mm\big)\big]$, $\E^{\nu_t}[f] = \E\big[f\big(\Solm^\nu\big)\big]$, as highlighted in \S\ref{setp} and detailed in~\S\ref{SFlow}. The information about the initial conditions is propagated through the solutions of (\ref{F_Kol}) or, more explicitly, via the transition evolutions as $\mmt = \mathcal{P}^{\mm*}_{t_0,t}\,\mmo$ and $\nu_t = \mathcal{P}^{\nu*}_{t_0,t}\,\nu_{t_0}$, $\mu_t\in \PP(\XXt)$, $\nu_t\in \PP(\Xt)$ (see Definition \ref{trans_evo}).

\smallskip
Below, we present the main results of this section which are followed by a discussion of links between these general bounds and some known inequalities for specific choices of $\varphi$-divergence.     

\vspace{-0.1cm}\begin{theorem}[\textbf{\textit{Information bounds}}]\label{iNF_IN}
 Let $\mu,\nu\in \PP(\XX)$ be probability measures on a Polish space $(\XX, \Bb(\XX))$ such that  $\Df(\mu\|\nu)\,{<}\,\infty$, $\varphi\in \mathcal{C}^2(\Rp)$ is strictly convex, and it satisfies  (\ref{Normality}).  Then, for any   $\ff\in L_{\varphi*}(\XX,\nu)$ there exist  $\mathfrak{B}_{\varphi,+}, \mathfrak{B}_{\varphi,-} \in \R$ such that  
\begin{align}\label{Goal}
\mathfrak{B}_{\varphi,-}(\mu\|\nu;\ff)\leqslant \E^{\mu}[\ff] - \E^{\nu}[\ff]\leqslant \mathfrak{B}_{\varphi,+}(\mu\|\nu;\ff),
\end{align}
where 
\begin{equation}
\mathfrak{B}_{\varphi,\pm}(\mu\|\nu;\ff):= \pm\inf_{\lambda > 0}\bigg\{\frac{1}{\lambda}\int_{\Xt} \varphi^*\big(\pm\lambda \big(\ff-\E^{\nu}[\ff]\big)\big)d\nu +\frac{1}{\lambda} \Df(\mu\|\nu)\bigg\},
\end{equation}
and
\begin{align}\label{Goal2}
\mathfrak{B}_{\varphi,+}(\mu\|\nu;\ff) = \mathfrak{B}_{\varphi,-}(\mu\|\nu;\ff) = 0, \;\; \textrm{iff}\; \mu = \nu; \; \textrm{ or if $\ff$\! is constant} \; \nu\text{-a.s.}
\end{align}
\noindent {\it Proof}: \rm See Appendix \ref{app_thm_iNF_IN}; the proof utilises a repeated use of the Legendre-Fenchel transform (\ref{Phistar}), the definition of $\Df$, and Fenchel-Young inequality which follows from (\ref{Phistar}). 
\end{theorem}

\vspace{-0.3cm}
\begin{rem}\rm
The above result generalises the `goal-oriented information inequality' for KL-divergence and $\chi^2$ developed in \cite{dupuis13, dupuis16, Jing12}, to a class of all information monotone divergences. Note that, when $\varphi(u) \,{=}\, u\log u-u+1, \,u\,{>}\,0,$ the Orlicz subspace $L_{\varphi^*}(\XX, \nu)$ is simply the class of all  cumulant generating functions (a.k.a.~logarithmic moment generating functions). The results in \cite{dupuis13, dupuis16, Jing12} are based on the regularity of cumulant generating functions. Our generalisation relies on a convex-analytic approach under the normality conditions (\ref{Normality}) imposed on $\varphi$ or $\varphi^*$.
\end{rem}

\vspace{-0.3cm}\begin{prop}[\textbf{\textit{Representation formula for}} $\mathfrak{B}_{\varphi,\pm}(\mu\|\nu;\ff)$]\label{InformaRe}
Given the bound (\ref{Goal}) and the assumptions of Theorem \ref{iNF_IN},  consider the convex function $\lambda\mapsto \mathcal{G}_{\varphi,\nu}(\lambda; \ff)$, $\lambda\in \Rp$,   defined~by 
 \begin{align*}
\mathcal{G}_{\varphi,\nu}(\lambda; \ff) = \int_{\Xt}\varphi^*\Big(\lambda(\ff-\E^{\nu}[\ff] )\Big)d\nu,
\end{align*}
where $\ff\in L_{\varphi*}(\XX,\nu)$, $\ff\neq \E^{\nu}[\ff], \;\; \nu\textrm{-a.s.}$
\btm[leftmargin = 0.7cm]

\vspace{.1cm}
\item[(1)]  Then
\begin{equation*}
\mathfrak{B}_{\varphi,+}(\mu\|\nu;\ff) = \mathcal{K}_{\varphi,\sff}^\nu\big(\Df(\mu\|\nu)\big),
\end{equation*}

\vspace{.1cm}
\noindent where 
$\mathcal{K}_{\varphi,\sff}^\nu(R^2)=\inf\big\{\delta\geqslant 0: \mathcal{G}^*_{\varphi,\nu}(\delta;\ff)>R^2\big\}$, and  $\mathcal{G}^*_{\varphi,\nu}(\delta;\ff)$ is the Legendre--Fenchel conjugate of $\mathcal{G}_{\varphi,\nu}(\lambda;\ff)$ defined by 
\vspace{.1cm}\begin{align*}
\mathcal{G}^*_{\varphi,\nu}(\delta;\ff) = \sup_{\lambda\geqslant0}\big\{\lambda \delta -\mathcal{G}_{\varphi,\nu}(\lambda;\ff)\big\}.
\end{align*}

\vspace{.1cm}
\noindent Similarly, the lower bound $\mathfrak{B}_{\varphi,-}$ admits the representation 
\vspace{.1cm}\begin{align*}
\mathfrak{B}_{\varphi,-}(\mu\|\nu;\ff) = \hat{\mathcal{K}}_{\varphi,\sff}^\nu(-\Df(\mu\|\nu)),
\end{align*}
where $\hat{\mathcal{K}}_{\varphi,\sff}^\nu(-R^2) \,{=}\,\sup\{\delta\geqslant 0: -\mathcal{G}^*_{\varphi,\nu}(-\delta;\ff){<}-R^2\}\,{=}-\inf\{\delta\geqslant 0: \mathcal{G}^*_{\varphi,\nu}(-\delta;\ff)>R^2\}$.

\vspace{.3cm}\item[(2)] If in addition  $\Df(\mu\|\nu)<\infty$, then 
\vspace{.1cm}\begin{align}
\begin{cases}\label{Df_rep_bnd} 
\mathfrak{B}_{\varphi,+}(\mu\|\nu;\ff) = \nabla\mathcal{G}_{\varphi,\nu}\left(\tilde{\mathcal{H}}_{+,\sff}^{-1}\big(\Df(\mu\|\nu)\big);\ff\right),\\[.4cm]
\mathfrak{B}_{\varphi,-}(\mu\|\nu;\ff) = \nabla\mathcal{G}_{\varphi,\nu}\left(-\tilde{\mathcal{H}}_{-,\sff}^{-1}\big(\Df(\mu\|\nu)\big);\ff\right),
\end{cases}
\end{align}

\vspace{.1cm}
\noindent where $\tilde{\mathcal{H}}^{-1}_{+,\sff}$ is the pseudo-inverse of a function $\mathcal{H}_{+,\sff}: [0, \infty)\rightarrow \R$ defined by 
\begin{align}\label{HFxn}
\mathcal{H}_{+,\sff}(\lambda) = -\mathcal{G}_{\varphi,\nu}(\lambda;\ff)+\lambda\nabla\mathcal{G}_{\varphi,\nu}(\lambda;\ff), 
\end{align}
which is strictly increasing on $(0,\infty)$, while  $\tilde{\mathcal{H}}^{-1}_{-,\sff}$ is the pseudo-inverse of $\mathcal{H}_{-,\sff}: (-\infty, 0]\rightarrow \R$  which is  strictly decreasing on $(-\infty, 0)$.
\etm
\end{prop}
\vspace{-0.1cm}\noindent {\it Proof.} See Appendix \ref{InformaRe_ap}; the proof utilises  Legendre--Fenchel transform (\ref{Phistar}), Fenchel--Young inequality, and the implicit function theorem. 

\bigskip
The following explicit representation of the $\varphi$-information bounds $\mathfrak{B}_{\varphi,\pm}$ in~(\ref{Goal}) can be deduced  via linearisation about $\Df(\mu\|\nu)=0$. One, but not the only, case when this result is useful is when the probability measure $\nu$ is a small perturbation of~$\mu$ (in the sense that $\Df(\mu\|\nu)\ll1$). 

\vspace{-0.1cm}\begin{cor}[\textbf{\textit{Linearisation of $\varphi$-information bounds}}]\label{UQ_lineariz}
Let $\mu, \nu\in \PP(\XX)$  and consider  a strictly convex function $\varphi\in \mathcal{C}^2(\Rp)$  satisfying (\ref{Normality}). If  $\ff\in L_{\varphi*}(\XX,\nu)$ with $\E^\nu[\ff]\neq \ff,\;\nu$-a.s. and $\Df(\mu\|\nu)<\infty$, then 
\vspace{.0cm}\begin{align}
\mathfrak{B}_{\varphi,\pm}(\mu\|\nu;\ff) = \pm\sqrt{\vphantom{\big\|}2\nabla^2\varphi^*(0)\text{Var}_\nu(\ff)}\sqrt{\Df(\mu\|\nu)}+\mathcal{O}(\Df(\mu\|\nu)),
\end{align} 
and 
\begin{align}\label{lin_phi}
\big\vert \E^\mu[\ff]-\E^\nu[\ff]\big\vert\leqslant \sqrt{\vphantom{\big\|}2\nabla^2\varphi^*(0)\text{Var}_\nu(\ff)}\sqrt{\Df(\mu\|\nu)}+\mathcal{O}(\Df(\mu\|\nu)).
\end{align}

\vspace{.2cm}
\noindent The term $\mathcal{O}(\Df(\mu\|\nu))$ can be resolved further if $\varphi^*\in \mathcal{C}^{N+2}(\R)$ for all $N\geqslant 1.$
\end{cor}
\vspace{-0.cm}\noindent {\it Proof.} See Appendix \ref{UQ_lineariz_ap}; the proof utilises the Legendre--Fenchel transform, strict convexity of $\varphi$ and $\varphi^*$, and the standard Taylor expansion around $\Df(\mu\|\nu)=0$.

\begin{rem}\rm A number of observations deserves a mention in regards to Theorem~\ref{iNF_IN}. \mbox{}
\begin{itemize}[leftmargin=0.7cm]
\item[{\bf (i)}] The information inequality  (\ref{Csiszar}) generalises  the inequality proved long ago by Csisz\'ar~\cite{Csiszar72}  in terms of the total variation distance $\Vert \mu-\nu\Vert_{\textsc{tv}}$ but it applies to a larger class of observables $\ff\in L_{\varphi*}(\XX;\nu)$ and it is tighter; see Corollary~\ref{csipin}. In our notation, Csisz\'ar's result concerns  the  existence of $\tilde{\mathcal{K}}: \Rp\rightarrow\Rp$,  $\tilde{\mathcal{K}}(s)\rightarrow 0$ as $s\rightarrow 0$, such that
\begin{align}\label{Cszar2}
\Vert \mu-\nu\Vert_{\textsc{tv}}\leqslant \tilde{\mathcal{K}}\big(\D_\textsc{kl}(\mu\|\nu)\big).
\end{align}
Noting that the variational representation of $\Vert \mu-\nu\Vert_{\textsc{tv}}$ is given by 
\begin{align}\label{Var_TV}
\Vert \mu-\nu\Vert_{\textsc{tv}} = \sup_{\Vert \sff\Vert_{\infty}\leqslant 1}\Big\{ \E^\mu[\ff]-\E^\nu[\ff]\Big\},
\end{align}
and combining (\ref{Cszar2}) with  (\ref{Var_TV}) yields a version of  the Csisz\'ar--Pinsker--Kullback inequality
 \begin{align}\label{Cszar3}
\big\vert \E^\mu[\ff] - \E^\nu[\ff]\big\vert\leqslant \Vert \ff\Vert_{\infty} \,\tilde{\mathcal{K}}\big(\D_\textsc{kl}(\mu\|\nu)\big).
\end{align}
As shown in Corollary \ref{csipin}, the bound (\ref{Goal}) can be symmetrised and it yields the inequality 
\begin{align}\label{ours}
\big\vert \E^\mu[\ff] - \E^\nu[\ff]\big\vert\leqslant \Vert \ff\Vert_{\infty} \,\widetilde{\mathcal{K}}_{\varphi,f}^\nu\big(\D_\varphi(\mu\|\nu)\big).
\end{align}
One can easily  verify that $L^\infty(\XX; \nu)\subset L_{\varphi*}(\XX,\nu)$,  which implies that (\ref{ours}) holds for a larger class of observables than (\ref{Cszar3}). 

\medskip
\item[{\bf (ii)}] The information bound (\ref{Goal}) also generalises another well-known bound on the error in observables in terms of $\chi^2$-\,divergence, which is given by\,\footnote{\,The bound (\ref{robbins}) leads to the {\it Hammersley-Chapman-Robbins inequality} \cite{chapman51}  when $f$ is taken to be an unbiased estimate of some functional $\mathfrak{f}(\nu)$, i.e.,  $\E^\nu[f] = \mathfrak{f}(\nu)$, which is widely used in statistical estimation (e.g.,~\cite[p. 114]{Lehmann}).} 
\vspace*{.1cm}\begin{align}\label{robbins}
\qquad \big\vert \E^\mu[\ff]-\E^\nu[\ff]\big\vert \leqslant \sqrt{\text{Var}_\nu(\ff)}\sqrt{\chi^2(\mu\|\nu)}, \quad \qquad \text{Var}_\nu(\ff) := \E^\nu\big[(\E^\nu[\ff] -\ff)^2\big]. 
\end{align}

\smallskip
\noindent The bound (\ref{robbins}) is more useful than (\ref{Cszar3}) or (\ref{ours}) when  $\ff\in L^2(\XX;\nu)$ and not necessarily $\ff\in L^{\infty}(\XX;\nu)$.  The link between (\ref{Goal}) and (\ref{robbins}) is discussed in Corollary \ref{chi_link}. Furthermore, we note that for $\varphi(u) = (u-1)^2, \; u>0$, so that $\varphi^*(\xi) = \frac{1}{4}\xi^2+\xi$, the bound (\ref{lin_phi}) in Corollary~\ref{UQ_lineariz} becomes an equality and $\mathcal{O}(\Df(\mu\|\nu))\equiv 0$; this highlights the fact that, while (\ref{lin_phi}) is particularly useful when $\Df(\mu\|\nu)\ll1$, the result in (\ref{lin_phi}) is not restricted to such cases and the terms $\mathcal{O}(\Df(\mu\|\nu))$ can be resolved further for sufficiently regular $\varphi^*$. 

\medskip
 \item[{\bf (iii)}] In \cite{MBUda18} an information-theoretic/probabilistic  measure of  average finite-time expansion rates of solutions  of non-autonomous SDEs on $\XX=\Xt$, with $\Xt = \R^d$ or  $\Xt\,{=}\,\bar{\T}^\ell$ (see \S\ref{setp}),  was defined via the  KL-divergence (see Table~\ref{phidlist}).  Here, we focus on ODE's for brevity. It turns out that for the probability measure $\mu_t^x$ evolving from  $\mu^x_{t_0}$ concentrated on a neighbourhood of $x\in \Xt$  the map $x\mapsto\DK\big({\mu}^x_t\|{\mu}^x_{t_0}\big)$ is  linked to so-called  finite-time Lyapunov functionals 
$$\Lambda^{t-t_0}_{t_0}(x,y) = \scaleobj{.8}{\frac{1}{|t-t_0|}}\log \frac{|Y_{t_0,t}^x|}{|y|}, \qquad y\ne0,\;y\in \Xt,$$
which are commonly used to assess  the growth of perturbation about solutions of ODEs with the initial condition $x\in \Xt$; i.e., $Y^{t_0,x}_t = D\phi_{t_0,t}(x+y)-D\phi_{t_0,t}(x)$. It is worth noting that, similar to (i) and (ii), the bound on the Lyapunov exponents derived in \cite{MBUda18} and given by 
\vspace*{.2cm}\begin{align}\label{Lyapfuncti}
\big| \E^{\mu^x_{t_0}}[\Lambda_{t_0}^{t-t_0}(x,\ccdot)]\big|\leqslant \scaleobj{.8}{\frac{1}{|t-t_0|}}\,\widetilde{\mathcal{K}}^{\mu^x_{t_0}}\big(\DK\big({\mu}_t^x\|{\mu}^x_{t_0}\big)\big), \qquad t>t_0, \;t, t_0\in \Ic,
\end{align}

\medskip
\noindent is an instance of a `coarsened' $\varphi$-information inequality  (\ref{Csiszar}) for $\varphi(u) = u\log u-u+1$,  $f(u)=\log u$, $u>0$, $\mu = \mu_t^x$, $\nu = \mu^x_{t_0}$, and $\widetilde{\mathcal{K}}^{\mu^x_{t_0}}\geqslant $, 
since (\ref{Lyapfuncti}) can be rewritten as 
\vspace*{.1cm}\begin{align*}
\big| \E^{\mu^x_{t_0}}[\Lambda_{t_0}^{t-t_0}(x,\cdot)]\big| &= \left|\int_\Xt \Lambda_{t_0}^{t-t_0}(x,y)\mu^x_{t_0}(dy)\right|\notag\\[.2cm]
&= \scaleobj{.8}{\frac{1}{|t-t_0|}}\left| \int_\Xt \log|y^x_t|\mu^x_{t_0}(dy)- \int_\Xt \log|y|\mu^x_{t_0}(dy)\right|\notag\\[.2cm]
&= \scaleobj{.8}{\frac{1}{|t-t_0|}}\left| \int_\Xt \log|y|\mu^x_{t}(dy)- \int_\Xt \log|y|\mu^x_{t_0}(dy)\right|\notag\\[.4cm]
&=\scaleobj{.8}{\frac{1}{|t-t_0|}}\big| \E^{\mu^x_{t}}\big[\log |y|\big] -\E^{\mu^x_{t_0}}\big[\log |y|\big]\big| \leqslant \scaleobj{.8}{\frac{1}{|t-t_0|}}\,\widetilde{\mathcal{K}}^{\mu^x_{t_0}}\big(\DK\big({\mu}_t^x\|{\mu}^x_{t_0}\big)\big),
\end{align*}

\medskip
\noindent where we used the properties of solutions of the Liouville equation; see, e.g., \cite[Theorem~4.8]{Pavliotis}. Tighter versions of (\ref{Lyapfuncti}), based on (\ref{Goal}), will be considered in future applications.
 
\end{itemize}
\end{rem}

\vspace{-0.4cm}

\begin{cor} \label{csipin}
For $\varphi(u) = u\log u-u+1$, $u>0$,  and $\ff\in L^\infty(\XX;\nu)$, the $\varphi$-information bound~(\ref{Goal}) implies  the Csisz\'ar--Pinsker--Kullback inequality (\ref{Cszar3}); however, the  bound (\ref{Goal}) expressed in terms of $\mathfrak{B}_{\varphi,\pm}(\mu\|\nu;\ff)$ is tighter.

\smallskip
\noindent Proof. \rm
We combine the results from Theorem \ref{iNF_IN} and Proposition \ref{InformaRe}. From the definition  of the convex function $\mathcal{G}_{\varphi,\nu}(\lambda;\ff)$ in Proposition \ref{InformaRe}, we have that  for $\lambda\geqslant 0$
\begin{align}\label{gradG_1}
\nabla_\lambda \mathcal{G}_{\varphi,\nu}(\pm\lambda;\ff) = \pm\int_{\Xt}\Big(\ff-\E^\nu(\ff)\Big)\nabla\varphi^*\Big(\pm\lambda\big(\ff-\E^\nu(\ff)\big)\Big)d\nu.
\end{align}
If $\ff\in L^\infty(\XX;\nu)$, then based on (\ref{gradG_1}) and the representation formulas (\ref{Df_rep_bnd}),  we have 
\vspace*{.2cm}\begin{align*}
\mathfrak{B}_{\varphi,+}(\mu\|\nu;\ff) &\leqslant \|\ff\|_{\infty}\int_{\Xt} 2\nabla\varphi^*\left(\tilde{\mathcal{H}}_{+,\sff}^{-1}\big(\Df(\mu\|\nu)\big)\right)d\nu,\\[.2cm]
\mathfrak{B}_{\varphi,-}(\mu\|\nu;\ff) &\geqslant -\|\ff\|_{\infty}\int_{\Xt} 2\nabla\varphi^*\left(\tilde{\mathcal{H}}_{-,\sff}^{-1}\big(\Df(\mu\|\nu)\big)\right)d\nu.
\end{align*}

\smallskip
\noindent The normality conditions (\ref{Normality}) imply that  there exists $\widetilde{\mathcal{K}}_{\varphi,\sff}^\nu: (-\infty, \infty)\rightarrow\Rp$ depending on $\ff\!, \nu$ and $\varphi$ with $\widetilde{\mathcal{K}}_{\varphi, \sff}^\nu(s)\rightarrow 0$ as $s\downarrow 0$,  such that 
\vspace*{.1cm}\begin{align}\label{BetterB}
\mathfrak{B}_{\varphi,+}(\mu\|\nu;\ff) \leqslant \|\ff\|_{\infty}\,\widetilde{\mathcal{K}}_{\varphi,\sff}^\nu\big(\Df(\mu\|\nu)\big), \quad \mathfrak{B}_{\varphi,-}(\mu\|\nu;\ff) \geqslant -\|\ff\|_{\infty}\,\widetilde{\mathcal{K}}_{\varphi,\sff}^\nu\big(\Df(\mu\|\nu)\big).
\end{align}

\smallskip
\noindent For $\varphi(u) = u\log u-u+1$, $u>0$, combining (\ref{Goal}) with (\ref{BetterB}) leads   to the Csisz\'ar--Pinsker--Kullback inequality (\ref{Cszar3}); the bound (\ref{Goal}) utilising $\mathfrak{B}_{\varphi,\pm}(\mu\|\nu;\ff)$ is tighter by construction. \qed
\end{cor}

\begin{cor}\label{chi_link} For $\varphi(u) = (u-1)^2, \; u>0$, and $\ff\in L^2(\XX;\nu)$, the $\varphi$-information bound (\ref{Goal}) implies  the   bound (\ref{robbins}); however, the  bound (\ref{Goal}) expressed in terms of $\mathfrak{B}_{\varphi,\pm}(\mu\|\nu;\ff)$ is tighter.

\medskip
\noindent Proof. \rm
Similar to the proof of Corollary \ref{csipin}, we combine the results from Theorem \ref{iNF_IN} and Proposition \ref{InformaRe}.
Given the convex function $\mathcal{G}_{\varphi,\nu}(\lambda;\ff)$ in Proposition \ref{InformaRe}, we have that  for $\lambda\geqslant 0$
\vspace*{.2cm}\begin{align}\label{gradG_2}
\nabla_\lambda \mathcal{G}_{\varphi,\nu}(\lambda;\ff) = \int_{\Xt}\Big(\ff-\E^\nu(\ff)\Big)\nabla\varphi^*\Big(\lambda\big(\ff-\E^\nu(\ff)\big)\Big)d\nu.
\end{align}

\smallskip
\noindent If $\ff\in L^2(\XX;\nu)$, then based on (\ref{gradG_2}) and (\ref{Df_rep_bnd}),  we  have by H\"older's inequality
\vspace*{.3cm}\begin{align*}
\mathfrak{B}_{\varphi,+}(\mu\|\nu;\ff) &\leqslant \scaleobj{1}{\sqrt{\text{Var}_\nu(\ff)}}\left(\int_{\Xt} \left(\nabla\varphi^*\left(\tilde{\mathcal{H}}_{+,\sff}^{-1}\big(\Df(\mu\|\nu)\big)\right)\right)^2d\nu\right)^{1/2},\\[.3cm]
\mathfrak{B}_{\varphi,-}(\mu\|\nu;\ff) &\geqslant  -\scaleobj{1}{\sqrt{\text{Var}_\nu(\ff)}}\left(\int_{\Xt} \left(\nabla\varphi^*\left(\tilde{\mathcal{H}}_{-,\sff}^{-1}\big(\Df(\mu\|\nu)\big)\right)\right)^2d\nu\right)^{1/2},
\end{align*}

\medskip
\noindent where $\text{Var}_\nu(\ff) := \E^\nu\big[\big(\E^\nu[\ff\hspace{.03cm}] -\ff\big)^2\big]$. The normality conditions (\ref{Normality}) imply that  there exists $\widetilde{\mathcal{K}}_{\varphi,\sff}^\nu: (-\infty, \infty)\rightarrow\Rp$ depending on $\ff\!, \nu$ and $\varphi$ with $\widetilde{\mathcal{K}}_{\varphi, \sff}^\nu(s)\rightarrow 0$ as $s\downarrow 0$, such that 
\vspace*{.2cm}\begin{align}\label{BetterB2}
\mathfrak{B}_{\varphi,+}(\mu\|\nu;\ff) \leqslant \scaleobj{.9}{\sqrt{\text{Var}_\nu(\ff)}}\,\widetilde{\mathcal{K}}_{\varphi,\sff}^\nu\big(\Df(\mu\|\nu)\big), \quad \mathfrak{B}_{\varphi,-}(\mu\|\nu;\ff) \geqslant -\scaleobj{.9}{\sqrt{\text{Var}_\nu(\ff)}}\,\widetilde{\mathcal{K}}_{\varphi,\sff}^\nu\big(\Df(\mu\|\nu)\big).
\end{align}

\medskip
\noindent For $\varphi(u) = (u-1)^2, \; u>0$,  s.t.~$\varphi^*(\xi) = \frac{1}{4}\xi^2+\xi$, the bounds (\ref{BetterB2}) lead  to (\ref{robbins}); however, the $\varphi$-information bound (\ref{Goal}) utilising $\mathfrak{B}_{\varphi,\pm}(\mu\|\nu;\ff)$ is tighter by construction, as claimed.\qed

\end{cor}


\section{Preliminaries on stochastic flows and martingale solutions}\label{SFlow}

In the remainder of this work, several technical concepts will be necessary; in particular, the notions of a stochastic flow and a Lebesgue a.e.~martingale solutions.  In order to make the presentation self-contained, we recall some important definitions and background results. 

\smallskip
Throughout this section $\XX$ is a smooth $n$-dimensional  differentiable manifold with no  boundary, and $\kappa\in \PP(\XX)$ denotes a probability measure on $(\XX, \Bb(\XX))$; i.e., the exposition here  is general and it is not restricted to  the specific manifolds $\XXt$ and $\Xt\subseteq \XXt$ associated with the original dynamics and its approximation discussed in \S\ref{intro}-\ref{setup_sec}, which are considered subsequently in~\S\ref{Main_Bounds}.  

\begin{definition}[\textbf{\textit{Stochastic flow}} \cite{Kunitanote, Kunitabook}]\label{stchf}\rm 
For  $\phi_{s,t}:\XX\times\Om_\mathcal{H}\rightarrow \XX$, $s,t\in \Ic=[t_0, t_0+T]$, let $x\mapsto \phi_{s,t}(x,\om)$  be a continuous random field on the probability  space $\big(\Om_\mathcal{H}, \mathcal{H}, \Ms\big)$. The map $\phi_{s,t}(\,\cdot\,,\om)$  defines a {\it stochastic flow of homeomorphisms} if there exists a null set $\mathcal{N}\subset \Om_\mathcal{H}$ such that for any $\om\notin \mathcal{N},$ the family of continuous maps $\big\{\phi_{s,t}(\,\cdot\,,\om)\,{:}\;\, s,t\in\Ic\big\}$~satisfies the following: 
\btm[leftmargin = 1cm]
\item[(i)] $\phi_{s,t}(\ccdot,\om) = \phi_{u,t}\big(\phi_{s,u}(\ccdot\,,\om),\om\big)$ holds for any $s,t,u\in \Ic$,

\vspace{.1cm}\item[(ii)] $\phi_{s,s}(\ccdot\,,\om) = {\rm id}_{\XX}$, for all $s\in \Ic$,

\vspace{.1cm}\item[(iii)] the map $\phi_{s,t}(\ccdot,\om): \XX\rightarrow \XX$ is a homeomorphism for any $s,t\in \Ic$.
\etm
The map $x\mapsto \phi_{s,t}(x,\om)$ defines  a {\it stochastic flow of $\mathcal{C}^l$-diffeomorphisms}  if it  is $l$-times continuously differentiable w.r.t.~$x\in \XX$  $\forall\,s,t\in \Ic$, a.a.~$\om\in \Om$, with continuous derivatives in $(s,t).$ 
\end{definition}

\begin{definition}[\textbf{\textit{Transition probability kernel} }\cite{bogachev10,Roc-Kry}]\label{martrans}\rm
A map $(s,x,t,B)\mapsto P(s,x;t,B)$, $B\in \Bb(\XX)$,   is called a {\it transition probability kernel} on $(\XX, \Bb(\XX))$ if the following hold:
\begin{itemize}[leftmargin= .9cm]
\item[(i)]  $B\mapsto P(s,x;t,B)$ is a probability measure on $\Bb(\XX)$ for any $s,t\in \Ic,x\in \XX$. 

\vspace{.1cm}\item[(ii)]  $x\mapsto P(s,x;t,B)$ is $\Bb(\XX)$-measurable for any  $s,t\in \Ic,B\in \mathcal{B}(\XX)$.
\vspace{.1cm}\item[(iii)] The following holds for any $s,t,u\,{\in}\, \Ic$, $s\leqslant u\leqslant t$, and for all $x\,{\in}\, \XX$, $B\in \Bb(\XX)$  
\begin{equation}\label{ChK}
P(s,x;t,B) = \int_{\mathcal{M}} P(u,y;t,B)P(s,x;u,dy). 
\end{equation}

\item[(iv)] $P(s,x;s,B) = \I_{B}(x)$ for all $x\in \XX$ and $B\in \Bb(\XX)$.
\end{itemize}
\end{definition}

\begin{definition}[\textbf{\textit{Transition evolutions and their duals}}]\label{trans_evo}\rm
Any transition probability kernel on $(\XX, \Bb(\XX))$ defines a family of linear {\it transition evolution operators}  $(\mathcal{P}_{s,t})_{t\geqslant s}$,  $s,t\in \Ic$, as follows 
\begin{align}\label{P}
\big(\mathcal{P}_{s,t}f\big)(x): = \int_{\mathcal{M}}f(y)P(s,x;t,dy), \qquad f\in \mathbb{M}(\XX).
\end{align}
The {dual} of $\mathcal{P}_{s,t}$, which acts on probability measures $\kappa\in \PP(\XX)$, is defined by 
\begin{align}\label{P*}
\big(\mathcal{P}^*_{s,t}\,\kappa\big)(B) := \int_{\mathcal{M}}P(s,x;t, B)\kappa(dx),  \qquad B\in \Bb(\XX).
\end{align}
For any $s, t\in \Ic$, $s\leqslant t$,  we define $\kappa_{s,t} := \mathcal{P}^*_{s,t}\,\kappa_{t_0,s}$; for $s=t_0$ we simply write  $\kappa_{t} := \mathcal{P}^*_{t_0,t}\,\kappa_{t_0}$. 

\end{definition}

\begin{rem}\rm A stochastic flow,  $\phi_{s,t}\!: \XX\times\Omega_\mathcal{H}\rightarrow\XX$, generates a transition probability kernel via 
\begin{align}\label{P_phi}
P(s,x;t,B) := \Ms\big(\{\om\in \Om_\mathcal{H}\!:\, \phi_{s,t}(x,\om)\in B \}\big), \quad  s,t\in \Ic, \;t\geqslant s, \quad B\in\Bb(\XX).
\end{align}
Consequently, the transition evolution~(\ref{P}) and its dual (\ref{P*}) are induced by $\phi_{s,t}$ as follows  
\begin{align}\label{PP}
(\mathcal{P}_{s,t}f)(x) &= \int_{\XX}f(y)P(s,x;t,dy)=\mathbb{E}\big[f\big(\phi_{s,t}(x)\big)\big], \quad t\geqslant s, \quad f\in \mathbb{M}(\XX),\\[.1cm]
(\mathcal{P}_{s,t}^*\kappa_s)(B) &= \int_{\XX}\mathbb{E}\big[\I_B\big(\phi_{s,t}(x)\big)\big]\kappa_s(dx), \hspace{2.35cm} t\geqslant s, \quad B\in\Bb(\XX).
\end{align}
\end{rem}

\begin{theorem}[\textbf{\textit{Representation of solutions to SDE's via stochastic flows}}]\label{SDE_flow}\mbox{}

\noindent Let  $b^\kappa(t,\ccdot):\XX\rightarrow\XX$ and $\sigma^\kappa(t,\ccdot): \XX\rightarrow \XX^{\otimes m}$, $m\geqslant 1$, $t\in \Ic$, be measurable $t$-continuous functions satisfying 
\begin{alignat}{2}
&|\langle b^\kappa(t,x),x\rangle|+\|\sigma^\kappa(t,x)\|^2_\textsc{hs}\leqslant C(1+|x|^2),& \hspace{1cm} &x\in \XX, t\in \Ic, \;C>0, \label{k_growth}\\[.1cm]
&\vert b^\kappa(t, x) - b^\kappa(t,y)\vert + \Vert \sigma^\kappa(t,x)-\sigma^\kappa(t,y)\Vert_{\textsc{hs}}\leqslant L_K\vert x-y\vert, && x,y\in K, t\in \Ic, \;L_K>0, \label{Lipsch}
\end{alignat}
where $K$ is any compact subset of $\XX$. 
Then, for  an $m$-dimensional Wiener process $W_t$ on the Wiener space $(\Om,\mathcal{F},\p)$, the It\^o SDE  
\begin{align}\label{SDE}
dX^\kappa_t = b^\kappa\big(t,X^\kappa_t\big)dt+ \sigma^\kappa\big(t,X^\kappa_t\big) dW_{t-t_0}, \qquad \;X^\kappa_{t_0} \sim \kappa_{t_0},\; \kappa_{t_0}\in \PP(\XX),
\end{align}
has a unique global solution on $\Ic=\big[t_0, \,t_0+T\big)$.
For the filtration $(\F^{\kappa}_t)_{t\geqslant t_0}$ induced by a given version of $W_t$ and $X_{t_0}^\kappa$ independent of $\F^{\kappa}_\Ic$ s.t.~$\E|X^\kappa_{t_0}|^2<\infty$, there exists a unique strong solution of (\ref{SDE}) which is adapted to $(\F^{\kappa}_t)_{t\geqslant t_0}$ 
and $\E|X^\kappa_{t}|^2<\infty$. The solutions to (\ref{SDE})  can be represented as a stochastic flow of homeomorphisms, i.e.,  
\begin{equation}\label{SDE_sol}
X^\kappa_t(x,\om) = \phi^\kappa_{t_0,t}(x,\om) \qquad \p\,\textrm{-\,a.s.,} \;\;\kappa_{t_0}\textrm{-\,a.a.}\;x\in \XX,
\end{equation}
Moreover,   $\phi_{t_0,t}^\kappa(\ccdot,\om)$ in  (\ref{SDE_sol}) is a $C^l$-diffeomorphism on $\XX$ over $\Ic$ if, in addition to (\ref{k_growth})-(\ref{Lipsch}), the $l$-th derivatives of $b^\kappa(t,\ccdot)$, and the $l$-th derivatives of $\sigma^\kappa(t,\ccdot)$ are locally bounded and $\delta$-H\"older continuous for all $t\in \Ic$. 
The same holds for solutions of the Stratonovich counterpart of (\ref{SDE}) with $l$~derivatives of $\brr(t,\ccdot)$ and $l+1$ derivatives of $\sigma(t,\ccdot)$ $\delta$-H\"older continuous for all $t\in \Ic$.
\end{theorem}
\noindent {\it Proof.} The first part is standard and well known; see, e.g., \cite{oksendal} or \cite[Theorem 3.4.6]{Kunitabook}. For the representation (\ref{SDE_sol}) see, e.g., \cite[Theorems 4.7.1, 3.4.6]{Kunitabook} for homeomorphisms and  \cite[Theorem 4.7.2, 3.4.6]{Kunitabook}) for the diffeomorphism representation.

\begin{rem}\label{cf_remk}\rm
The function spaces containing $b^\kappa(t,\ccdot)$,  $\sigma^\kappa_k(t,\ccdot)$, with $\sigma^\kappa_k$ the columns of $\sigma^\kappa$,  for which (\ref{SDE_sol}) holds are, respectively, $\tilde{\mathcal{C}}^{l,\delta}(\XX,\XX)$, $\bar{\mathcal{C}}^{l,\delta}(\XX,\XX)$, $l\in \mathbb{N}_0$, $0<\delta\leqslant 1$ (see the Glossary).  Consequently, the solutions of (\ref{SDE})  are represented by a flow of $C^l$-diffeomorphisms for $t\mapsto b^\kappa(t,\ccdot)$, $t\mapsto \sigma^\kappa(t,\ccdot)$ continuous and integrable over $\Ic$ and such that 
\begin{equation}\label{reg_growth_cond}
b^\kappa(t,\ccdot)\in \tilde{\mathcal{C}}^{l,\delta}(\XX,\XX), \qquad \sigma^\kappa_k(t,\ccdot)\in \bar{\mathcal{C}}^{l,\delta}(\XX,\XX), \qquad l\geqslant 2, \; k=1,\dots,m, \; t\in \Ic,
\end{equation}
where $\brr^\kappa_i = b^\kappa-\frac{1}{2}\sigma^\kappa_{jk}\partial_{x_j}\sigma^\kappa_{ik}$ in the Stratonovich counterpart of (\ref{SDE}).
If the conditions of Theorem~\ref{SDE_flow} hold,  finiteness of the $n$-th moment of (\ref{SDE_sol}) for a bounded $\Ic\subset \R$  can be shown in a standard fashion by utilising It\^o's formula (e.g., \cite[Theorem 3.4.6]{Kunitabook}); for unbounded $\Ic$ additional `dissipative' growth constraints may have to be imposed; e.g., \cite{MBUda20}.  
\end{rem}

\bn [\textbf{\textit{Solution of forward Kolmogorov equation}} \cite{Stroock79,Figali}]\label{weak_kol}\rm
Consider the forward Kolmogorov equation (with ${b}^\kappa$ the Stratonovich-corrected drift (\ref{gen}) and $a^\kappa = \sigma^\kappa (\sigma^\kappa)^*$\,)
\begin{align}\label{PDE}
\partial_t\kappa_t = \mathcal{L}^{\kappa*}_t \kappa_t = -\sum_{i=1}^n\partial_{x_i}({b}_i^\kappa\kappa_t)+\frac{1}{2}\sum_{i,j=1}^n\partial_{x_ix_j}(a_{ij}^\kappa\kappa_t) , \;\; \textrm{in }\; \Ic\times\XX,\quad 
\kappa_{t_0} \,{\in}\, \PP(\XX). 
\end{align}
We say that $\kappa_t \,{\in}\, \PP(\XX)$ is the {\it time-marginal measure solution} to the forward Kolmogorov equation~(\ref{PDE}) in the distributional (or weak$-^*$) sense if for all $t\in \Ic$, and any $ f\in \mathcal{C}^{\infty}_c(\XX)$
\begin{align*}
\frac{d}{dt}\int_{\XX}f(x)\kappa_t(dx)= \!\int_{\XX} \!(\mathcal{L}^\kappa_t f)(x)\kappa_t(dx) =\! \int_{\XX}\!\bigg(\!\sum_i^n {b}^\kappa_i(t,x)\partial_{x_i}f(x)+\frat\!\sum_{i,j=1}^n \!a_{ij}^\kappa(t,x)\partial_{x_ix_j}f(x)\bigg)\kappa_t(dx), 
\end{align*}
 and $\kappa_t$ converges narrowly to $\kappa_{t_0}$ as $t\downarrow t_0$; i.e., 
\begin{align}
\lim_{t\,\downarrow \,t_0}\int_{\XX}f(x)\kappa_t(dx) = \int_{\XX}f(x)\kappa_{t_0}(dx) \qquad \forall\,f\in \mathcal{C}^\infty_c(\XX).
\end{align}
Given that  (\ref{PDE}) is in divergence form, it is well-posed 
(e.g.,~\cite[p.~111]{Figali}) provided that, for any relatively compact subset $K\subseteq \XX$,
\begin{align}\label{well_pos_kolm}
\int_\Ic\int_{K}\left(|{b}^\kappa(t,x)|+\|\sigma^\kappa(t,x)\|^2_{\textsc{hs}}\right)\kappa_t(dx)dt<\infty.
\end{align} 
\en

\begin{rem}\label{Feller}\rm
The following well-established facts (e.g.,~\cite{Yor, Str03, Stroock79,MBUda20}) will be needed in subsequent proofs. By assumptions on the coefficients $(b^\kappa,\sigma^\kappa)$ in  (\ref{SDE}), it can be shown that $\mathcal{P}_{s,t}f\in \mathcal{C}_c(\XX)$ for $f\in \mathcal{C}_c(\XX)$, i.e., $(\mathcal{P}_{s,t})_{t\geqslant s}$, $t,s\,{\in}\,\Ic$, is a Feller evolution\footnote{\,The family of operators $(\mathcal{P}_{s,t})_{t\geqslant s}$, $t,s\,{\in}\,\Ic$ on $\XX$ is a Feller evolution if $\mathcal{P}_{s,t}: \mathcal{C}_\infty(\XX)\rightarrow\mathcal{C}_\infty(\XX)$ for all $t,s\,{\in}\,\Ic$.}. Specifically, there exists  $L_{b^\kappa,\sigma^\kappa}\,{<}\,\infty$ depending on $\sup_{x\in \XX,\; s\leqslant r\leqslant t}\big(1+\vert x\vert^2\big)^{-1}\left[\vert \langle b^\kappa(s,x),x\rangle\vert\vee \Vert \sigma^\kappa(s,x)\Vert^2_{\textsc{hs}}\right]$, such that 
 \btm[leftmargin = 0.7cm]
\item[(i)] For each  $f\in \mathcal{C}_c(\XX),$ we have $\mathcal{P}_{s,t}f\in \mathcal{C}_c(\XX)$.

\vspace{.1cm}\item[(ii)] For each $f=|x|^p$, $p\geqslant 2$,  $\LG^\kappa_t f\in \mathcal{C}(\XX)$,  
$ \vert \LG^\kappa_t f(x)\vert \leqslant L_{b^\kappa,\sigma^\kappa,p}\left( 1+\vert x\vert^p\right)$, with $\mathcal{L}^\kappa_t$ given by~(\ref{gen}).
If $f\in \mathcal{C}_c^2(\XX)$ is supported on the ball {\it $\textsf{B}_\textsf{R}(0)$}, then {\it $\sup_{x\in \Scd}\vert \LG^\kappa_t f(x)\vert\leqslant L_{b^\kappa,\sigma^\kappa} \left( 3+ 8\textsf{R}^{\,2}\right)\Vert f\Vert_{\infty}.$}
\etm
\end{rem}

\bn[\textbf{\textit{The martingale problem}}\label{mrt_prb} \cite{Stroock79}]\label{mart_sol}\rm
Consider the  path space $\mathcal{W}_n:=C(\Ic,\XX)$  of  continuous maps $ \Ic\ni t\mapsto \gamma_t\in \XX$ equipped with the Borel $\mathfrak{S}$\,-\,algebra generated by open sets in the compact-open topology, as  for the Wiener space $\PS$. Identify $ \mathcal{W}_n \simeq  \Om$ via  $\Om\ni\om\mapsto \gamma_{\scaleobj{.8}{(\cdot)}}(\om)\in \XX$, s.t.~$(\mathcal{W}_n, \mathcal{B}(\mathcal{W}_n))\simeq(\Om, \mathcal{F})$; e.g., \cite{Arnold1}.  A probability measure  $\Ps_{t_0,x}\in \PP(\mathcal{W}_n)$ is a solution to the martingale problem for the operator $\mathcal{A}_t$ starting from $x\,{\in}\, \XX$,  $t_0\in \Ic$ if 
\btm[leftmargin = 0.75cm]

\vspace{-.0cm}\item[(i)] $\Ps_{t_0,x}\big(\{\om :\, \gamma_{t_0}(\om) = x\}\big) =1$. 

\vspace{.2cm}\item[(ii)] For any $f\in \mathcal{C}_c^{\infty}(\XX),$ the process $\Mt^f_t:=f\big(\gamma_t(\om)\big) \,{-}\, f\big(\gamma_{t_0}(\om)\big) \,{-}\int_{t_0}^t\mathcal{A}_u f\big(\gamma_u(\om)\big)du$,\\
is a $\Ps_{t_0,x}$\,-\,martingale w.r.t.~the filtration $(\mathcal{H}_t)_{t\geqslant t_0}$, $\mathcal{H}_t := \mathfrak{S}\{ \gamma_u: t_0\leqslant u\leqslant t\}$ on $\mathcal{W}_n$.
\etm

\noindent {The martingale problem is said to be well-posed if $\Ps_{t_0,x}$ exists uniquely for any $(t_0,x)\in \Ic\times\XX$.  } 
\en

\vspace{-.4cm}
\begin{rem}\rm
If the law $\Ms^\kappa_{t_0,x}$ of the solutions  (\ref{SDE_sol}) of (\ref{SDE}) are absolutely continuous w.r.t~ the Wiener measure $\p$, with Borel sets in, respectively, $\mathcal{B}(\mathcal{W}_n)$ and $\mathcal{B}(\Om)$, identified via the map $\Om\ni\om\mapsto \phi^\kappa_{t_0,\cccdot}(x,\om)\in \mathcal{W}_n$, then $\Ms^\kappa_{t_0,x}$ is a solution to the martingale problem for the generator  $\LG_t^\kappa$ (see Definition~\ref{weak_kol}) starting from $x\,{\in}\,\XX$ at $t_0\,{\in}\,\Ic$ (e.g.,~\cite{Malliavin97, Stroock79}). This holds, in particular, for $(b^\kappa,\sigma^\kappa)$ as in (\ref{reg_growth_cond}).
\end{rem}

\begin{prop}[\hspace{-0.03cm}\textbf{\textit{Martingale solutions and the forward Kolmogorov equation}} \cite{Stroock79,Figali}]\label{mart_sln}

\noindent  - \textsc{Existence}. 
\noindent The following are equivalent:
 \begin{itemize}[leftmargin = 0.8cm] 
 
 \item[(i)] Let  $(\kappa_t)_{t\in\Ic}$, $\kappa_t\,{\in}\, \PP(\XX)$ denote the solutions of (\ref{PDE}) and assume (\ref{well_pos_kolm}) holds. Then, there exists a measurable family of probability measures $(\Ps_{t_0,x}^\kappa)_{x\in \XX}$  such that $\Ps_{t_0,x}^\kappa \in \PP(\mathcal{W}_n)$ is a martingale solution of the SDE (\ref{SDE}) for $\mathcal{L}^\kappa_t$ starting at $t_0\,{\in}\, \Ic$ from $\kappa_{t_0}\text{-\,a.e.}$ $x\,{\in}\, \XX$, and 
 \begin{align}\label{k_lift}
 \int_{\XX}f(x)\kappa_t(dx) = \int_{\XX}\int_\Om f\big(\Solm^\kappa(x,\om)\big)\Ps_{t_0,x}^\kappa(d\om)\kappa_{t_0}(dx), \qquad f\in \mathcal{C}_c^{\infty}(\XX).
 \end{align}
 
 \item[(ii)] Assume that (\ref{k_growth})-(\ref{Lipsch}) hold, and  let $(\Ps_{t_0,x}^\kappa)_{x\in \XX}$  be a measurable family of probability measures on $\mathcal{W}_d$ such that $\Ps_{t_0,x}^\kappa\in \PP(\mathcal{W}_n)$ is a martingale solution of the SDE (\ref{SDE}) for the operator  $\mathcal{L}^\kappa_t$ starting $t_0\,{\in}\, \Ic$ from $\kappa_{t_0}\text{-\,a.a.}$ $x\in \XX$,  $\kappa_{t_0}\in \PP(\XX)$.  
 Then, $\kappa_t:=\Ps^{\hspace{0.01cm}\kappa}_{t_0}\circ\Solm^{\kappa,-1}\in \PP(\XX)$  satisfying (\ref{k_lift}) solves the forward Kolmogorov equation~(\ref{PDE}).
\end{itemize}
  \noindent - \textsc{Uniqueness}. The following are equivalent for $B$ a Borel set in $\XX$:
\begin{itemize}[leftmargin=0.8cm]
\item[(i)] Assuming that  (\ref{well_pos_kolm}) holds, time marginals of martingale solutions $(\Ps^{\kappa}_{s,x})_{x\in \XX}$ of the SDE~(\ref{SDE}) for the operator $\mathcal{L}^\kappa_t$ are unique for any $x\in B$, $B\in \mathcal{B}(\XX)$.
\item[(ii)] Assuming that (\ref{k_growth})-(\ref{Lipsch}) hold, finite non-negative measure-valued solutions of (\ref{PDE}) are weakly unique for any Borel probability measure $\kappa_{t_0}\in \PP(\XX)$ concentrated on $B\in \mathcal{B}(\XX)$.
\end{itemize}

\end{prop}

\begin{definition}[\textbf{\textit{Lebesgue a.e.~martingale solution}}\footnote{\,In \cite{Figali} Lebesque a.e.~martingale solutions are referred to as “Stochastic Lagrangian Flows”. We avoid this notion due to the potential confusion with the stochastic flows introduced in Definition \ref{stchf}.} ]\label{Lmartsol}\rm
Given a time-marginal probability  measure $\kappa_{t_0}\,{\in}\, \PP(\XX)$, $\kappa_{t_0}(dx) = \rho^\kappa_{t_0}(x)m_n(dx) $, $\rho^\kappa_{t_0}\in  L^1(\XX,m_n)\cap L^\infty(\XX;m_d)$, a measurable family of measures $(\Ps_{t_0,x}^\kappa)_{x\in \XX}$, $\Ps_{t_0,x}^\kappa\in \PP(\mathcal{W}_n)$, is a {\it Lebesgue a.e. \!$\kappa_{t_0}$-\,martingale solution} starting at $t_0\in \Ic$~if:
\btm[leftmargin=0.75cm]
\item[(i)] For $\kappa_{t_0}\textrm{-\,a.e.}$ $x,$ $\Ps_{t_0,x}^\kappa$ is a martingale solution of (\ref{SDE}) starting from $x\,{\in}\, \XX$,  $t_0\,{\in}\, \Ic$.

\vspace{.0cm}\item[(ii)] For any $t\in \Ic$ there exists   $\kappa_t(dx) = \rho^\kappa_t(x)m_n(dx)$, $\rho^\kappa_{t_0}\in L^1_+(\XX,m_n)\cap L^\infty(\XX;m_n)$ such that  
\begin{equation*}
\kappa_t= \Ps^{\kappa}_{t_0}\circ \Solm^{\kappa,-1}\ll m_n, \qquad \big(\Ps^{\kappa}_{t_0}\circ \Solm^{\kappa,-1}\big)(B) := \Ps^{\kappa}_{t_0}\big(\om\,{:}\,\, \Solm^\kappa(\ccdot,\om)\,{\in}\, B\big), \quad B\in \Bb(\XX),
\end{equation*}
\etm
where $\Ps^{\hspace{0.01cm}\kappa}_{t_0} \,{:=}\, \int_{\XX} \Ps_{t_0,x}^\kappa\, \kappa_{t_0}(dx)\in \PP(\mathcal{W}_n)$, and $\Solm^{\kappa,-1}$ is the inverse of $\Solm^{\kappa}$ solving~(\ref{SDE}). 
\end{definition}

Intuitively, Lebesgue a.e.~martingale solutions consist of those solutions to the martingale problem whose time marginals are absolutely continuous w.r.t.~the Lebesgue measure. Absolutely continuous measure solutions of the forward Kolmogorov equation coincide with time marginals of such martingale solutions. For $\kappa_{t_0}= \delta_x$,  $\Ps^{\kappa}_{t_0}\circ \Solm^{\kappa,-1}$ is the transition kernel~defined in~(\ref{P_phi}). 

\begin{prop}[\!\textbf{\textit{Existence and uniqueness of Lebesgue a.e.~martingale solution}~\cite{Figali}}]\label{lebae_exist} \mbox{}

\vspace{-.6cm}\begin{itemize}[leftmargin=0.8cm]

\item[(i)] Suppose that the forward Kolmogorov equation has solutions belonging to the convex subset $\mathcal{E}_{+}\subset L^1(\Ic\times\XX; m_n\otimes m)$ defined by\,\footnote{$f\in \mathcal{C}(\Ic; w^*-L^{\infty}(\XX; m_n)),$ means that $t\mapsto f(t)$ is $\text{weak}-^*$ continuous in $L^{\infty}(\XX; m_n).$ } 
\begin{align*}
\hspace{.8cm}\mathcal{E}_{+}(\Ic\times\XX) = \Big\{f\in L^{\infty}\big(\Ic; L_{+}^1(\XX;m_n)\big)\cap L^{\infty}\big(\Ic; L^{\infty}(\XX;m_n)\big)\!: \;f\in \mathcal{C}\big(\Ic; w^*-L^\infty(\XX;m_n)\big)\Big\}.
\end{align*}
Then, there exists  $(\Ps_{t_0,x}^\kappa)_{x\in \XX}$, s.t.~$\Ps_{t_0,x}^\kappa\in  \PP(\mathcal{W}_n)$, which is Lebesgue a.e.~$\kappa_{t_0}$-martingale solution of the SDE (\ref{SDE}) for $\kappa_{t_0}(dx) \,{=}\, \rho^\kappa_{t_0}(x)m_n(dx)$,  $\rho^\kappa_{t_0}\in L^1_{+}(\XX;m_n)\cap L^\infty(\XX;m_n)$.

\vspace{.2cm}\item[(ii)] If $(\tilde{\Ps}{\mathstrut}^\kappa_{t_0,x})_{x\in \XX}$ is Lebesgue a.e.~$\tilde{\kappa}_{t_0}$-martingale solution of (\ref{SDE}) for  $\tilde{\kappa}_{t_0}(dx) \,{=}\, \tilde{\rho}^\kappa_{t_0}(x)m_n(dx)$, $\tilde{\rho}^\kappa_{t_0}\,{\in}\, L_{+}^1(\XX;m_n)\cap L^\infty(\XX;m_n)$, 
then $\Ps_{t_0,x}^\kappa \,{=}\, \tilde{\Ps}{\mathstrut}^\kappa_{t_0,x}$ for $m_n \,\text{-\,a.e.}\; x\,{\in}\, \text{supp}(\tilde{\kappa}_{t_0})$.
\end{itemize}
\end{prop}

\begin{rem}\label{abs_con}\rm
The existence and uniqueness of Lebesgue a.e.~martingale solution depends on the growth and regularity of the coefficients, and regularity of the initial probability measure in~(\ref{SDE}). Consideration of Lebesgue a.e.~martingale solutions will be important in the reconstruction of the solutions of (\ref{F_Kol}a) in terms of the solutions of (\ref{F_Kol}b). Here, we focus on SDE's generating flows of $\mathcal{C}^{l}$-diffeomorphisms with  $t$-continuous coefficients  (see Glossary and Remark \ref{cf_remk})
\begin{alignat}{4}
b^\mm(t,\,\cdot\,) &\in \tilde{\mathcal{C}}^{l,\delta} (\XXt,\XXt), \quad &\sigma^\mm_k(t,\,\cdot\,)&\in \bar{\mathcal{C}}^{l,\delta}(\XXt,\XXt), &\qquad &1\leqslant k\leqslant m,&&\;l\geqslant 3, \;0<\delta\leqslant 1,  \;t\in\Ic, \notag\\
b^\nu(t,\,\cdot\,)  &\in \tilde{\mathcal{C}}^{l,\delta'} (\Xt,\Xt), \quad &\sigma^\nu_{k'}(t,\,\cdot\,)&\in \bar{\mathcal{C}}^{l,\delta'}(\Xt,\Xt), &\qquad &1\leqslant k'\leqslant m',& &\;l\geqslant 3, \;0<\delta'\leqslant 1, \;t\in \Ic,\notag
\end{alignat}
which are integrable over $\Ic$ and such that $\sigma^\mu(\sigma^\mu)^*$, $\sigma^\nu(\sigma^\nu)^*$ are uniformly elliptic.
  It is known from the Jacobi theorem (e.g., \cite{Amb99}) that a diffeomorphism is quasi-invariant w.r.t.~the Lebesgue measure on a finite-dimensional smooth manifold.  Thus, the flow of $\mathcal{C}^l$-\,diffeomorphisms induced by an SDE is a Lebesgue a.e.~martingale solution. In this setting, issues associated with  the potential lack or loss of the absolute continuity\footnote{\,The lack or loss of absolute continuity w.r.t.~to the Lebesgue measure may occur, for example,   due to finite-time explosion of solutions or, in the absence of uniform ellipticity or hypoellipticity in (\ref{F_Kol}).} are avoided, and considering $\varphi$-divergences (\ref{phi_dA}) between the time-marginal measures associated with the SDE dynamics~(\ref{SDE1}) and its approximation~(\ref{SDE2}) is well-posed. This follows from the fact that for probability measures with strictly positive (Lebesgue) densities $\rrto,\rho^\nu_{t_0}$, and for the dynamics generated by flows of $C^{l}$-diffeomorphisms, the solutions of the forward Kolmogorov equations~(\ref{F_Kol}) are absolutely continuous w.r.t.~the Lebesgue measure and they have strictly positive densities (see, e.g.,~\cite{Stroock79, Roc-Kry, Figali, Gmeas}). 
 \end{rem}

\section{Bounds on loss of information in path-based predictions}\label{Main_Bounds}

In this section, we consider bounds on discrepancies between laws of two different stochastic flows in terms of $\varphi$-divergences (\S\ref{phi_def}). These bounds can be subsequently utilised in the information inequalities derived in \S\ref{Info_ineq} in order to bound the error in estimating path-based observables in terms of the uncertainties in the Eulerian fields generating the underlying dynamics; see \S\ref{out_main} for the outline.  An example illustrating the utility of these results is presented in \S\ref{tests}.

First, in \S\ref{Inf_IP}, we derive an information bound via an appropriate reconstruction of  the generator of the  SDE (\ref{SDE1}) in terms of the generator of  (\ref{SDE2}). 
For a sufficiently non-degenerate SDE dynamics, this  approach allows to express information bounds in terms of differences between coefficients of the two SDEs and it provides an analytically tractable connection between the Eulerian (field-based) error and the uncertainty in Lagrangian (path-based) predictions. In \S\ref{ftdr_sec}, we  derive  a   bound on uncertainty in Lagrangian predictions based on the difference between so-called {\it finite-time divergence rate} ($\varphi$-FTDR) fields~\cite{MBUda18} which utilise a recently developed  framework for quantifying expansion rates in arbitrary stochastic flows. Importantly, this bound is not restricted to SDE/ODE dynamics and  it can be exploited  within a computational framework, in both the stochastic and deterministic settings, to mitigate the error in Lagrangian predictions by tuning the $\varphi$-FTDR fields in simplified models to those  generated by the original dynamics.  Finally, in \S\ref{Path-space}, we extend  the results of \S\ref{Inf_IP}--\ref{ftdr_sec} to the path space by means of  a projection of certain $\varphi$-admissible path space measures associated with the original dynamics onto path space measures associated with the approximation.

 \subsection{Bound on information loss in Lagrangian predictions via generator reconstruction}\label{Inf_IP}
Here, we consider the relationship between solutions of forward Kolmogorov equations for the SDEs (\ref{SDE1}) and  (\ref{SDE2}),  following an approach  recently developed in \cite{Rockner16} in the context of the KL-divergence and the total variation distance. The main idea is to represent one of the forward Kolmogorov equations in terms of the other one. This type of `reconstruction' is standard when the SDEs have the same diffusion coefficients; for different diffusion coefficients, the reconstructed Kolmogorov equation may be not as regular as the original one even in  the uniformly  elliptic case. However, we shall demonstrate in \S\ref{Recon_sol} that, under some non-degeneracy assumptions on the coefficients in (\ref{SDE1}) and (\ref{SDE2}), the reconstructed Kolmogorov equation generates Lebesgue a.e.~martingale~solutions. This, in turn, allows to derive a bound on the $\varphi$-divergence between time-marginal measures induced by the original dynamics and its approximation; such a bound can then be used in conjunction with (\ref{Csiszar}) to mitigate errors in path-based observables.  
 
 In order to simplify derivations, we first consider the case when the original and approximate dynamics evolve on the same domain; i.e., $\Xt=\XXt$. The case of $\Xt\subset \XXt$ is discussed in \S\ref{MM-M}.   

\subsubsection{\bf Generator reconstruction when  $\Xt=\XXt$}\label{Xt=XXt}
 In this case $\mmt = \mu_t$ in (\ref{margs0}), $\mu_t\in \PP(\Xt)$, and we set $\brr^\mm \equiv \brr^\mu$, $\sigma^\mm \equiv \sigma^\mu$ to highlight this fact. Suppose that the coefficients $(\brr^\mu,\sigma^\mu)$ of  (\ref{SDE1}), and $(\brr^\nu,\sigma^\nu)$ of  (\ref{SDE2}) are such that $(\Ms_{t_0,x}^\mu)_{x\in\Xt}$ is the Lebesgue a.e.~$\mu_{t_0}$-martingale solution for the generator $\mathcal{L}^\mu_t$, and  $(\Ns_{t_0,x}^\nu)_{x\in\Xt}$ is the Lebesgue a.e.~$\nu_{t_0}$-martingale solution for the generator  $\mathcal{L}^\nu_t$ (cf.~Definition~\ref{Lmartsol} and Remark~\ref{abs_con}). Then, the families of probability densities $(\rho_{t}^{\mu})_{t\in \Ic}$ and $(\rho^{\nu}_{t})_{t\in \Ic}$ w.r.t.~the Lebesgue measure on $\Xt$ exist  \cite{Figali, Roc-Kry, Stroock79}  and satisfy (in the weak sense on $\mathcal{E}_+$)  
\begin{align}
\partial_t\rho^{\mu}_{t} = \LG^{\mu*}_{t}\rho^{\mu}_{t} = \textstyle \frac{1}{2}\partial^2_{x_ix_j}(a^{\mu}_{ij}\rho^{\mu}_{t}) - \partial_{x_i}({b}_i^{\mu}\rho^{\mu}_{t}), \qquad \rho_{t}^\mu\in L^1_+(\Xt,dx)\cap L^\infty(\Xt,dx),\label{fke1}\\[.2cm]
\partial_t\rho^{\nu}_{t} = \LG^{\nu*}_{t}\rho^{\nu}_{t} = \textstyle\frac{1}{2}\partial^2_{x_ix_j}(a^{\nu}_{ij}\rho^{\nu}_{t}) - \partial_{x_i}({b}_i^{\nu}\rho^{\nu}_{t}), \qquad \rho_{t}^\nu\in L^1_+(\Xt,dx)\cap L^\infty(\Xt,dx),\label{fke2}
\end{align} 
with  ${b}_i^{\scriptscriptstyle (\ccdot)} \,{=}\, \brr_i^{\scriptscriptstyle (\ccdot)}\,{+} \frac{1}{2} \sigma_{jk}^{\scriptscriptstyle (\ccdot)}\partial_{x_j} \sigma_{ik}^{\scriptscriptstyle (\ccdot)}$, $a_{ij}^{\scriptscriptstyle (\ccdot)}\,{=}\,\sigma_{ik}^{\scriptscriptstyle (\ccdot)}\sigma_{jk}^{\scriptscriptstyle (\ccdot)}$, and summation implied over repeated indices.

\begin{lem}[\textbf{\textit{Reconstructed Kolmogorov equation}}]\label{THT_def}
  Given the solutions $(\rho_t^\mu)_{t\in \Ic}$, $(\rho_t^\nu)_{t\in \Ic}$ of the forward Kolmogorov equations (\ref{fke1}) and (\ref{fke2}), we have the following reconstructed equation
\begin{align}\label{RecEQ}
\partial_t\Sden = \LG^{\nu*}_{t}\Sden -\nabla_x\,{\cdot} \,\big(\varTheta_{\mu\nu}\Sden\big),
\end{align}
with solutoins and derivatives understood in the distributional sense, and  where the reconstructed field $\betaLk$ is given by 
\begin{equation}\label{Per_vF}
\varTheta_{\mu\nu}(t,x) := \textstyle\frac{1}{2}\big(a^{\nu}(t,x) - a^{\mu}(t,x)\big)\nabla_x\log\Sden(x)- \big(h^{\nu}(t,x) - h^{\mu}(t,x)\big),
\end{equation}
with 
\begin{equation*}
 h^{\mu}_i(t,x) := \textstyle{b}^{\mu}_i(t,x) - \frac{1}{2}\partial_{x_j}a^{\mu}_{ij}(t,x), \qquad h_i^{\nu} (t,x):= {b}^{\nu}_i(t,x)-\frac{1}{2}\partial_{x_j}a^{\nu}_{ij}(t,x). \notag
\end{equation*}
\end{lem}

\noindent {\it Proof.} This is derived directly as follows (with derivatives understood in the distributional sense)
\begin{align*}
\hspace{.8cm}\partial_t\Sden - \LG^{\nu*}_{t}\Sden &= \LG^{\mu*}_{t}\Sden-\LG^{\nu*}_{t}\Sden = \textstyle \partial_{x_i}\Big(\frac{1}{2}\partial_{x_j}(a^{\mu}_{ij}\Sden) - \frac{1}{2}\partial_{x_j}(a^{\nu}_{ij}\Sden) + {b}_i^{\nu}\Sden-{b}_i^{\mu}\Sden \Big)\\
&\textstyle =-\partial_{x_i}\left(\frac{1}{2}\partial_{x_j}\big((a^{\nu}_{ij} - a^{\mu}_{ij})\rho^{\mu}_t\big)+({b}_i^{\mu} - {b}_i^{\nu})\rho^{\mu}_t\right)\\
 &\textstyle = -\partial_{x_i}\left(\frac{1}{2}(a^{\nu}_{ij}-a^{\mu}_{ij})\partial_{x_j}\Sden -\big(({b}_i^{\nu}-\frac{1}{2}\partial_{x_j}a_{ij}^{\nu})-({b}_i^{\mu}- \frac{1}{2}\partial_{x_j}a^{\mu}_{ij})\big)\Sden\right)\\ &=
-\nabla_x \cdot(\varTheta_{\mu\nu}\Sden). \hspace{9.5cm}\qed
\end{align*}

\medskip
\noindent The main result of this section, which implicitly relies on the weak solvability of (\ref{RecEQ}) is as follows:

\vspace{-0.2cm}\begin{theorem}[\textbf{\textit{Information bound for time-marginal probability measures of SDEs}}]\label{Info_ineq2}\mbox{}

\noindent Consider the dynamics induced by the SDE  (\ref{SDE1}) and its approximation (\ref{SDE2}) on $\Xt$, and their respective time-marginal measures $\mu_t$, $\nu_t\in \PP(\Xt)$.  Assume that the following conditions hold:
\begin{itemize}[leftmargin=0.8cm]
\item[(i)] The coefficients in (\ref{SDE1}), (\ref{SDE2}) are such that $\brr^{\,\mu}, \brr^{\,\nu}\in \mathcal{C}\big(\Ic; \tilde{\mathcal{C}}^{3,\delta}(\Xt; \Xt)\big)$, and the columns of $\sigma^\mu$, $\sigma^\nu$,  are  $\sigma_{k}^{\,\mu}, \sigma_{k'}^{\,\nu}\in {\mathcal{C}}\big(\Ic; \bar{\mathcal{C}}^{4,\delta}(\Xt;  \Xt)\big)$. Moreover, the right inverses  $\tilde\sigma^{\mu,-1}, \tilde\sigma^{\nu,-1}$ of $\sigma^\mu, \sigma^\nu$ are uniformly bounded on $\Xt$ and strictly positive (i.e., (\ref{SDE1}), (\ref{SDE2}) have uniform ellipticity), and $b^\mu-b^\nu, a^\mu-a^\nu, \partial_{x_j}a_{ij}^\mu, \partial_{x_j}a_{ij}^\nu\in L^\infty(\Ic\times\Xt)$\footnote{\,The conditions on $b^\mu,b^\nu,a^\mu,a^\nu$ can be considerably weakened but we defer such generalizations to future~work.}.

\vspace{0.2cm}\item[(ii)] The probability measures on the initial conditions in (\ref{SDE1}) and (\ref{SDE2})  have all moments finite and are such that\,\footnote{\,Here and below we set $dx\equiv m_d(dx)$ with $m_d$ the Lebesgue measure on $\Xt$ to simplify notation.} 
 $$\mu_{t_0}(dx) = \nu_{t_0}(dx) = \rho_{t_0}(x)dx, \qquad \rho_{t_0}\in L_{+}^1(\Xt,dx)\cap L^\infty(\Xt,dx).$$ 
\end{itemize}
Then, for $\varphi\in \mathcal{C}^2(\Rp)$ a strictly convex function satisfying  (\ref{Normality}) the following holds  
\begin{align}\label{InformaBB}
\notag \Df(\mu_t\|\nu_t) &= \int_{\Xt}\varphi\left(\eta_t(x)\right)\rho_t^{\nu}(x)dx\\[.1cm] 
& \hspace{.0cm} \leqslant \frac{1}{2}\int_\Ic\int_{\Xt}\big\vert (\tilde\sigma^{\nu,-1}\varTheta_{\mu\nu})(s,x)\big\vert^2 \varphi^{\prime\prime}(\eta_s(x))\eta_s^2(x)\rho_s^{\nu}(x)dxds, \quad\, t>t_0,
\end{align} 
where the $\varphi$-divergence $\Df$ is as defined in (\ref{phi_dA}),  $\betaLk$ is as defined in (\ref{Per_vF}), and  
$\mu_t, \nu_t\in \PP(\Xt)$ are such that $\mu_t(dx) = \rho^\mu_t(x)dx$, $\nu_t(dx) = \rho^\nu_t(x)dx$, with  $\eta_t(x) := \rho^\mu_t(x)/\rho_t^{\nu}(x)<\infty$, $t\in\Ic$, where $\rho^\mu_t, \rho^\nu_t>0$ 
solve the forward Kolmogorov equations (\ref{fke1}) and  (\ref{fke2}).
\end{theorem}
\vspace{-.1cm}\noindent {\it Proof}. See \S\ref{5.3proof} and comments below; preparatory results are obtained in \S\ref{Recon_sol}.

\begin{cor}
Setting  $\varphi(u) = u\log u-u+1, \; u>0$, in (\ref{InformaBB}) leads to a simplified  bound on the lack of information in $\nu_t$ relative to $\mu_t$ in terms of the KL-divergence; namely  
\begin{align}
\DK(\mu_t||\nu_t) &\leqslant \frac{1}{2}\int_\Ic\int_{\Xt}\big\vert \tilde\sigma^{\nu,-1}\varTheta_{\mu\nu}(s,x)\big\vert^2\rho^\mu_s(x)dxds, \qquad t>t_0. 
\end{align}

\end{cor}

\begin{rem}\rm The following observations are worth pointing out:
\begin{itemize}[leftmargin=.7cm]

\item[(i)] Evidently, the vector field $ \varTheta_{\mu\nu}$ in the information bound in Theorem~\ref{Info_ineq2} gives a connection between Eulerian (field-based) model error and uncertainty in the Lagrangian (path-based) predictions. It is analytically tractable (see \S\ref{tests} for an example) as the bound is based on the coefficients, $(\brr^\mu,\sigma^\mu)$ and $(\brr^\nu,\sigma^\nu)$, of the respective SDEs;  this is crucial for our purpose in combination with the path-based version of the bound (\ref{Csiszar}), since it is generally not possible to derive explicit forms of $\mu_t$ and $\nu_t$, even for SDE's with simple coefficients.

\item[(ii)] The following instance  of the vector field $\betaLk $ is well-studied in theory and applications (e.g., \cite{Roc-Kry, Yor}).  If the diffusion coefficients in (\ref{SDE1}),  (\ref{SDE2}) are such that  $\sigma^{\nu} = \sigma^{\mu} = \sigma>0,$ we have 
\begin{align*}
\betaLk(t,x) =  b^{\nu}(t,x)-b^{\mu}(t,x),
\end{align*}
which is simply the difference between the drift terms of the original and approximate dynamics. Even in this case, minimisation of the loss of information in (\ref{InformaBB}) between $\mu_t$ and $\nu_t$ involves the $L^2(\Xt,\mu_t)$ norm of difference of the two  (Eulerian) fields. 
\end{itemize}
\end{rem}

\begin{rem} \rm The following comments are in order:
\begin{itemize}[leftmargin=0.7cm]

\item[(i)]The uniform ellipticity assumptions in Theorem \ref{Info_ineq2} could be relaxed to allow for the hypoelliptic case but an even more lengthy proof would involve dealing with Malliavin covariance, Malliavin integration by parts, and  the generalised It\^o isometry.   The issue of approximating deterministic dynamics is much more subtle in this framework, since for $\sigma^\mu=0$ the associated transition evolutions $(\mathcal{P}^\mu_{s,t})_{t\geqslant s}$ will not, in general, have the smoothing property (e.g.,~\cite{Daprato08}). The smoothing property, and $\sigma^\nu>0$, are  necessary for  deriving the upper bound on $\Df(\mu_t\|\nu_t)$ in~(\ref{InformaBB}); see, in particular, the proof of Proposition \ref{log_gr} and Corollary \ref{sqrint} below. Model tuning in the deterministic case can be considered in this framework through viscosity solutions of (\ref{fke1}), e.g., \cite{dolcetta95}.  We postpone such generalisations to a separate publication.  

\item[(ii)] It can be shown in a way analogous to \cite[Proposition~4.4]{MBUda18} that, under the assumptions of Theorem \ref{Info_ineq2}, $\Df(\mu_t\|\nu_t)<\infty$ and thus $\varphi(\eta_t)\in L^1(\Xt,\nu_t)$. However, the bound in (\ref{InformaBB}) is non-uniform in $T$. This fact does not prevent one from minimising it in terms of the coefficients $(\brr^\nu,\sigma^\nu)$ of the approximating SDE (\ref{SDE2}).  

\end{itemize}
\end{rem}

The proof of Theorem \ref{Info_ineq2} requires some preparation and is postponed to \S\ref{5.3proof}.  The main technical issue, dealt with in \S\ref{Recon_sol}, which is implicitly required in the main proof concerns establishing the weak solvability of (\ref{RecEQ}). Sufficient regularity and growth conditions of the vector field $\varTheta_{\mu\nu}$  (\ref{Per_vF}) for  weak solvability of (\ref{RecEQ}) with the operator $\LG^{\nu*}_t-\nabla (\varTheta_{\mu\nu} \,\ccdot\,)$  are  not immediately obvious from the properties of  the coefficients of the associated SDEs,  even if $\LG^{\nu}_t$ is non-degenerate. This is due to the presence of  the logarithmic gradient $\nabla_x \log\rho^{\mu}_t$ in the vector~field~$\varTheta_{\mu\nu}$. However, we show that whenever $\LG_t^{\nu}$ is a one-point generator of a stochastic flow of diffeomorphisms, one can construct a Lebesgue a.e.~martingale solution (Definition \ref{Lmartsol}) from the reconstructed operator $\LG^{\nu}_t+ \varTheta_{\mu\nu}\nabla$; this is equivalent to  establishing the weak solvability of the Cauchy problem associated with  $\LG^\nu_t+\varTheta_{\mu\nu}\nabla$ or its $L^2(\Xt;dx)$ adjoint $\LG^{\nu*}_t-\nabla\,{\cdot}\,(\varTheta_{\mu\nu}\,\,\cdot\,\,).$ The main steps which are necessary for solvability of~(\ref{RecEQ}) rely on the fact that, if assumptions of Theorem \ref{Info_ineq2} are satisfied, we have 
 \begin{equation*}
  \int_{\Ic}\int_{\Xt}\vert (\tilde\sigma^{\nu,-1}\varTheta_{\mu\nu})(s,x)\vert^n \nu_s(dx)ds  <\infty,  \qquad  n\in \mathbb{N}_1,\;\;t>t_0.
 \end{equation*}
 Then, the Girsanov theorem is used  to construct the transition evolution (cf.~Definition \ref{trans_evo}) for an It\^o process with coefficients $\big({b}^{\nu}+\betaLk, \sigma^\nu\big)$,  
which implies solvability of the reconstructed backward Kolmogorov equation for $t\in \Ic$
\begin{equation*}
\partial_t u(t,x) + \mathcal{L}^{\nu}_tu(t,x)+ \varTheta_{\mu\nu}\nabla_x u(t,x)=0,\qquad u(t_0+T,x)  = f(x)\in \mathcal{C}^2_\infty(\Xt).
\end{equation*}
The remainder of the proof is relatively straightforward and is given in \S\ref{5.3proof}.

\bigskip
 \subsubsection{Solutions of the reconstructed Kolmogorov equation}\label{Recon_sol}
Here, we investigate the weak solvability of the reconstructed equation (\ref{RecEQ}) which is needed in the proof of Theorem \ref{Info_ineq2}  discussed in \S\ref{5.3proof}. The main argument relies on Girsanov's theorem but the derivation requires a few preparatory results which are discussed first.

\vspace{-0.cm}\begin{prop}\label{Lim_procedure}\label{log_gr} 
Assume that the conditions (i)-(ii) of Theorem \ref{Info_ineq2} are satisfied.  Then
 \begin{align*}
\notag &\big\vert \nabla_x\log\rho^\mu_{t}(x)\big\vert^2 \leqslant\frac{\mathfrak{C}_{\sigma^\mu}}{(t-t_0)^2}\sup_{y\in\Xt}\E\left[\int_{t_0}^t\Vert D_x\phi^\mu_{t_0,\xi}(x)\Vert_\textsc{hs}^2 \,ds \Big|\phi^\mu_{t_0,t}(x)=y\right]<\infty, \quad t>t_0, 
\end{align*}  
where  $0<\mathfrak{C}_{\sigma^\mu}<\infty$ is independent of time,  $\rho_t^\mu$ weakly solves  (\ref{fke1}), and $D_x\phi^\mu_{t_0,t}(x,\om)$ is the derivative flow of $\phi^\mu_{t_0,t}(x,\om)$ which is itself a flow of $\mathcal{C}^2$-diffeomorphisms associated with the SDE~(\ref{SDE1}). Moreover, $\nabla_x\log\rho^\mu_{t}(x)$ has bounded first and second derivatives. 
\end{prop}

\noindent {\it Proof.} See Appendix \ref{Uniell_prof2}; the proof relies on the Bismut-Elworthy-Li formula.

\bigskip
\begin{cor}\label{sqrint}
Proposition \ref{Lim_procedure} implies that under appropriate assumptions on the coefficients of (\ref{SDE1}) and (\ref{SDE2}), as in Theorem \ref{Info_ineq2}, the reconstructed vector field $\betaLk$~(\ref{Per_vF}) is such that 
\begin{align}
 \int_{\Ic}\int_{\Xt}\vert (\tilde{\sigma}^{\nu,-1}\betaLk)(s,x)\vert^n \nu_s(dx)ds<\infty, \quad n\in \mathbb{N}_1,\;\;t>t_0. \label{sq_int}
\end{align}
\end{cor}

\medskip
\noindent{\it Proof.}
This fact, which will be relevant in Proposition \ref{Summary_prop},  can be established by recalling that 
\begin{align*}
\betaLk = \frat\big( a_{ij}^\nu -a_{ij}^\mu\big)\nabla_x\log\rho^\mu -\big( {b}_i^\nu -\partial_{x_j}a_{ij}^\nu -{b}_i^{\mu} +\partial_{x_j}a^{\mu}_{ij}\big),
\end{align*}
 where we skip the explicit dependence on $(t,x)$. Given the growth and regularity of the coefficients in (\ref{SDE2}) assumed in Theorem \ref{Info_ineq2}, there exists a Lebesgue a.e.~martingale solution starting at $t_0\in \Ic$, $x\in \Xt$, for the generator $\mathcal{L}^\nu_t$ (cf.~Definition~\ref{Lmartsol}, Proposition~\ref{lebae_exist}, and Remark~\ref{abs_con}). Thus, there exists a time-marginal probability measure $\nu_t$ solving (\ref{PDE}) on $\Ic\times\Xt$, which is such that  $\nu_t(dx)= \rho^\nu_t(x)dx,\, \rho^\nu_t>0$ (Propositions \ref{mart_sln}, \ref{lebae_exist} and Remark \ref{cf_remk}). Application of the  Cauchy-Schwarz inequality to the left-hand side of (\ref{sq_int})  gives the desired result, since all that is needed is  the bound in Proposition~\ref{log_gr}, and  existence of  the moments $\int_\Xt |x|^n\rho^\nu_t(x)dx$ for $t\in \Ic$; these follow from the fact that 
$$\int_\Ic\int_\Xt \big\vert (\tilde{\sigma}^{\nu,-1}\betaLk)(s,y)\big\vert^n\nu_s(dy)ds=\int_\Ic\int_\Xt\!\E\left[ \big\vert (\tilde{\sigma}^{\nu,-1}\betaLk)(s, \phi^\nu_{t_0,s}(x))\big\vert^n ds\right]\nu_{t_0}(dx);$$
the existence of the right-hand side is satisfied due  to the growth conditions imposed on $(\brr^\mu,\sigma^\mu)$ and $(\brr^\nu,\sigma^\nu)$ which are present in $\betaLk$  (see also Remark \ref{cf_remk} and Proposition \ref{Lim_procedure})).

The above construction yields the following version of Girsanov's transformation:

\vspace{-0.2cm}\begin{prop}\label{Summary_prop}
Assume that the conditions (i)-(ii) of Theorem \ref{Info_ineq2} hold and consider 
\begin{align}
d X^\nu_t &= {b}^\nu\big(t,X_t^\nu\big)dt+\sigma^\nu\big(t,X_t^\nu\big)dW_{t-t_0}, \;\;\; {X}^\nu_{t_0}\sim \nu_{t_0},\label{mrk} \\[.2cm]
d \tilde X^\nu_t &= \tilde{b}^\nu\big(t,\Yt^\nu\big)dt+\sigma^\nu\big(t,\Yt^\nu\big)d\tilde W_{t-t_0}, \;\;\; \tilde{X}^\nu_{t_0}\sim \nu_{t_0}, \label{Rec_Markov}
\end{align}
where (\ref{mrk}) is an It\^o diffusion,  and (\ref{Rec_Markov}) is an It\^o process with $\tilde{b}^\nu = {b}^{\nu}+\betaLk$ and $\betaLk$  defined in (\ref{Per_vF}). Then, the transition evolution $(\tilde{\mathcal{P}}^{\nu}_{t_0,t})_{t\in \Ic}$ associated with (\ref{Rec_Markov}) can be represented as  
\begin{align*}
\tilde{\mathcal{P}}^\nu_{t_0,t}f(x) = \E\big[f\big({X}^\nu_{t}\big)g(t_0,t,x)\big], \;\;\; f\in\mathcal{C}_\infty(\Xt),
\end{align*}
where the martingale $\big(g(t_0,t,x,\omega)\big)_{t\geqslant t_0}$ is given by 
\begin{align*}
g(t_0,t,x,\omega) = \exp\left\{-\frac{1}{2}\int_{t_0}^t\left\vert (\tilde{\sigma}^{\nu,-1}\betaLk)\big(\xi, {X}^\nu_{\xi}(x,\om)\big)\right\vert^2\!d\xi +\int_{t_0}^t\Big\langle (\tilde{\sigma}^{\nu,-1}\betaLk)\big(\xi, {X}^\nu_{\xi}(x,\om)\big), dW_{\xi}\Big\rangle\right\}.
\end{align*}

\end{prop}

\noindent {\it Proof.} This is standard (e.g., \citep{Yor, Daprato08, Stroock79}) given the results outlined earlier in this section.

\begin{cor}\label{summ_cor}
Proposition \ref{Summary_prop}  yields the weak solvability of the backward and forward Kolmogorov equations (e.g., \cite[Theorem 10.8]{Daprato08} with some modifications on the growth conditions)  
\begin{alignat}{3}
\partial_t u(t,x) &+ \mathcal{L}^{\nu}_tu(t,x)+ \varTheta_{\mu\nu}\nabla_x u(t,x)=0,&\qquad &u(t_0+T,x)  = f(x)\in \mathcal{C}^2_\infty(\Xt), \label{weak_sol_Xt}\\
\partial_t\rho_t^{\mu} &= \LG^{\nu*}_{t}\rho_t^{\mu}-\nabla_x\cdot(\betaLk\rho_t^{\mu}), &\qquad &\rho_{t_0}^{\mu}(x) \in L_+^1(\Xt;dx)\cap L^\infty(\Xt;dx).
\end{alignat}
\end{cor}

\subsubsection{\bf Proof of Theorem \ref{Info_ineq2}}\label{5.3proof}
With the results derived in \S\ref{Recon_sol}, we are ready to prove our result on information bound between Lebesgue a.e.~flows of solutions of two stochastic differential equations on $\Xt.$ First, we give three preparatory lemmas which will facilitate the proof.

\vspace{-0.2cm}\begin{lem}\label{Ide_Lem}
Let $\LG_t^*$ be the $L^2$ dual of the operator $\LG_t$ with the coefficients $(b,a)$, $a = \sigma \sigma^*$. Then, for  $f,g\in \mathcal{C}^{2}(\Xt)$ and $\varphi\in \mathcal{C}^2(\Rp),$ 
\begin{align*}
\mathcal{L}_t^{*}\varphi(f) &= \varphi^{\prime}(f)\LG_t^*f+ \frat\varphi^{\prime\prime}(f)\langle a\nabla f, \nabla f\rangle + \left(f\varphi^{\prime}(f) - \varphi(f)\right)\nabla\cdot (b-\frat\nabla a),\\
\LG_t^*(fg)& = f\LG_t^*g+g\LG_t^*f + \langle a \nabla f, \nabla g\rangle + fg \nabla\cdot(b-\frat\nabla a).
\end{align*}
\end{lem}
\noindent {\it Proof.} See Appendix \ref{Ide_Lem_ap}; the proof follows from the chain and product rules for differentiation.

\vspace{-0.1cm}\begin{lem}\label{Pre_lem2}
Let $\varphi\in \mathcal{C}^2(\Rp)$ and assume the   (i)-(ii) of conditoins Theorem \ref{Info_ineq2} hold. Then 
\begin{align}\label{Pre_eq}
\partial_t(\varphi(\eta_t)\rho_t^\nu) = \LG_t^{\nu*}(\varphi(\eta_t)\rho_t^\nu)-\frat\rho_t^\nu\varphi^{\prime\prime}(\eta_t)\langle a^\nu\nabla\eta_t, \nabla\eta_t\rangle -\varphi^{\prime}(\eta_t)\nabla\cdot(\betaLk\rho_t^\mu), 
\end{align} 
where $0<\eta_t := \rho^\mu_t/\rho^\nu_t<\infty$. 
\end{lem}

\vspace{-0.2cm}
\noindent {\it Proof.} See Appendix \ref{Pre_lem2_ap}; the proof follows by direct application of Lemma \ref{Pre_lem2}, and it utilises weak solvability of the forward Kolmogorov equation in Corollary \ref{summ_cor}, as well as properties of solutions for the forward Kolmogorov equation (analogous to \cite[Lemma 2.4]{Rockner16}).

\vspace{-0.cm}\begin{lem}\label{Integrand}
 Assume that the conditions (i)-(ii) of Theorem \ref{Info_ineq2} hold. Then, for $f\in \mathcal{C}_c^{\infty}(\Xt)$ and any compact interval $[\tau, \,t]\subseteq \Ic$, we have 
\begin{align}\label{GW19}
\notag &\int_{\Xt}\varphi(\eta_t(x))f(x){\nu}_t(dx) + \frac{1}{2}\int_{\tau}^{t}\int_{\Xt}\big\vert( \sigma^{\nu})^*(s,x)\nabla_x\eta_s(x)\big\vert^2\varphi^{\prime\prime}(\eta_s(x))f(x)\nu_s(dx)ds\\[.1cm]
\notag &\hspace{4cm}= \int_{\Xt}\varphi(\eta_\tau(x))f(x)\nu_\tau(dx)+\int_{\tau}^t\int_{\Xt}\varphi(\eta_s(x))\LG_t^\nu f(x)\nu_s(dx)ds\\[.1cm]
\notag &\hspace{4.6cm}+ \int_{\tau}^t \!\!\int_{\Xt}\!\bigg[ \big\langle \betaLk(s,x)\eta_s(x), \nabla_x\eta_s(x)\big\rangle\varphi^{\prime\prime}(\eta_s(x))f(x)\\
&\hspace{7.cm} + \big\langle \betaLk(s,x), \nabla f(x)\big\rangle\varphi^{\prime}(\eta_s(x))\eta_s(x)\bigg] \nu_s(dx)ds.
\end{align}

\vspace{-0.2cm}
\noindent{\it Proof.} \rm See Appendix \ref{Integrand_ap}; this follows by multiplying  by  integration of (\ref{Pre_eq}) against $f\in \mathcal{C}_c^{\infty}(\Xt)$, and the standard Newton-Leibnitz formula. 
\end{lem}

\noindent {\bf Proof of Theorem \ref{Info_ineq2}.}
First, re-write the term $\langle \betaLk \,\eta_s, \nabla\eta_s\rangle$ of the integrand in the equality~(\ref{GW19}) of Lemma (\ref{Integrand}) using the Young inequality to obtain 
\begin{align}\label{GW23}
 \big\langle \tilde{\sigma}^{\nu,-1}\betaLk\,\eta_s, (\sigma^{\nu})^*\nabla\eta_s\big\rangle \leqslant \frat\big\vert \tilde{\sigma}^{\nu,-1}\betaLk\big\vert^2\eta_s^2+\frat\big\vert (\sigma^\nu)^*\nabla\eta_s\big\vert^2.
\end{align}
 where $\tilde\sigma^{\nu,-1}$ is the (uniformly bounded) right inverse of $\sigma^\nu$. Combining (\ref{GW23}) with  (\ref{GW19})  leads~to  
\begin{align}\label{GW24}
\notag \int_{\Xt}\varphi(\eta_t(x))f(x){\nu}_t(dx) 
 & \leqslant  \int_{\Xt}\varphi(\eta_\tau(x))f(x)\nu_\tau(dx)+\int_{\tau}^t\int_{\Xt}\varphi(\eta_s(x))\LG_t^\nu f(x)\nu_s(dx)ds\\
\notag &\hspace{.5cm}+ \int_{\tau}^t \int_{\Xt}\langle \betaLk(s,x), \nabla f(x)\rangle\varphi^{\prime}(\eta_s(x))\eta_s(x)\nu_s(dx)ds\\
&\hspace{1cm}+ \frac{1}{2}\int_{\tau}^t\int_{\Xt}\vert (\tilde{\sigma}^{\nu,-1}\betaLk)(s,x)\vert^2\varphi^{\prime\prime}(\eta_s(x))\eta_s^2(x)f(x)\nu_s(dx)ds.
\end{align}
Next, we show that for $f\in \mathcal{C}_c^{\infty}(\Xt)$ with $f\geqslant 0,$ 
\begin{align}\label{GW25}
\lim_{\tau\rightarrow t_0}\int_{\Xt}\varphi(\eta_\tau(x))f(x)\nu_\tau(dx) =0.
\end{align}
From the strict convexity of $\varphi\in \mathcal{C}^2(\R^+)$ and the normality condition (\ref{Normality}), we know that $\varphi$ is locally bounded Lipschitz continuous; in particular, given that $\varphi(1)=0$, $0<\eta_t<\infty$ (cf. \cite[Lemma~2.4]{Rockner16} and Lemma~\ref{Pre_lem2}), we have  $C_{\varphi}>0$ such that 
 \begin{align*}
\varphi(\eta) =  \varphi(\eta) -\varphi(1) \leqslant C_{\varphi}\vert \eta -1\vert, 
\end{align*}

\vspace{-.3cm}
\noindent and, consequently,  
\begin{align*}
\int_{\Xt}\varphi(\eta_\tau(x))f(x)\nu_\tau(dx)\leqslant\int_{\Xt}\vert \varphi(\eta_\tau(x))-\varphi(1)\vert f(x)\rho_\tau^\nu(x)dx\leqslant C_\varphi\int_{\Xt}\vert \rho_\tau^\mu(x)-\rho_\tau^\nu(x)\vert f(x)dx.
\end{align*}
Since $\mu_t(dx) = \rho^\mu_t(x)dx$, $\nu_t(dx) = \rho_t^\nu(x)dx$ for $t\in \Ic$, $\rho^\mu_{t_0}=\rho^\nu_{t_0}>0$,  the smoothness of the coefficients of (\ref{SDE1}), (\ref{SDE2})  ensures that $\rho_t^\nu, \rho_t^\mu$ are strictly positive and regular, and  by \cite[Lemma~2.1]{Rockner16} 
\begin{align}\label{Rockner_limit}
\lim_{\tau \rightarrow t_0} \int_{\Xt}\vert \rho_\tau^\mu(x)-\rho^\nu_\tau(x)\vert f(x)dx =0, \qquad f\in \mathcal{C}_c^{\infty}(\Xt), \;f\geqslant 0.
\end{align}
The limit (\ref{GW25})  follows from (\ref{Rockner_limit}). 

The second term on the right-hand side of (\ref{GW24}) is bounded as follows. Consider a cut-off function $\chi\in \mathcal{C}_c^\infty(\Xt)$ such that 
$$\chi(x) = 1 \;\; \textrm{for} \;\;|x|< 1, \quad \textrm{and} \quad \chi(x) = 0 \;\; \textrm{for} \;\; |x|>  1,$$ 
and the sequence of functions $(f_n)_{n\in\N_1}\subset\mathcal{C}_c^\infty(\Xt)$ with $f_n\geqslant 0,$ where $f_n(x) := \chi(n^{-1}|x|^p), \,p\geqslant 2$. We see that $ f_n\rightarrow 1$, and $\nabla f_n\rightarrow 0$ and $\LG_t^\nu f_n\rightarrow 0$ as $n\rightarrow\infty.$
Next, recall from item (ii) in Remark~\ref{Feller} that 
for each $f = |x|^p,\,p\geqslant 2$,   $\LG_t^{\nu} f\in \mathcal{C}(\Xt)$, and there exists $L_{b^\nu,\sigma^\nu,p}>0$ such that 
\begin{equation*}
\vert \LG_t^{\nu} f(x)\vert \leqslant L_{b^\nu,\sigma^\nu,p} \left( 1+\vert x\vert^p\right);
\end{equation*}
in particular, 
\begin{align}\label{GW26}
\vert \LG_t^\nu f_n(x)\vert\leqslant L_{b^\nu,\sigma^\nu}\left( 1+\vert x\vert^2\right).
\end{align}
Next, we use the local boundedness of $\varphi$ to obtain for $t\in \Ic$
\begin{align}\label{GW27}
\left\vert \int_{t_0}^t\int_{\Xt}\varphi(\eta_s(x))\LG_t^\nu f_n(x)\nu_s(dx)ds\right\vert \leqslant C_{\varphi}\int_{t_0}^t\int_{\Xt}\left\vert \LG_t^\nu f_n(x)\right\vert\nu_s(dx)ds.
\end{align}
By the bounds (\ref{GW26})-(\ref{GW27}), and the fact that $\int_\Xt|x|^2\nu_t(dx)<\infty$ (cf. Corollary \ref{sqrint}) due to the assumed regularity of the coefficients $(\brr^\nu,\sigma^\nu)$  in (\ref{SDE2}), we have  
\begin{align}\label{GW28}
\lim_{n\rightarrow \infty}\int_{t_0}^t\int_{\Xt}\varphi(\eta_s(x))\LG_t^\nu f_n(x)\nu_s(dx)ds =0
\end{align}
by the dominated convergence theorem. 

In order to bound the third term on the rhs of (\ref{GW24}), we observe that from the definition of~$f_n,$ we have 
\begin{align}\label{chibnd}
\vert \nabla f_n(x)\vert = n^{-1}\vert \nabla \chi(n^{-1}x)\vert \leqslant \Vert \nabla \chi\Vert_{\infty}<\infty,
\end{align}

\vspace{-.3cm}
\noindent so that 
\begin{align}\label{thtchibnd}
\big\langle \betaLk(s,x), \nabla f_n(x)\big\rangle \varphi^{\prime}(\eta_s(x))\eta_s(x)\leqslant \Vert \nabla \chi\Vert_{\infty} \big\vert \betaLk(s,x)\big\vert \big\vert \varphi^{\prime}(\eta_s(x))\big\vert \eta_s(x).
\end{align}
Consequently, using (\ref{thtchibnd}), (\ref{chibnd}), the uniform boundedness of the right inverse of $\sigma^\nu$ (so that $0<\|\sigma^\nu\|_\textsc{hs}\leqslant C_{\sigma^\nu}$), and the Young's inequality 
\begin{align*}
&\int_{t_0}^t\int_{\Xt}\big\langle \betaLk(s,x), \nabla f_n(x)\big\rangle\varphi^{\prime}(\eta_s(x))\eta_s(x)\nu_s(dx)\\ 
&\hspace{3cm} \leqslant \Vert \nabla \chi\Vert_{\infty} \int_{t_0}^t\int_{\Xt}\vert \tilde\sigma^{\nu,-1}\betaLk(s,x)\vert \vert \varphi^{\prime}(\eta_s(x))\vert \|\sigma^\nu(s,x)\|_\textsc{hs} \eta_s(x)\nu_s(dx)ds\\
&\hspace{3.cm}\leqslant \frat \, \Vert \nabla \chi\Vert_{\infty}\int_{t_0}^t\int_{\Xt} \Big(\vert \tilde \sigma^{\nu,-1}\betaLk(s,x)\vert^2 + C_{\sigma^\nu}\vert \varphi^{\prime}(\eta_s(x))\vert^{2}\, \eta^2_s(x)\Big) \nu_s(dx)ds<\infty,
\end{align*}
which, given the assumptions on $(b^\mu, \sigma^\mu)$,  $(b^\nu, \sigma^\nu)$, follows from Corollary \ref{sqrint} and local boundedness of $\varphi^{\prime}\in \mathcal{C}^1(\R)$ with $\eta_t<\infty$ (cf. Lemma~\ref{Pre_lem2}).  Thus, by the  dominated convergence theorem, we arrive~at 
\begin{align}\label{GW29}
\lim_{n\rightarrow\infty}\int_{t_0}^t\int_{\Xt}\langle \betaLk(s,x), \nabla f_n(x)\rangle\varphi^{\prime}(\eta_s(x))\eta_s(x)\nu_s(dx) = 0.
\end{align} 
Finally, put (\ref{GW25}), (\ref{GW28}), (\ref{GW29}) into the estimate (\ref{GW24}); since $f_n\underset{n\rightarrow \infty}{\longrightarrow} 1$, we have for $\varphi\in \mathcal{C}^2(\R^+)$
\begin{align*}
\hspace{1cm}\int_{\Xt}\varphi(\eta_t(x))\nu_t(dx) \leqslant \frat\int_{t_0}^t\int_{\Xt}\left\vert (\tilde{\sigma}^{\nu,-1}\betaLk)(s,x)\right\vert^2\varphi^{\prime\prime}(\eta_s(x))\eta_s^2(x)\nu_s(dx)ds,\hspace{2cm}
\end{align*}
which follows from Corollary \ref{sqrint} and the dominated convergence theorem. \qed

\subsubsection{\bf Bound on information loss via generator reconstruction when $\Xt\subset \XXt$}\label{MM-M}\mbox{}

In many applications involving so-called reduced-order models (obtained via averaging or homogenisation, or data-driven techniques), the original dynamics evolves on a higher-dimensional domain than its approximation; i.e., $\Xt\subset \XXt$ in (\ref{SDE2}). We consider $\Xt$ to be a linear subspace of $\XXt$ and set  $\mathfrak{M} = \Xt\times \mathcal{Y}$, with $\mathcal{Y}$ some finite-dimensional smooth manifold without boundary (usually $\mathcal{Y} = \R^{\ell-d}$, or $\mathcal{Y} = \bar{\mathbb{T}}^{\ell-d}$; $\ell = \dim\XXt$, $d = \dim\Xt$). By default, the time-marginal probability measures $\nu_t$ solving the forward Kolmogorov equation (\ref{fke2}) are defined on  $\Xt\subseteq\XXt$. 
Thus,  in order to carry out meaningful analysis when  $\Xt\subset \XXt$, and following the setup outlined in \S\ref{setp}, one has to consider projections of the time-marginal probability measures $\mmt\in \PP(\XXt)$ associated with the original dynamics onto probability measures  $\mu_t\in \PP(\Xt)$ which are obtained via marginalisation of their (Lebesgue) densities; i.e., for  $\xxt=(x,y)\in \Xt\times \mathcal{Y}$,  we set  
 \begin{equation}\label{margs}
 \mmt(dxdy) = \rrt(x,y)dxdy, \qquad  \mu_t(dx) = \rho_t^\mu(x)dx, \quad \rho^\mu_t(x) = \int_\mathcal{Y}\rrt(x,y)dy.
 \end{equation}
Below, we outline modifications of the results  presented in \S\ref{Xt=XXt}  for $\Xt=\XXt$ to the case $\Xt\subset \XXt$. 

\smallskip
Throughout, we assume that the coefficients $(\brr^\mm, \sigma^\mm)$ of the SDE (\ref{SDE1}) on $\XXt$, and the coefficients $(\brr^\nu, \sigma^\nu)$ of the SDE (\ref{SDE2}) on $\Xt\subset\XXt$ satisfy the same conditions as in \S\ref{Xt=XXt}, so that the respective solutions exist globally on $\Ic$ and are represented by flows of $\mathcal{C}^3$-diffeomorphisms. Thus, analogously to \S\ref{Xt=XXt}  the probability densities $(\rrt)_{t\in \Ic}$ and $(\rho^{\nu}_{t})_{t\in \Ic}$ w.r.t.~the Lebesgue measures on, respectively, $\XXt$ and $\Xt$ exist  and satisfy (in the distributional  sense on $\mathcal{E}_+$)  
\begin{align}
\partial_t\hspace{0.05em}\rrt = \LG^{\mm*}_{t}\rrt = \textstyle \frac{1}{2}\partial^2_{\xxt_i \xxt_j}(a^{\mm}_{ij}\hspace{0.09em}\rrt) - \partial_{\xxt_i}({b}_i^{\mm}\hspace{0.09em}\rrt), \qquad \rrt\in L^1_+(\XXt,d\xxt)\cap L^\infty(\XXt;d\xxt),\label{fke111}\\[.2cm]
\partial_t\rho^{\nu}_{t} = \LG^{\nu*}_{t}\rho^{\nu}_{t} = \textstyle\frac{1}{2}\partial^2_{x_ix_j}(a^{\nu}_{ij}\rho^{\nu}_{t}) - \partial_{x_i}({b}_i^{\nu}\rho^{\nu}_{t}), \qquad \rho_{t}^\nu\in L^1_+(\Xt,dx)\cap L^\infty(\Xt;dx).\label{fke222}
\end{align}  

First, note that the evolution of the density $\rho^\mu_t$ of the time-marginal measure $\mu_t\in \PP(\Xt)$ can still be formally written in the form (\ref{fke111}), i.e.,
\begin{equation}\label{nfke}
\partial_t\rho^\mu_t(x) = \mathcal{L}^{\mu*}_{\rho,t} \,\rho^\mu_t(x), \qquad \rho_{t_0}^\mu\in L^1_+(\Xt,dx)\cap L^\infty(\Xt;dx),
\end{equation}
 with the operator  $\mathcal{L}^{\mu*}_{\rho,t}$ defined via the coefficients  $(b^\mu_\rho, \sigma^\mu_{\rho})$ given by 
\begin{align*}
b^\mu_\rho(t,x) = \int_\mathcal{Y} b^\mm\big(t,(x,y)\big)\hspace{.05em}\rrt(y|x)dy, \quad \sigma_{\rho}^\mu(t,x)= \int_\mathcal{Y}\sigma^\mm\big(t,(x,y)\big)\rrt(y|x)dy, \quad x\in \Xt, \;y\in \mathcal{Y},
\end{align*}
where $\rrt(y|x)$ is determined from the joint density $\rrt(x,y) = \rrt(y|x)\rho^\mu_t(x)$ solving the forward Kolmogorov equation (\ref{fke111}) for the original dynamics of (\ref{SDE1}) on $\XXt$. If  the branches of the conditional density $x\mapsto \rrt(y|x)$ are continuously differentiable;  the regularity of  $b^\mm(t,(x,y))$ and $\sigma^\mm(t,(x,y))$ ensures that  $b_\rho^\mu(t,x)$, $\sigma_{\rho}^\mu(t,x)$  are sufficiently regular, which is what is required in Theorem~\ref{Info_ineq2}. In other words,  the time-marginal measure $\mu_t\in \PP(\Xt)$  does not have to be generated  by  a flow of solutions of  an SDE, as long as  
\begin{equation}\label{weak_sol_Xt2}
\partial_t\rho_t^{\mu} = \LG^{\nu*}_{t}\rho_t^{\mu}-\nabla_x\cdot(\tilde{\varTheta}_{\mu\nu}\rho_t^{\mu}), \qquad \rho_{t_0}^{\mu} \in L_+^1(\Xt;dx)\cap L^\infty(\Xt;dx), 
\end{equation}
obtained analogously to (\ref{RecEQ}) in Lemma \ref{THT_def} with 
\begin{align}\label{Per_vF2}
\tilde{\varTheta}_{\mu\nu} = \frat\left( a_{ij}^\nu -(a^\mu_\rho)_{ij}\right)\nabla_x\log\rho^\mu_t -\left( {b}_i^\nu -\partial_{x_j}a_{ij}^\nu -({b}^{\mu}_\rho)_i +\partial_{x_j}(a^{\mu}_\rho)_{ij}\right),
\end{align}
is solvable on $\Xt$. Similar to Corollary~\ref{summ_cor},  the solvability of  (\ref{weak_sol_Xt2}) is controlled by the regularity of the coefficients $(\brr^\nu,\sigma^\nu)$ of the approximate SDE dynamics (\ref{SDE2}), which is defined on $\Xt$ by default, and the boundedness of the logarithmic gradient in the reconstructed field $\tilde{\varTheta}_{\mu\nu}$~(\ref{Per_vF2}). Given that  $\nabla_{\xxt}\log\rrt(\xxt)\equiv\nabla_{(x,y)}\log\rrt(x,y)= \nabla_{(x,y)}\log\rrt(y|x) +\nabla_{x}\log\rho^\mu_t(x)$,  the boundedness of $\nabla_x\log\rho^\mu_t(x)$ on $\Xt$ follows from the boundedness of $\nabla_{\xxt}\log\rrt(\xxt)$ on $\XXt$ which was established in Proposition~\ref{log_gr} (which only depends on the dynamics on $\XXt$). 
 
 Note further that, similar to the setup in \S\ref{Xt=XXt}, for our purposes of obtaining the bound~(\ref{InformaBB}) and minimising it in terms of the coefficients $(\brr^\nu,\sigma^\nu)$ of (\ref{SDE2}), the densities $(\rrt)_{t\in \Ic}$ can be assumed to be known. In a more general case,  (\ref{nfke}) becomes a nonlinear Fokker-Planck-Kolmogorov equation (e.g., \cite{Roc-Kry}) associated with a mean-field, McKean-Vlasov type SDE with coefficients $(b^\mu_\rho, \sigma^\mu_\rho)$.  Although such a generalisation is not needed here,  an analogue of Proposition \ref{log_gr} could be derived by means of a Bismut-Elworthy-Li formula for mean-field SDE dynamics (e.g.,~\cite{banos18}).

Proposition \ref{Summary_prop} in \S\ref{Xt=XXt}, establishing the solvability of (\ref{weak_sol_Xt}), relies on the regularity of the coefficients $(\brr^\nu,\sigma^\nu)$ of the SDE (\ref{SDE2}) defined on $\Xt$,  and the square-integrability of $|\tilde\sigma^{\nu,-1}\betaLk|$ discussed in Corollary \ref{sqrint}. Thus, once the boundedness of $\nabla_x\log\rho^\mu_t(x)$ is asserted (see above), the square-integrability, boundedness, and smoothness of  $|\tilde\sigma^{\nu,-1}\tilde{\varTheta}_{\mu\nu}|$ follows by the same steps as those in Corollary \ref{sqrint}.

Lemmas \ref{Ide_Lem}-\ref{Integrand} rely on manipulating differential operators acting on functions on $\Xt$, and on the existence and finiteness of $\rho^\mu_t(x)/\rho^\nu_t(x)$, $t>t_0$, which follows from the properties of solutions, $\rrt(\xxt)$, $\rho^\nu_t(x)$, of the forward Kolmogorov equations (\ref{fke111}) and (\ref{fke222}) on the original domains, respectively, $\XXt$ and $\Xt$. In particular, $0<\rrt(\xxt)\equiv\rrt(x,y)<\infty$ implies $0<\rho_t^\mu(x)<\infty$. 

Finally, the proof of Theorem \ref{Info_ineq2} follows analogously to that in \S\ref{Xt=XXt} with $\betaLk$ in (\ref{Per_vF}) replaced by~$\tilde{\varTheta}_{\mu\nu}$ in (\ref{Per_vF2}). We consider an example associated with such a configuration in \S\ref{toysec}.

\subsection{Bound on information loss in Lagrangian predictions via $\varphi$-divergence rate fields}\label{ftdr_sec}
Here,  we derive another bound on $\Df(\mu_t\|\nu_t)$, $\mu_t, \nu_t\in \PP(\Xt)$, which  can be combined with
\begin{align}\label{repr_bound}
\hat{\mathcal{K}}_{\varphi,f}^{\nu}\big(-\Df(\mu_t\|\nu_t)\big)\leqslant  \E^{\mu_t}[f] - \E^{\nu_t}[f] \leqslant \mathcal{K}_{\varphi,f}^{\nu}\big(\Df(\mu_t\|\nu_t)\big),  
\end{align}
derived in \S\ref{Info_ineq} (see (\ref{Csiszar})) in order to minimise the error in the estimates $\E^{\nu_t}[f]$ of $\E^{\mu_t}[f]$. Similar to \S\ref{Info_ineq}, these results do not rely on the dynamics being generated by SDE's/ODE's, and they generalise to path space probability measures, as discussed in \S\ref{Path-space}. Using the same approach as in the previous sections, if the approximate dynamics evolves on the subspace $\Xt\subset \XXt$, the time-marginal probability measure $\mu_t\in \PP(\Xt)$ is obtained from the original probability measure $\mmt\in \PP(\XXt)$ by marginalisation of its density as in (\ref{margs}). Specific links to path-based observables are achieved from the relationships  $\E^{\mu_t}[f] = \E[f\big(\pi^\nu_\mm\circ\Solm^\mm\big)]$, $\E^{\nu_t}[f] = \E[f\big(\Solm^\nu\big)]$, with $\pi^\nu_\mm:\XXt\rightarrow \Xt$ the natural projection, as highlighted in \S\ref{setp} and detailed in \S\ref{SFlow}. The information about the initial conditions is propagated under the action of the transition evolutions; namely $\mmt = \mathcal{P}^{\mm*}_{t_0,t}\,\mmo$ and $\nu_t = \mathcal{P}^{\nu*}_{t_0,t}\,\nu_{t_0}$ for all $t\in \Ic=\big[t_0, \,t_0+T\big)$; see Definition \ref{trans_evo}. 

\smallskip
The bound obtained below is expressed in terms of the difference between $\Df(\mu_t\|\mu_{\hspace{.02cm}t_0})$ and $\Df(\nu_t\|\nu_{\hspace{.02cm}t_0})$, where $\mu_{t_0} = \nu_{t_0} $, in the form
\begin{align}\label{FTDR_intro}
\Df(\mu_t\|\nu_t)\leq \big\vert \Df(\mu_t\|\mu_{\hspace{.02cm}t_0})-\Df(\nu_t\|\nu_{\hspace{.02cm}t_0})\big\vert, \quad  \forall\,t\in \Ic=\big[t_0, \,t_0+T\big).
\end{align}
Importantly, (\ref{FTDR_intro}) is useful in a computational framework aimed at Lagrangian (path-based) uncertainty quantification and model tuning, and it can be evaluated, in both the stochastic and deterministic setting, by minimising the discrepancy between two scalar fields of {\it Finite-Time $\varphi$-Divergence Rates} ($\varphi$-FTDR) defined~by 
 $$x\mapsto \Df(\mu^x_t\|\varkappa^{\hspace{.04cm}x}_{\hspace{.02cm}t_0}), \qquad x\mapsto \Df(\nu^x_t\|\varkappa^{\hspace{.04cm}x}_{\hspace{.02cm}t_0}), \qquad x\in \Xt,$$
for $\mu^x_t$ and $\nu^x_t$ evolving from the $\varkappa^{\hspace{.04cm}x}_{\hspace{.02cm}t_0}\in \PP(\Xt)$ concentrated on a neighbourhood of $x\in \Xt$ \cite{MBUda18}. 

\begin{theorem}\label{difference_bound}
Let $\varphi\in \mathcal{C}^2(\Rp)$ be a strictly convex function satisfying the normality conditions~(\ref{Normality}).  Let $\mu_t,\nu_t\in  \PP(\Xt)$   evolve from $\nu_{t_0} \,{=}\, \mu_{t_0} $. If $\Df(\mu_t\|\nu_t)<\infty$, $\Df(\mu_t\|\mu_{\hspace{.02cm}t_0})<\infty$, and $\Df(\nu_t\|\nu_{\hspace{.02cm}t_0})<\infty$ for all $t\in\Ic=[t_0, \,t_0+T)$, then 
\begin{align}\label{ftdr_bnd}
\Df(\mu_t\|\nu_t)\leqslant \big\vert \Df(\mu_t\|\mu_{\hspace{.02cm}t_0})-\Df(\nu_t\|\nu_{\hspace{.02cm}t_0})\big\vert \quad \forall \,t\in \Ic.
\end{align}
\noindent {\it Proof.}\rm  \;\;See Appendix \ref{difference_bound_app}; this result relies on joint convexity of $\Df$ in its arguments and it does not require the time-marginal measures $\mu_t$, $\nu_t$,  to be associated with a Markov process. 
\end{theorem}

\begin{rem}\rm The following comments are in order:
\begin{itemize}[leftmargin=0.7cm]
\item[(i)] The bound in (\ref{ftdr_bnd}) is non-uniform in $T$, unless the underlying dynamics have a stationary or a cyclo-stationary measure.
\item[(ii)] When the original and approximate dynamics are generated by the SDEs (\ref{SDE1}) and (\ref{SDE2}), and $\mu_t$, $\nu_t$, with $\mu_{t_0} = \nu_{t_0}$,  solve their respective forward Kolmogorov equations,
the sufficient condition for $\Df(\mu_t\|\nu_t)<\infty$  $\forall \,t\in \Ic$ is that the fields $\brr^\mm(t,\,\cdot\,)$, $\sigma^\mm(t,\,\cdot\,)$, $\brr^\nu(t,\,\cdot\,)$, $\sigma^\nu(t,\,\cdot\,)$  satisfy conditions of Theorem \ref{Info_ineq2}.
These conditions also lead to $\Df(\mu_t\|\mu_{t_0}),\Df(\nu_t\|\nu_{t_0})<\infty$. For deterministic dynamics, i.e., $\sigma^\mm=\sigma^\nu=0$, the evolution of the uncertainty in  Lagrangian predictions is given by the push forward of the measures, $\mmo$, $\nu_{t_0}$, on the initial conditions  under, respectively, $\phi_{t_0,t}^\mm(\xxt,\om) = \psi_{t_0,t}^\mm(\xxt)$ and $\phi_{t_0,t}^\nu(x,\om) = \psi_{t_0,t}^\nu(x)$.

\item[(iii)] The general bound in (\ref{repr_bound}) combined with  (\ref{ftdr_bnd}) allows us to optimize the accuracy of Lagrangian predictions which are based on the approximate dynamics. If the time-marginal measures are associated with the laws of paths  $t\mapsto \phi^\mm_{t_0,t}$, $t\mapsto \phi^\nu_{t_0,t}$ (where $\phi^{\scriptscriptstyle(\scaleobj{1.3}{\cdot})}_{t_0,t}$ is not necessarily a stochastic flow), $\varphi$-FTDR fields provide a probabilistic way of quantifying local, finite-time expansion rates between neighbouring trajectories in both deterministic and stochastic cases, as introduced and analysed in \cite{MBUda18}.  Thus,  minimising the discrepancy between  $\varphi$-FTDR fields generated on $\Xt$ by the original dynamics and its approximation allows for minimising the uncertainty/error in the Lagrangian observables $x\mapsto \E\big[f(\phi^\nu_{t_0,t}(x)\big]$ relative to   $x\mapsto \E\big[f(\pi^\nu_\mm\circ\phi^\mm_{t_0,t}(x)\big]$ with  $\pi^\nu_\mm\,{:} \;\XXt\rightarrow \Xt$ the natural  projection onto $\Xt\subseteq\XXt$. 
\end{itemize}

\noindent In order to provide a concise example, assume that $\Xt=\XXt$, and the maps $\phi^\mm_{t_0,t}\equiv\phi^\mu_{t_0,t}$, $\phi^\nu_{t_0,t}$ are given by stochastic flows. Then, following \S\ref{SFlow}, the evolution of the considered time-marginal measures is given by  $\mu_t^x = \mathcal{P}^{\Phi^x*}_{t_0,t}\mu_{B_\varepsilon(x)}$ and $\nu_t^x = \mathcal{P}^{\Psi^x*}_{t_0,t}\mu_{B_\varepsilon(x)}$, where the initial probability measure is taken to be supported on the ball $B_\varepsilon(x)$  of radius $\varepsilon$ centred at $x\in \Xt$ and regularised by convolving it with an arbitrary Gaussian. Here, $\mathcal{P}^{\Phi^x}_{t_0,t}$  is the transition evolution (\ref{P}) induced by the two-point motion $\Phi^x_{t_0,t}(v,\omega): = \phi^\mu_{t_0,t}(x+v,\om)-\phi^\mu_{t_0,t}(x,\om)$ with $\phi^\mu_{t_0,t}$ the stochastic flow associated with the SDE (\ref{SDE1}). Similarly,  $\mathcal{P}^{\Psi^x}_{t_0,t}$  is the transition evolution induced by the two-point motion  $\Psi^x_{t_0,t}(v,\omega): = \phi^\nu_{t_0,t}(x+v,\om)-\phi^\nu_{t_0,t}(x,\om)$ with $\phi^\nu_{t_0,t}$ the stochastic flow associated with the SDE (\ref{SDE2}). Then, the criterion for optimising Lagrangian predictions on $\Xt$ over some time interval $\Ic$ can be obtained from the bound  (\ref{ftdr_bnd}) in the form 
\begin{align}
&\int_\Ic\int_\Xt\Df\Big(\mathcal{P}^{\Phi^x*}_{t_0,t}\mu_{B_\varepsilon(x)}\|\mathcal{P}^{\Psi^x*}_{t_0,t}\mu_{B_\varepsilon(x)}\Big)dxdt \notag\\
&\hspace{2cm}\leqslant \int_\Ic\int_\Xt\Big| \Df\Big(\mathcal{P}^{\Phi^x*}_{t_0,t}\mu_{B_\varepsilon(x)}\|\mu_{B_\varepsilon(x)}\Big)-\Df\Big(\mathcal{P}^{\Psi^x*}_{t_0,t}\mu_{B_\varepsilon(x)}\|\mu_{B_\varepsilon(x)}\Big)\Big| dxdt, \notag
\end{align}
where $x\mapsto\Df\big(\mathcal{P}^{\Phi^x*}_{t_0,t}\mu_{B_\varepsilon(x)}\|\mu_{B_\varepsilon(x)}\big)$,  $x\mapsto\Df\big(\mathcal{P}^{\Psi^x*}_{t_0,t}\mu_{B_\varepsilon(x)}\|\mu_{B_\varepsilon(x)}\big)$, correspond to $\varphi$-FTDR fields discussed in \cite{MBUda18}. Crucially, such a bound is amenable to computational treatment utilising algorithms for fast approximations of the `transfer operator' (e.g., \cite{dellnitz99b}) represented here by the duals  $\mathcal{P}^{\Phi^x*}_{t_0,t}$,  $\mathcal{P}^{\Psi^x*}_{t_0,t}$,  of $\mathcal{P}^{\Phi^x}_{t_0,t}$,  $\mathcal{P}^{\Psi^x}_{t_0,t}$; specific applications will be considered in future~work.
\end{rem}

\subsection{Information bounds on path space }\label{Path-space}
Here, we  extend the results derived in the previous sections to probability measures induced by the underlying dynamics on the path space. In order to simplify exposition, 
we restrict the discussion to the case of $\XXt=\Xt$ in (\ref{SDE1}) and (\ref{SDE2}); i.e., the original dynamics and its approximation are defined on the same domain, and $\mathcal{W}_d=C(\Ic;\Xt)$.

\smallskip
Following the notation, definitions and background results of \S\ref{SFlow} (see also Glossary),  we consider path space measures $\Ps^{\mu}_{t_0}, \Ns^{\nu}_{t_0}\in \PP(\mathcal{W}_d)$ defined via the solutions $\Ps_{t_0,x}^\mu$, $\Ns_{t_0,x}^\nu$ of the martingale problem (see \S\ref{SFlow} and Definition \ref{mrt_prb}) generated by the SDEs (\ref{SDE1}) and (\ref{SDE2}), so that  
$$\Ps^{\mu}_{t_0}(d\om) \,{=} \,\int_{\Xt}\Ps_{t_0,x}^\mu(d\om)\mu_{t_0}(dx), \qquad  \Ns^{\nu}_{t_0}(d\om)\,{=}\,\int_{\Xt}\Ns_{t_0,x}^\nu(d\om)\nu_{t_0}(dx),\qquad \mu_{t_0}, \nu_{t_0}\in \PP(\Xt),$$ 
 with finite-dimensional distributions $\Ps^\mu_{\Ic_{n}}, \Ns^\nu_{\Ic_{n}}$  defined on $\otimes_{i=1}^nA_i$, $A_i\in \Bb(\Xt)$, by 
\begin{align*}
\Ps^\mu_{\Ic_{n}}(\otimes_{i=1}^nA_i) &:= \int_\Xt\!\Ps^{\mu}_{t_0,x}\big\{\om: \phi^\mu_{s_1,t_0}(x,\om)\in A_1,\cdots, \phi^\mu_{s_n,t_0}(x,\om)\in A_n\big\}\mu_{t_0}(dx), \\[.2cm]
\Ns^\nu_{\Ic_{n}}(\otimes_{i=1}^nA_i) &:= \int_\Xt\!\Ns^{\nu}_{t_0,x}\big\{\om: \phi^\nu_{s_1,t_0}(x,\om)\in A_1,\cdots, \phi^\nu_{s_n,t_0}(x,\om)\in A_n\big\}\nu_{t_0}(dx),
\end{align*}
where  $t_0<s_1<\cdots<s_n\leqslant t_0+T$, and $X^\mu_{t_0,t}(\om) = \phi^\mu_{t_0,t}(x,\om)$, $X^\nu_t = \phi^\nu_{t_0,t}(x,\om)$ solve, respectively, (\ref{SDE1}) and (\ref{SDE2}) on $\Ic{=}\,\big[t_{0},\,t_0\,{+}\,T\big)$ with coefficients satisfying the same conditions as those in  Theorem \ref{Info_ineq2}, so that global solutions on $\Ic$ exist and are represented by a flow of $\mathcal{C}^3$-diffeomorphisms.

\smallskip
First, note that the information bound derived in \S\ref{Info_ineq} applied to arbitrary probability measures, provided that the original measure was absolutely continuous with respect to the approximating measure. Thus, in particular, the bound (\ref{Csiszar}) applies to path space probability measures $\Ps^{\mu}_{t_0}, \Ns^{\nu}_{t_0}\in \PP(\mathcal{W}_d)$ (in lieu of $\mu_t, \nu_t\in \PP(\Xt)$)  as long as $\Ps^{\mu}_{t_0}\ll \Ns^{\nu}_{t_0}$, which is discussed below. 

\smallskip
We also derive an  identity which yields a unique projection $\Ms^{\mu\nu}_{t_0}$ of the path space probability measure $\Ns^\nu_{t_0}$ onto a closed convex subset of $\PP(\Wd)$ of path space probability measures with  time-marginals $(\mu_{t})_{{t\in \Ic}}$ solving the forward Kolmogorov equation associated with (\ref{SDE1}); namely 
\begin{align*}
\Df\big(\Ps^{\mu\nu}_{t_{0}}\|\Ns^\nu_{t_{0}}\big):&=\inf\left\{ \Df\big(\Ps\|\Ns^\nu_{t_{0}}\big): \Ps\ll \Ns^{\nu}_{t_{0}},\;\;\Ps\circ \phi_{t_0,t}^{\nu,-1} = \mu_t, \;\; \forall\,t\in \Ic\right\}\notag\\[.2cm]
& \hspace{0cm}= \sup\bigg\{ \int_{\Xt}\sum_{i=1}^n f_i(x)\mu_{s_i}(dx) -\int_{\Xt}\varphi^*\Big(\sum_{i=1}^n f_i(x)\Big) \nu_{s_i}(dx);\notag\\[.0cm] 
&\hspace{3cm} \forall \;1\leqslant n <\infty ,\;f_{1}, \dots, f_{n}\in \mathbb{M}_\infty(\Xt), \;   s_1,\dots,s_{n}\in \Ic\bigg\},
\end{align*}
where the supremum  is over all $n$-tuple partitions of $\Ic$, and $n$-tuples of functions in $\mathbb{M}_\infty(\Xt)$\footnote{\,This could be extended to $f_i\in L_{\varphi*}(\XX;\nu_{s_i}),\; \forall \;1\leqslant n <\infty$.}, and  the marginal measures $(\mu_{s_i})_{i=1}^n$ solving (\ref{F_Kol}a), and $(\nu_{s_i})_{i=1}^n$ solving (\ref{F_Kol}b).  Importantly, we show that the $\varphi$-projection, defined via $\varphi$-divergence between solutions of the corresponding forward Kolmogorov equations with $\mu_{t_0}=\nu_{t_0}$, can be linked to $\varphi$-FTDR fields (\S\ref{ftdr_sec}) via  
\begin{align}\label{Inpath}
\Df\big(\Ps^{\mu\nu}_{t_{0}}\|\Ns^\nu_{t_{0}}\big)\leqslant\sup\bigg\{\sum_{i=1}^n\big\vert \Df(\mu_{s_i}\|\mu_{t_0})-\Df(\nu_{s_i}\|\nu_{t_0})\big\vert, \;\;\forall \;1\leqslant n <\infty , \;   s_1,\dots,s_n\in \Ic\bigg\},
\end{align}
where the supremum in (\ref{Inpath}) is over all $n$-tuple partitions of $\Ic$. Moreover, for finite-dimensional distributions, one has a more practically useful bound  
\begin{align}\label{Inphs}
\Df\big(\Ps^{\mu\nu}_{\Ic_n}\|\Ns^\nu_{\Ic_n}\big)\leqslant \sum_{i=1}^n\big\vert \Df(\mu_{t_i}\|\mu_{t_0})-\Df(\nu_{t_i}\|\nu_{t_0})\big\vert, \quad  t_1,\dots,t_n\in \Ic.
\end{align}
The information identities (\ref{Inpath}) and (\ref{Inphs}) yield a suitable  way of quantifying  model error for stochastic flows generated by SDEs in path space and phase space, respectively.  Note that while the path space bounds are more difficult to deal with in practice, they take into account more information than the bounds in \S\ref{Inf_IP}\,--\,\ref{ftdr_sec} based on families of one-point marginals $(\mu_t)_{t\in \Ic}$, $(\nu_t)_{t\in \Ic}$.

\medskip
We start by recalling some facts about absolute continuity of probability measures on path~space. 

\begin{lem}[\cite{Pinsky04}]\label{Absestm}
Let $\Ms, \Ns\in \PP(\mathcal{W}_d)$ and $(\F_n)_{n\in\N}$ be a filtration on $\mathcal{W}_d$ such that for each $A_n\in \F_n$ we have $\mathfrak{S}\left(\bigcup_{n\in\N}A_n\right) = {\mathfrak{S}}(\{ \phi_{t_0,\cccdot}\}) =\F_t$. If $\Ms|_{\F_n}\ll \Ns|_{\F_n}$ for all $n\in \N$, then 
\btm
\item[(a)] $\Ms\ll \Ns$ iff \;$\limsup_{n\rightarrow \infty}\frac{\vphantom{\big|}\scaleobj{1.3}{ d\Ms}|_{\F_n}}{\vphantom{\big|}\scaleobj{1.3}{d\Ns}|_{\F_n}}<\infty\quad \Ns\; \text{-a.s.}$,

\vspace{.2cm}\item[(b)] $\Ms\perp \Ns$ iff \;$ \limsup_{n\rightarrow \infty}\frac{\vphantom{\big|}\scaleobj{1.3}{d\Ms}|_{\F_n}}{\vphantom{\big|}\scaleobj{1.3}{d\Ns}|_{\F_n}}=\infty \quad \Ns\; \text{-a.s.}$,

\vspace{.2cm}\item[(c)] if $\Ms\ll \Ns$ and $\varphi\left( \frac{\vphantom{\big|}\scaleobj{1.3}{d\Ms}}{\vphantom{\big|}\scaleobj{1.3}{d\Ns}}\right)\in L^1(\mathcal{W}_d; \Ns),$ then $\Df(\Ms\|\Ns) = \lim_{n\rightarrow\infty}\Df\big(\Ms|_{\F_n}\|\Ns|_{\F_n}\big).$

\etm
\end{lem}
\noindent {\it Proof.} See Appendix \ref{Absestm_Ap}; (a) and (b) follow from Lebesgue decomposition, and (c) follows from Fatou's lemma applied to $\Df\big(\Ms|_{\F_n}\|\Ns|_{\F_n}\big)$, strict convexity of $\varphi$, and Jensen's inequality.

\begin{lem}[\cite{Pinsky04, Yor}]\label{Martinest}
Let $\Mt_t$ be a continuous local martingale w.r.t.~$(\mathcal{W}_d, \F_t)$ and let $\langle \Mt\rangle_t$\footnote{\,Throughout this section $\langle \ccdot \rangle_t$ denotes the quadratic variation at time $t$, while $\langle \ccdot, \ccdot\rangle$ denotes the inner product.} be the corresponding quadratic variation process such that $\langle \Mt\rangle_{\infty} := \lim_{t\rightarrow\infty}\langle \Mt\rangle_t.$ 

\smallskip
If $\mathcal{E}(\Mt_t) = \exp\left(\Mt_t-\frac{1}{2}\langle \Mt\rangle_t\right)$,
then
\begin{align*}
\Big\{\om: \lim_{t\rightarrow\infty}\mathcal{E}(\Mt_t) =0\Big\} = \Big\{\om: \langle \Mt\rangle_{\infty} =\infty\Big\} \qquad \mathbb{P}\text{\,-\,a.s.}
\end{align*}
\end{lem}
\noindent {\it Proof.}  See Appendix \ref{Martinest_Ap}; this  follows from Fatou's lemma for conditional expectation and Doob's theorem, the fact that $\mathcal{E}(\Mt_t)$ is a supermartingale, and the identity $\mathcal{E}(-\Mt_t)\!= \!\mathcal{E}(\Mt_t)^{-1}\!\exp\left(\langle \Mt\rangle_t\right)$.

\medskip
The following proposition is reminiscent of the reconstruction procedure developed  in \S\ref{Inf_IP}. It also yields a formula for $\Df\big(\Ms_{t_0,x}^{\mu\nu}\|\Ns_{t_0,x}^\nu\big)$ where $\Ms_{t_0,x}^{\mu\nu}, \,\Ns_{t_0, x}^\nu\in \PP(\Wd)$ are extremal martingale solutions for two It\^o SDEs with coefficients\footnote{\,See \S\ref{Inf_IP} for details on the reconstructed field $\Theta_{{\mu\nu}}$.} $(b^\nu+\Theta_{{\mu\nu}}, \sigma^\nu)$  and $( b^\nu, \sigma^\nu)$ satisfying conditions of Theorem \ref{Info_ineq2}; see \S\ref{SFlow} and the Glossary for  background results and definitions.  The general procedure is not entirely new, but we adapt it to the current setup, since it is not immediately obvious how absolute continuity in the phase space $\Xt$ relates to the absolute continuity in the path space~$\mathcal{W}_d.$   
\begin{prop}\label{DifFO}
Assume that the coefficients in the SDE's (\ref{SDE1}) and (\ref{SDE2}) satisfy conditions of Theorem \ref{Info_ineq2},
 and let $\big(( \Ms_{t_0,x}^{\mu\nu})_{x\in \Xt}, \brr^\nu+\betaLk, \sigma^\nu\big)$, $\big((\Ns_{t_0,x}^\nu)_{x\in\Xt}, \brr^\nu, \sigma^\nu\big)$, with $\betaLk$ defined in (\ref{Per_vF}), be families of extremal martingale solutions on $\Ic\times\Xt$  starting from the same initial condition  $x\in\Xt$ at time $t_0\in \Ic$; furthermore, set $\betaL(t,x) = \big({a}^{\nu,-1}\betaLk\big)(t,x)$, where ${a}^{\nu} = \sigma^\nu(\sigma^\nu)^*$. Then 
\btm[leftmargin = .9cm]
\item[(i)] $\Ms_{t_0,x}^{\mu\nu}\perp \Ns_{t_0,x}^\nu$ iff
\begin{align*}
\int_{t_0}^t\left\langle\betaL(s,x), a^\nu\betaL(s,x)\right\rangle\big(s,\phi^\nu_{t_0,s}(s,\om)\big)ds =\infty \qquad \Ns_{t_0,x}^\nu\,\text{-\,a.s.,} \quad t\in\Ic.
\end{align*}
\item[(ii)] $\Ms_{t_0,x}^{\mu\nu}\ll \Ns_{t_0,x}^\nu$ iff 
\begin{align*}
\int_{t_0}^t\left\langle\betaL(s,x), a^\nu\betaL(s,x)\right\rangle\big(s,\phi^\nu_{t_0,s}(s,\om)\big)ds <\infty \qquad \Ns_{t_0,x}^\nu\,\text{-\,a.s.,}\quad t\in \Ic.
\end{align*}
\item[(iii)] Moreover, assume that $\varphi\in \mathcal{C}^2(\Rp)$ is a strictly convex function satisfying the normality conditions (\ref{Normality}), and  such that $\varphi\left(d\,\Ms_{t_0,x}^{\mu\nu}\big/d\,\Ns_{t_0,x}^\nu\right)\in L^1(\mathcal{W}_d;\Ns_{t_0,x}^\nu)$. Then, 
 \begin{align}\label{Formul1}
\Df\left(\Ms_{t_0,x}^{\mu\nu}\|\Ns_{t_0,x}^\nu\right) =  \frac{1}{2}\E^{\Ns_{t_0,x}^\nu}\bigg( \int_{t_0}^t\big\langle \betaL, a^\nu\betaL\big\rangle\big(s,\phi_{t_0,s}^\nu(x)\big) \varphi^{\prime\prime}(D_s)D_s^2ds\bigg), \quad t\in\Ic,
\end{align}
where $D_s = d\Ms^{\mu\nu,s}_{t_0,x}/d\Ns^{\nu,s}_{t_0,x}$ with $\Ms_{t_0,x}^{\mu\nu,s} :=\Ms_{t_0, x}^{\mu\nu}|_{\F_{s}}$, $\Ns_{t_0,x}^{\nu,s}:=\Ns_{t_0,x}^\nu|_{\F_{s}}$, we skipped the explicit $(s,x)$-dependence in $\betaL$, and  $\E^{\Ns_{t_0,x}^\nu}[f(x)\phi_{t_0,t}(x)]:=\int_\Om f(x)\phi_{t_0,s}(x,\om) \Ns_{t_0,x}^\nu(d\om)$.
 \etm
\end{prop}
\noindent {\it Proof.}  See Appendix \ref{DifFO_app}. Parts (i)-(ii) follow from Lemmas \ref{Absestm} and \ref{Martinest}. Part (iii) is more involved and it relies on a localisation procedure applied to $\E^{\Ns^\nu_{t_0,x}}\big[\varphi(D_{t})\big]$.

\begin{definition}[\textit{\textbf{$\varphi$-admissible flows of probability measures  and $\varphi$-projection}}]
\mbox{}\rm

\noindent Let $\{s_i\}_{i=0}^n\subset \Ic$ be any $n$-tuple partition of $\Ic$ for all $n\in \mathbb{N}$, and  consider the set 
 \begin{align}\label{Admm}
\mathbb{C}_{\varphi,\Ic}^{{\mu\nu}} :=\Big\{\Ms\in\PP\big(\mathcal{W}_d\big): \;\;\Df\big(\Ms\|\Ns_{t_0}^\nu\big)<\infty, \;\; \Ms\circ \phi^{\nu,-1}_{s_i,t_0} \,{=}\, \mu_{s_{i}} \,{\in}\,\PP(\Xt), \;\;\Ns_{t_0}^\nu\in \PP(\mathcal{W}_d)\Big\}. 
\end{align} 
\btm[leftmargin = 0.7cm]
\item[(a)] We say that the family of time-marginal probability measures  $(\mu_{t})_{t\in \Ic}$ on $\Xt$ is  $\varphi$-admissible if $\mathbb{C}_{\varphi,\Ic}^{{\mu\nu}}$ is a nonempty subset of $\PP(\mathcal{W}_d)$.

\vspace{.2cm}\item[(b)] A measure $\Ms^{\mu\nu}_{t_0}\in \mathbb{C}_{\varphi,\Ic}^{{\mu\nu}}$ is called a $\varphi$-projection of $\Ns_{t_{0}}^\nu$ onto $\mathbb{C}_{\varphi,\Ic}^{{\mu\nu}}$ if 
\begin{align*}
\Df\big(\Ms^{\mu\nu}_{t_{0}}\|\Ns_{t_{0}}^\nu\big) = \inf\Big\{\Df\big(\Ms\|\Ns_{t_{0}}^\nu\big)\,{:}\;\;\; \;\Ms\in \mathbb{C}_{\varphi,\Ic}^{{\mu\nu}}, \;\;\Ns_{t_0}^\nu\in \PP(\mathcal{W}_d)\Big\}.
\end{align*}
\etm
\end{definition}
\begin{theorem}\label{path_proj}
Let  $\{\mu_{t}\}_{t\in \Ic}, \{\nu_{t}\}_{t\in\Ic}$, denote families of time-marginal probability measures on $\Xt$  solving weakly the forward Kolmogorov equations associated, respectively, with the SDEs (\ref{SDE1}) and (\ref{SDE2}), where $\mu_{t_0}(dx) \,{=}\, \nu_{t_0}(dx)=\rho_{t_0}(x)dx$,  $\rho_{t_0}\,{\in}\, L^{1}_{+}(\Xt;dx)\cap L^{\infty}(\Xt;dx)$. 

Assume further that the coefficients in (\ref{SDE1}) and (\ref{SDE2}) satisfy conditions of Theorem \ref{Info_ineq2}, and let $\big((\Ns_{t_0,x}^\nu)_{x\in\Xt}, b^\nu, \sigma^\nu\big)$ be Lebesgue a.e.~martingale solution of  (\ref{SDE2}) satisfying 
\begin{align*}
  \int_{\XX}f(x)\nu_{t}(dx) &= \int_{\Xt}\int_\Om f\big(\phi^\nu_{t_0,t}(x,\om)\big)\Ns_{t_0,x}^\nu(d\om)\nu_{t_0}(dx) \qquad \forall \;t\in \Ic, \; f\in \mathbb{M}_\infty(\Xt).
\end{align*}
 Assume also that  $\varphi\in\mathcal{C}^2(\Rp)$  is strictly convex and  it satisfies the normality conditions~(\ref{Normality}). Then, there exists a unique $\varphi$-projection, $\Ms^{\mu\nu}_{t_{0}}\in \mathbb{C}_{\varphi,\Ic}^{{\mu\nu}}\subset\PP(\Wd)$ 
with a Markovian version, which satisfies the following variational identity
\begin{align}\label{InfEq}
&\Df\big(\Ms^{\mu\nu}_{t_{0}}\| \Ns_{t_{0}}^\nu\big)= \sup\bigg\{ \int_{\Xt}\sum_{i=1}^n f_i(x)\mu_{s_i}(dx) -\int_{\Xt}\varphi^*\Big(\sum_{i=1}^n f_i(x)\Big) \nu_{s_i}(dx);\notag\\ &\hspace{5cm} \forall \;1\leqslant n <\infty ,\;f_{1}, \dots, f_{n}\in\mathbb{M}_\infty(\Xt), \;   s_1,\dots,s_n\in \Ic\bigg\},
\end{align}
where the supremum  is over all $n$-tuple partitions of $\Ic$, and $n$-tuples of functions in $\mathbb{M}_\infty(\Xt)$\footnote{\,As noted earlier this could be extended to $f_i\in L_{\varphi*}(\XX;\nu_{s_i}),\; \forall \;1\leqslant n <\infty$.} for all $n\in \mathbb{N}$.
Moreover,  for finite-dimensional distributions, $\Ms^{\mu\nu}_{\Ic_{n}}, \,\Ns^\nu_{\Ic_{n}}\in \mathbb{C}_{\varphi,\Ic_{n}}^{{\mu\nu}}\subset \mathbb{C}_{\varphi,\Ic}^{{\mu\nu}}$, with some fixed partition $\Ic_n:=\{t_i\}_{t=0}^n\subset\Ic$, the following holds 
\begin{align}\label{Inphs2}
\Df\big(\Ms^{\mu\nu}_{\Ic_{n}}\| \Ns_{\Ic_{n}}^\nu\big) &\leqslant \sum_{i=1}^n\Big\vert \Df(\mu_{t_i}\|\mu_{t_0})-\Df(\nu_{t_i}\|\nu_{t_0})\Big\vert, 
\end{align} 
where $\Df(\mu_{t_i}\|\mu_{t_0})$ and $\Df(\nu_{t_i}\|\nu_{t_0})$ denote $\varphi$-FTDR fields associated with sequences $(\mu_{t_i})_{i=0}^n$, $(\nu_{t_i})_{i=0}^n$, of time-marginal probability measures  induced by (\ref{SDE1}) and (\ref{SDE2}) such that $\mu_{t_0} = \nu_{t_0}$.

\vspace{.3cm}\noindent{\it Proof.} See Appendix \ref{path_proj_app}; the proof is quite long and it utilises all of the preceding results in this section. 

\end{theorem}

\section{Case study}\label{tests}
Here, we illustrate the results derived in the previous sections, in particular \S\ref{Info_ineq} and  \S\ref{Inf_IP}, applied to some dimensionally reduced approximations of a simple slow-fast system. Detailed analysis of the impact of various (Eulerian) approximations of the governing dynamics on the accuracy of path-based (Lagrangian) observables  is postponed to subsequent publications. 

\vspace{-0.1cm}\subsection{Dimensional reduction of multi-scale dynamics via averaging}
First, we briefly outline a useful framework for deriving approximations of slow-fast dynamics which is obtained via averaging over a subset of dynamical variables representing the fast dynamics.

Consider a slow-fast SDE on $\R^{d}\times\R^{l}, \; d\geqslant 1,\; l\geqslant 1$, given for  $t\in [0,\,T]$ by 
\begin{equation}\label{Slow fast}
\begin{cases}
dX_t = \varepsilon \hspace{0.09em}b_\tx( t, X_t, Y_t)dt+\sqrt{\varepsilon}\hspace{0.09em}\sigma_\textsc{x}(X_t)dB_t,
\\[.2cm]
dY^{\varepsilon}_t = b_{\hspace{0.04em}\textsc{y}}(X_t, Y_t)dt+ \sigma_\textsc{y}(X_t, Y_t)dW_t, 
\end{cases} \quad (X_0,Y_0)\sim \mm_0\in \PP(\R^{d}\times\R^{l}),\quad 0<\varepsilon\ll 1,
\end{equation}
where the vector fields generating  (\ref{Slow fast}) satisfy   $b_{\hspace{0.04em}\textsc{x}}(t,\ccdot, \ccdot)\,{\in}\,\tilde{\mathcal{C}}^{\infty}(\Rd\,{\times}\,\R^l; \Rd)$,  $(\sigma_\textsc{x})_{m}\,{\in}\, \bar{\mathcal{C}}^{\infty}(\Rd;\Rd)$, $1\leqslant m\leqslant d$,  $b_{\hspace{0.05em}\textsc{y}}\,{\in}\,\, \tilde{\mathcal{C}}^{\infty}(\Rd\,{\times}\,\R^l; \R^l)$, $(\sigma_\textsc{y})_k \,{\in}\,\,\bar{\mathcal{C}}^{\infty}\big(\Rd\times\R^l; \R^l\big)$, $1\leqslant k\leqslant l$, and $B_t$, $W_t$ independent Brownian motions of appropriate dimension. In line with the notation in previous sections we set $\mm_0(dxdy) = \rrr_0(x,y)dxdy$, $\rrr_0(x,y)>0$. We refer to $X_t\in \Rd$ as the slow variable and to $Y_t\in \R^l$ as the fast variable.  Under the above regularity of the coefficients, the SDE~(\ref{Slow fast}) generates a stochastic flow of $\mathcal{C}^\infty$-\,diffeomorphisms (e.g., \cite{Arnold1, Kunitabook}). Moreover, if $\sigma_\textsc{y}(\sigma_\textsc{y})^*$ is also coercive and bounded, the `fast' system with the slow variable fixed, i.e., 
\begin{align}\label{fast_x}
d\tilde Y_t = b_{\hspace{0.05em}\textsc{y}}(x, \tilde Y_t)dt + \sigma_\textsc{y}(x,\tilde Y_t)dW_t,
\end{align}
admits ergodic measures $\big\{\Pi_x: x\in\mathbb{K}\Subset \Rd\big\}$ such that $x\mapsto \Pi_x$ is bounded Lipschitz continuous in the narrow topology generated by $\PP\big(\R^{l}\big).$

We note that the structure of the slow-fast system in (\ref{Slow fast}) is sufficient for the present illustration but it is relatively simple. In particular, the coefficients $(b_{\hspace{0.05em}\textsc{x}}, \sigma_\textsc{x})$, $(b_{\hspace{0.05em}\textsc{y}}, \sigma_\textsc{y})$ in (\ref{Slow fast}) do not depend on time or the time-scale separation parameter $\varepsilon$.  Often, the explicit dependence of these coefficients on $\varepsilon$  has to be considered in applications, and the associated averaging techniques have to be more sophisticated; for example, a framework for the analysis of convergence of the Heterogeneous Multiscale Methods (HMM; \cite{we03}) has to take into account the explicit dependence of the invariant measure of the fast dynamics (\ref{fast_x}) on the time-scale separation  $\varepsilon$ (e.g., \cite{liu10}).

 \smallskip
Given the simple slow-fast system (\ref{Slow fast}), define the averaged vector field $\bar{b}_\textsc{x}$ as 
\begin{align*}
\bar{b}_{\hspace{0.09em}\textsc{x}}(t,x)&= \int_{\R^l}{b}_{\hspace{0.05em}\textsc{x}}(t,x,y)\Pi_x(dy). 
\end{align*}
The `averaged' dynamics, representing an approximation of the evolution of the slow variables, is then given by  (see, e.g., \cite{Pavliotis,Bouchet16} among many others)
\begin{equation}\label{Average}
d \bar{X}_t = \bar{b}_{\hspace{0.05em}\textsc{x}}(t, \bar{X}_t)dt+\sigma_\textsc{x}(\bar{X}_t)dB_t, \qquad  \bar{X}_0 \sim \nu_0\in \PP(\R^d),
\end{equation}
where $X_{t/\varepsilon}\rightarrow \bar{X}_t$   as $\varepsilon\rightarrow 0$ in probability (an instance of weak LLN; e.g.,~\cite{Freidlin1, Kifer03}).~Then, as a consequence of ergodicity of the branches of the invariant measures $\{\Pi_x: x\in \mathbb{K}\Subset\Rd\},$ we have 
\begin{align*}
\bar{b}_{\hspace{0.05em}\textsc{x}}(t,x)= \int_{\R^l} b_{\hspace{0.05em}\textsc{x}}(t,x,y)\Pi_x(dy)= \lim_{T\rightarrow\infty}\frac{1}{T}\int_0^T b_{\hspace{0.05em}\textsc{x}}\big(t+s,x, \tilde Y_{s/\varepsilon}\big)ds.
\end{align*}

\noindent Another CLT-type approximation of ${X}_t$ can be obtained  by  accounting for the leading-order effects of fluctuations between the averaged dynamics and the slow dynamics, leading to (e.g.,~\cite{Bouchet16})
\begin{align*}
d \simbar X_t = \bar{b}_{\hspace{0.05em}\textsc{x}}\big(t, \simbar X_t\big)dt+\sigma_\textsc{x}\big(\simbar X_t\big)dB_t+\sqrt{\varepsilon}\sigma\big(t,\simbar X_t\big)d\hat{B}_t, \qquad  \simbar{X}_0 \sim \nu_0\in \PP(\R^d),
\end{align*}
where $B_t, \hat B_t$ are independent Brownian motions, and  the additional diffusion due to the fluctuations is defined via 
\begin{align*}
(\sigma\sigma^*)(t,x)&=\int_0^{\infty}C_{\tilde{F}}(t+s,x)ds,
\end{align*}
with the time-correlation matrix, $C_{\mathfrak{B}}(t,x)$, of $\mathfrak{B}(t,x,y) := b_{\hspace{0.05em}\textsc{x}}(t,x,y)-\bar{b}_{\hspace{0.05em}\textsc{x}}(t,x)$ defined by 
\begin{align*}
\notag C_{\mathfrak{B}}(t+s,x)= \int_{\R^{l}}\E\left[\mathfrak{B}(t+s, x, \tilde Y_s)\mathfrak{B}^{*}(t, x, y)+ \mathfrak{B}(t, x, y)\mathfrak{B}^{*}(t+s, x, \tilde Y_{s})\Big| \tilde  Y_0=y\right]\Pi_x(dy).
\end{align*}

With the above notation in place, we list the two approximations of the slow subsystem arising from stochastic averaging (e.g.,~\cite{Culina11, Arnold02, Bouchet16, Freidlin1, Pavliotis, Kifer03, Just02}) in Table~\ref{table: Red_equs}.
 \renewcommand{\arraystretch}{1.4}
 \begin{table}[H]
 \vspace{.0cm}
\begin{tabular}{ | c | c|}
 \hline
  \multicolumn{2}{| c |}{Approximations $\{\mathfrak{I}, \mathfrak{F}\}$ of the slow-fast SDE (\ref{Slow fast})} \\
  \cline{1-2}
   Approximation type/assumptions&  Dynamics\\
  \hline
 {\small Infinite time-scale separation: $\varepsilon \rightarrow 0$. $(\mathfrak{I})$}& $d\bar{X}_t = \bar{b}_{\hspace{0.05em}\textsc{x}}\big(t, \bar{X}_t\big)dt+\sigma_\textsc{x}\big(\bar{X}_t\big)dB_t$\\
  \hline
 {\small Finite time-scale separation: $0<\varepsilon \ll 1$, $(\mathfrak{F})$}& $d\simbar X_t = \bar{b}_{\hspace{0.05em}\textsc{x}}\big(t, \simbar X_t\big)dt+\sigma_\textsc{x}\big(\simbar X_t\big)dB_t+\sqrt{\varepsilon}\hspace{0.05em}\sigma\big(t,\simbar X_t\big)d\hat{B}_t$\\
  \hline
 
\end{tabular} 

\vspace{.2cm}\caption{Two reduced equations for slow-fast SDE \ref{Slow fast} considered in \S\ref{toysec}.} \label{table: Red_equs}
\end{table}

 \vspace{-1cm} \subsection{Information inequalities for approximated dynamics} Here, we outline how the various information bounds expounded in \S\ref{Inf_IP} and \S\ref{Info_ineq} can be utilised to assess the validity of dimensionally reduced approximations of slow-fast systems. In order to keep this paper reasonably concise, we shall only briefly describe such approaches for a toy example specified in \S\ref{toysec}; a full treatment of various types of multi-scale systems  is postponed  to subsequent publications.

\subsection{Toy example }\label{toysec}
Let $d=l=1,$ $t\in [0,\,T]$, and consider the following slow-fast SDE:
\begin{align}\label{Toy}
\begin{cases}
dX_t = \left(-\beta X_t+Y_t^2\right)dt+\sigma_\textsc{x} dB_t,\\[.3cm]
\displaystyle dY_t = -\frac{1}{\varepsilon}\gamma Y_tdt+\frac{1}{\sqrt{\varepsilon}}\sigma_\textsc{y}dW_t,
\end{cases} \quad (X_0,Y_0)\sim \mm_0\in \PP(\R^2).
\end{align}
Here $ 0<\varepsilon\ll1$ is the time-scale separation between the slow variable  $X_t\in \R$ and the fast variable $Y_t\in\R$, $\mm_0(dxdy) = \rrr^\mm_0(x,y)dxdy$, $\rrr^\mm_0(x,y)>0$, and $W_t, \; B_t\in \R$ are independent standard Brownian motions; the parameters in (\ref{Toy}) satisfy $\beta>0,\; \gamma>0,\; \sigma_\textsc{x}, \,\sigma_\textsc{y}\neq 0$. 

\smallskip
The particular structure of  (\ref{Toy}) generates a stochastic flow $\big\{\Phi_t^\mm\big((x,y),\ccdot\big), t\geqslant 0\big\}$ on $\R^2$, where 
$$\Phi_t^{\mm}\big((x,y),\om\big) = \Big(\phi_t^{\mm}\big((x,y),\om\big), \textrm{\textpsi}_t^{\mm}\big(y,\om\big)\Big),$$
 and $\textrm{\textpsi}^\mm_t(y,\om)$ is generated solely by the fast subsystem of (\ref{Toy}) via 
 \begin{align*}
\tilde Y_t(\om) = Y_t(\om) = e^{-\gamma t}y+\sigma_\textsc{y}\!\!\int_0^te^{-\gamma(t-s)}dW_s(\om):=\textrm{\textpsi}^\mm_t(y,\om).  
\end{align*} 
 Moreover, the fast dynamics in (\ref{Toy}) admits an ergodic measure $\Pi_x(dy) $ given by 
\begin{align*}
\Pi_x(dy) = \Pi_x(y)dy =\frac{\sqrt{\gamma}}{\pi\sigma_\textsc{y}}\exp\left(-\frac{y^2\gamma}{\sigma_\ty^2}\right)dy; 
\end{align*}
we use the same symbol for this ergodic measure and its density to simplify notation.

 Next, we calculate the coefficients of the reduced models listed in Table \ref{table: Red_equs}. The averaged equation ($(\mathfrak{I})$ in Table~\ref{table: Red_equs}) for the slow dynamics in (\ref{Toy}) from the weak LLN is given by 
\begin{align}
\notag d\bar{X}_t &= \bigg(-\beta\bar{X}_t+\frac{\sigma_\ty^2}{2\gamma}\bigg)dt+ \sigma_\tx dB_t=:\bar b_\tx^{\alpha}(\bar{X}_t)dt+\sigma_\tx dB_t, \qquad  \bar{X}_0 \sim \nu_0\in \PP(\R)
\end{align}
where $\alpha:=(\beta,\sigma_\ty)$, and $\nu_0(dx) = \rho^\nu_0(x)dx$, $\rho^\nu_0>0$.
In order to derive the approximation of (\ref{Toy}) from the CLT (see $(\mathfrak{F})$ in Table~\ref{table: Red_equs}), we consider the fluctuation between the averaged equation and the slow component of the the slow-fast SDE (\ref{Toy}),   $\sigma(x)$ can be evaluated as follows
\begin{align*}
\hspace{0cm} (\sigma\sigma^*)(x)& = \int_{0}^\infty C_{\mathfrak{B}}(s,x)ds =4\int_{0}^{\infty} \frac{\sigma^4_\ty}{4\gamma^2}e^{-2\gamma s} ds = \frac{\sigma^4_\ty}{2\gamma^3}.
\end{align*}
With the above calculations in place, we consider the following simple reduced models for (\ref{Toy}): 

 \vspace{-.1cm}
\begin{table}[H]
\begin{tabular}{ | c | c| c|c|}
 \hline
  \multicolumn{4}{| c |}{Approximations $\{\fI, \fF\}$ of the toy slow-fast SDE (\ref{Toy})} \\
  \cline{1-4}
   Approximation &  \hspace{3cm} Dynamics \qquad $\alpha = (\beta, \sigma_\ty)$ & Prob. measure & Flow \\
  \hline
  ($\fI$)& $d\bar{X}_t = \bar b^{\alpha}_\tx(\bar{X}_t)dt+ \sigma_\tx dB_t$ & $\nu_t^\fI\in \PP(\R)$, $\nu_0^\fI = \nu_0$ & $\phi^{\nu^{\fI}}_t(x,\om)$ \\
  \hline
   ($\fF$)& $d\simbar X_t = \bar b_\tx^{\alpha}\big(\simbar X_t\big)dt+\sigma_\tx dB_t+ \sqrt{\varepsilon}\frac{\;\;\scaleobj{1.2}{\sigma_\ty^2}}{\sqrt{2\gamma^3}}d\tilde{B}_t$ & $\nu_t^\fF\in \PP(\R)$, $\nu_0^\fF = \nu_0$ & $\phi^{\nu^\fF}_t(x,\om)$\\
  \hline
\end{tabular} 
\vspace{.1cm}\caption{Reduced equations for example \ref{Toy}. }\label{table: Red_Eqs}
\end{table}

\vspace{-.6cm}
\noindent Analogous to the notation adopted in (\ref{margs}), the probability measure $\mm_t\in \PP(\R^2)$ solving the forward Kolmogorov equation associated with (\ref{Toy}) and its $x$-marginal are denoted~by 
$$\mm_t(dxdy) = \rrr_t^{\mm}(x,y)dxdy, \qquad \mu_t(dx) = \rho^{\mu}_t(x)dx, \quad \rho^{\mu}_t(x) = \int_{\R^2}\rrr^{\mm}_t(x,y)dy, $$ 
while $\nu_t^\fI, \nu_t^\fF\in \PP(\R)$, $\nu_t^\fr(dx) = \rho^{\fr}_t(x)dx$, $\mathfrak{r}\in \{\mathfrak{I}, \mathfrak{F}\}$, denote the probability measures  of the respective reduced models listed in Table \ref{table: Red_Eqs}. 

Now, we can use the information inequality (\ref{Csiszar}) 
\begin{align}\label{obs_bnd}
\hat{\mathcal{K}}_{\varphi,\sff}^{\nu}\big(-\Df(\mu_t\|\nu^\fr_t)\big)\leqslant  \E^{\mu_t}[\ff] - \E^{\nu^\fr_t}[\ff] \leqslant \mathcal{K}_{\varphi,\sff}^{\nu}\big(\Df(\mu_t\|\nu^\fr_t)\big)
\end{align}
to investigate the reduced-order approximations in Table \ref{table: Red_Eqs} from the point of view of their effect on the error in  predictions of the Lagrangian observables (see (\ref{k_lift}) and Definition \ref{Lmartsol}) 
$$\E^{\mu_t}[f] = \int\E\big[f\big(\phi^{\mm}_t(x,y)\big)\big]\mm_0(dxdy), \quad \textrm{from}\quad  \E^{\nu^\fr_t}[f] = \int\E\big[f\big(\phi^{\nu^\fr}_t(x)\big)\big]\nu_0(dx), \quad \mathfrak{r}\in \{\mathfrak{I}, \mathfrak{F}\}.$$ 
  As discussed in \S\ref{Main_Bounds}, the bounds on $\Df\big(\mu_t\|\nu_t^{\fI}\big)$ and $\Df\big(\mu_t\|\nu_t^{\fF}\big)$ either can be considered either in terms of the analytically tractable bound (\ref{InformaBB}) in  Theorem \ref{Info_ineq2} (cf.~\S\ref{Inf_IP}),  or  can be  studied  via differences between the corresponding the $\varphi$-FTDR fields via the bound (\ref{ftdr_bnd}) in Theorem \ref{difference_bound} (cf.~\S\ref{ftdr_sec}).  Below, we utilise the bound (\ref{InformaBB}) in Theorem \ref{Info_ineq2};  the $\varphi$-FTDR will be valuable in computational  considerations which are postponed to future publications.

 \smallskip
 In order to derive the desired information bound  let  $\LG^{\mu,\ve*}$  be defined by 
\begin{align}\label{Lslow}
\LG^{\mu,\ve*}\rho^{\mu}_t(x) = \int_{\R}\!\!\!\LG^{\mm,\varepsilon*}\rrr_t^{\mm}(x,y)dy, \qquad \partial_t \rrr_t^{\mm} = \mathcal{L}^{\mm,\ve*}\rrr_t^{\mm},
\end{align}
where $\LG^{\mm,\ve*}$ is the $L^2$ dual of $\LG^{\mm,\ve} := \LG_x+y^2\partial_x+\varepsilon^{-1}\LG_y$ associated with (\ref{Toy}), where 
\begin{align*}
\LG_x &:= -\beta x\partial_x +\frac{\sigma_\tx^2}{2}\partial_{xx}^2, \qquad\text{and}\;\; \LG_y := -\gamma y\partial_y+\frac{\sigma^2_\ty}{2}\partial_{yy}^2\;.
\end{align*}
Let $\nu_t^\fr$ be the time-marginal probability measure associated with the reduced model for  $\mathfrak{r}\in \{\mathfrak{I}, \mathfrak{F}\}$.  We want to derive the following information bound based on Theorem~\ref{Info_ineq2}:
\begin{align}
\notag \Df(\mu_t\|\nu_{t}^\fr) &\leqslant \frac{1}{2}\int_0^t\int_{\R}\left\vert (\tilde\sigma^{\fr,-1}\varTheta^\varepsilon_{\mu\nu^\fr})(s,x)\right\vert^2\varphi^{\prime\prime}(\eta_s(x))\eta_s^2(x)\rho^\fr_s(x)dxds, 
\end{align}
where $\varTheta^\varepsilon_{\mu\nu^\fr}$ is defined as in (\ref{Per_vF}), $\eta_t = \rho^{\mu}_t/\rho_t^\fr$.

Here, for brevity, we focus on the above bound in terms of the KL-divergence which is obtained for  $\varphi(u) = u\log u-u+1, \; u>0$; noticing that in such a case $\varphi^{\prime\prime}(u) = u^{-1}$,  we have 
\begin{align}
\DK(\mu_t\|\nu_t^\fr) &\leqslant \frac{1}{2}\int_0^{t}\int_{\R}\Big\vert \tilde\sigma^{\fr,-1}\varTheta^\varepsilon_{\mu\nu^\fr}(s,x)\Big\vert^2\rho^{\mu}_s(x)dxds. 
\end{align}
The coefficients $\varTheta^\ve_{\mu\nu^\fr}, \; \fr\in \{\fI, \fF\}$ can be obtained by recalling the reconstructed forward Kolmogorov equation (\ref{RecEQ}) which, for the case of the dynamics in (\ref{Toy}), is given by 
\begin{align}\label{RecToy}
\partial_t \rho^{\mu}_t = \LG^{\mu,\ve*}\!\rho^{\mu}_t, \qquad  \LG^{\mu,\ve*}\rho^{\mu}_t=\LG^{\fr\hspace{0.01cm}*}\rho^{\mu}_t-\partial_x(\varTheta^\varepsilon_{\mu\nu^\fr}\rho^{\mu}_t).
\end{align}
This implies that 
\begin{align*}
\notag -\partial_x(\varTheta^\varepsilon_{\mu\nu^\fr}\rho^{\mu}_t) &= \LG^{\mu*}\rho^{\mu}_t(x)-\LG^{\fr\hspace{0.01cm}*}\rho^{\mu}_t(x) = \int_{\R}\LG^{\mm,\ve*}\rrr^{\mm}_t(x,y)dy - \LG^{\fr\hspace{0.01cm}*}\rho^{\mu}_t(x).
\end{align*}
{\bf Averaged SDE ($\fI$):}\\  If we write $\mathcal{L}^{\mm,\ve*} = \mathcal{L}_x^*-\partial_x(y^2\,\cdot\,)+\varepsilon^{-1}\mathcal{L}_y^*$, and recall from the stochastic averaging that 
\begin{align*}
\LG^{\fI*}\rho^{\mu}_t(x) &= -\partial_x\Big(\rho^{\mu}_t(x)\int_\R \big[-\beta x+y^2\,\big]\Pi_x(y)dy\Big)+\frac{\sigma_\tx^2}{2}\partial_{xx}^2\Big(\rho^{\mu}_t(x)\int_\R\Pi_x(y)dy\Big) \\ & = \int_{\R} \Big[\LG^*_x\big(\rho^{\mu}_t(x)\pi_x(y)\big)-\partial_x\big(y^2\rho^{\mu}_t(x)\Pi_x(y)\big)\Big]dy,
\end{align*}
we obtain 
\begin{align*}
\LG^{\mu,\ve*}\rho^{\mu}_t(x)-\LG^{\fI*}\rho^{\mu}_t(x) &= \LG^{*}_x\Big(\rho^{\mu}_t(x)\int_\R\left[\rrr^{\mm}_t(y|x)-\Pi_x(y)\right]dy\Big){-}\,\partial_x\Big(\rho^{\mu}_t(x)\int_\R y^2\left[\rrr^{\mm}_t(y|x)-\Pi_x(y)\right]dy\Big)\\ &\hspace{3cm}+\varepsilon^{-1}\int_{\R}\LG_y^*\rrr_t^{\mm}(x,y)dy\\
& = -\partial_x\Big(K^\varepsilon(t,x)\rho^{\mu}_t(x)\Big) +\frac{\sigma_\tx^2}{2}\partial^2_{xx}\Big(C_1^\varepsilon(t,x)\rho^{\mu}_t(x)\Big) +\varepsilon^{-1}\partial_x\Big(C_3^\varepsilon(t,x)\Big),
\end{align*}
where $C_i^\varepsilon(t,x), \; i=1,2,3$, and $K^\varepsilon(t,x)$ are defined by 
\begin{alignat*}{2}
K^\varepsilon(t,x)&= -\beta x C^\varepsilon_1(t,x)+C_2^\varepsilon(t,x), \hspace{1.6cm}&& 
C_1^\varepsilon(t,x) = \int_{\R}\Big[\rrr^{\mm}_t(y|x)-\Pi_x(y)\Big]dy,\\
C_2^\varepsilon(t,x) &= \int_{\R}y^2\Big[\rrr^{\mm}(y|x)-\Pi_x(y)\Big]dy, &&
C_3^\varepsilon(t,x) = \int_{c}^x\Big[\int_{\R}\LG_y^*\rrr_t^{\mm}(\xi,y)dy\Big]d\xi,
\end{alignat*}
for some constant $c\in \R.$ This implies that 
\begin{align*}
\varTheta^\ve_{\mu\nu^\fI}(t,x) = K^\varepsilon(t,x) -\frac{\sigma_\tx^2}{2\rho^{\mu}_t(x)}\partial_x\big(C_1^\varepsilon(t,x)\rho^{\mu}_t(x)\big)-\varepsilon^{-1} \frac{1}{\rho^{\mu}_t(x)}C_3^\varepsilon(t,x).
\end{align*}
To write the KL-divergence bound in this case, recall that $\tilde \sigma^{\fI,-1} = \sigma_\tx^{-1}, $ $\varphi(u) = u\log u-u+1$, $u>0$  and $\varphi^{\prime\prime}(u) = u^{-1}$. Then, for any $t\in [0, \,T],$ we have 
\begin{align}\label{InfAV}
\DK(\mu_t\|\nu_t^\fI) &\leqslant \frac{1}{2}\int_0^{t}\int_{\R}\Big\vert \sigma_\tx^{-1}\varTheta^\varepsilon_{\mu\nu^\fI}(s,x)\Big\vert^2\rho^{\mu}_s(x)dxds.
\end{align}
In order to compare the KL-divergence bound in (\ref{InfAV}) to that of other reduced models, we expand $\varTheta^\varepsilon_{\mu\nu^\fI}(t,x)$ in $\varepsilon>0$. To this end, we assume that $\rrr_t^{\mm}(x,y) = \rrr^0_t(x,y)+\varepsilon \rrr^1_t(x,y)+O(\varepsilon^2)$, where $\int\rrr^1_t(x,y)dxdy=0,$ and we obtain from (\ref{Lslow})
\begin{align*}
\mathcal{L}_y^*\rrr^0_t &= 0,\\
\partial_t\rrr^0_t &= \LG_x^*\rrr^0_t-\partial_x(y^2\rrr^0_t)+\LG^*_y\rrr_t^1,\\
\partial_t\rrr^1_t &= \LG_x^*\rrr^1_t-\partial_x(y^2\rrr^1_t),
\end{align*}
where the first equation implies that $\rrr_t^0(x,y) = \rho^\mu_t(x)\Pi_x(y).$
With the above expansion at hand, we can expand the functions $C_i^\varepsilon(t,x), \; \; i=1,2, 3,$ in $\varepsilon>0,$ as follows
\begin{align*}
C_1^\varepsilon(t,x) &= \varepsilon C_1^1(t,x)+O(\varepsilon^2)= \varepsilon\int_{\R}\rrr_t^1(y|x)dy+O(\varepsilon^2),\\
C_2^\varepsilon(t,x) &= \varepsilon C_2^1(t,x)+O(\varepsilon^2)= \varepsilon\int_{\R}y^2\rrr_t^1(y|x)dy+O(\varepsilon^2),\\
\varepsilon^{-1} C_3^\varepsilon(t,x) &= C_3^1(t,x)\,{+}\, \varepsilon C_3^2(t,x)\,{+}\,O(\varepsilon^2) \!= \!\!\int_c^x\!\Big[\LG^*_y\rrr_t^1(\xi,y)dy\Big]d\xi\,{+}\, \varepsilon\!\!\int_c^x\Big[\LG^*_y\rrr_t^2(\xi,y)dy\Big]d\xi\,{+}\, O(\varepsilon^2).
\end{align*}
Substituting into $\varTheta^\varepsilon_{\mu\nu^\fI}(t,x)$ and recalling that $K^\varepsilon(t,x) = -\beta x C_1^\varepsilon(t,x)+ C_2^\varepsilon(t,x),$ we have 
\begin{align*}
\varTheta^\varepsilon_{\mu\nu^\fI}(t,x) &= K^\varepsilon(t,x)-\frac{\sigma_\tx^2}{2}\frac{1}{\rho^{\mu}_t(x)}\partial_x\big(C_1^\varepsilon(t,x)\rho^{\mu}_t(x)\big)-\varepsilon^{-1}\frac{1}{\rho^{\mu}_t(x)}C^\varepsilon_3(t,x)\\ 
&= -\frac{1}{\rho^{\mu}_t(x)}C_3^1(t,x)-\varepsilon\beta xC_1^1(t,x)+\varepsilon C_2^1(t,x)\\ &\hspace{1cm}-\varepsilon\frac{\sigma_\tx^2}{2}\frac{1}{\rho^{\mu}_t(x)}\partial_x\big( C_1^1(t,x)\rho^{\mu}_t(x)\big)-\varepsilon\frac{1}{\rho^{\mu}_t(x)}C_3^2(t,x)+O(\varepsilon^2)\\
&= -\frac{1}{\rho^{\mu}_t(x)}C_3^1(t,x) + O(\varepsilon).
\end{align*}
Finally, we can write the KL-divergence bound (\ref{InfAV}) expanded in $\varepsilon>0$, for $ t\in [0, \,T]$,  as follows
\begin{align}\label{InfAV2}
\DK(\mu_t^\ve||\nu_t^\fI) &\leqslant \frac{1}{2}\int_0^{t}\int_{\R}\Big\vert \sigma_\tx^{-1}C_3^1(s,x)\frac{1}{\rho^{\mu}_s(x)}\Big\vert^2\rho^{\mu}_s(x)dxds + O(\varepsilon).
\end{align}

\noindent {\bf ``Fluctuating" approximation ($\fF$):} The generator of the dynamics associated with this approximation is (see Table \ref{table: Red_Eqs}) 
\begin{align*}
\LG^{\fF*}\rho^{\mu}_t(x) &= \LG^{\fI*}\rho^{\mu}_t(x)+\varepsilon\frac{\sigma_\ty^4}{4\gamma^3}\partial^2_{xx}\rho^{\mu}_t(x)
\end{align*}
and 
\begin{align*}
\notag -\partial_x(\varTheta^\varepsilon_{\mu\nu^\fF}(t,x)\rho^{\mu}_t(x)) &= \int\LG^{\mm,\ve*}\rrr^{\mm}_t(x,y)dy- \LG^{\fF*}\rho^{\mu}_t(x)\\
&= -\partial_x\Big(\varTheta_{\mu\nu^\fI}^\varepsilon(t,x)\rho^{\mu}_t(x)\Big) - \varepsilon\frac{\sigma_\ty^4}{4\gamma^3}\partial^2_{xx}\rho^{\mu}_t(x),
\end{align*}
so that 
\begin{align*}
\varTheta^\varepsilon_{\mu\nu^\fF}(t,x) &=\varTheta_{\mu\nu^\fI}^\varepsilon(t,x)+\frac{\varepsilon\sigma^4_\ty}{4\gamma^3}\partial_x\log\rho^{\mu}_t(x).
\end{align*}
Recalling that $\tilde\sigma^{\fF,-1} = \big(\sigma_\tx, \sqrt{\varepsilon}\frac{\sigma_\ty^2}{\sqrt{2\gamma}}\big)^{-1}$,  we have for all $t\in [0, \,T],$
\begin{align}\label{InfN}
\DK(\mu_t\|\nu_t^\fF)\leqslant \frac{1}{2}\int_0^t\int_{\R}\Big\vert  \Big(\sigma_\tx, \sqrt{\varepsilon}\frac{\sigma_\ty^2}{\sqrt{2\gamma^3}}\Big)^{-1}\varTheta^\varepsilon_{\mu\nu^\fF}(s,x)\Big\vert^2\rho^{\mu}_s(x)dxds.
\end{align}
Furthermore, we recall that 
\begin{align*}
\varTheta^\varepsilon_{\mu\nu^\fF}(t,x)&= \varTheta^\varepsilon_{\mu\nu^\fI}(t,x)+\varepsilon\frac{\sigma^4_\ty}{4\gamma^3}\partial_x\log\rho^{\mu}_t(x)= \frac{1}{\rho^{\mu}_t(x)}C_3^1(t,x)+O(\varepsilon).
\end{align*}
Since, $(\sigma_\tx, \sqrt{\varepsilon}\frac{\sigma_\ty^2}{\sqrt{2\gamma}})^{-1}= \Big(\sigma^2_\tx+ \varepsilon\frac{\sigma^4_\ty}{2\gamma^2}\Big)^{-1}\Big(\sigma_\tx, \sqrt{\varepsilon}\frac{\sigma_\ty^2}{\sqrt{2\gamma^3}}\Big)^*$, we have 
\begin{align*}
\Big\vert (\sigma_\tx, \sqrt{\varepsilon}\frac{\sigma_\ty^2}{\sqrt{2\gamma^3}})^{-1}\varTheta^\varepsilon_{\mu\nu^\fF}(t,x)\Big\vert^2&= \Big(\varTheta^\varepsilon_{\mu\nu^\fF}(t,x)\Big)^2\Big(\sigma_\tx^2+\varepsilon\frac{\sigma_\ty^4}{2\gamma^3}\Big)^{-1}\\
&=\Big(\frac{1}{\rho^{\mu}_t(x)}C_3^1(t,x)+O(\varepsilon)\Big)^2\Big(\sigma_\tx^2+\varepsilon\frac{\sigma_\ty^4}{2\gamma^3}\Big)^{-1}\\
&=\Big(\frac{1}{\rho^{\mu}_t(x)^2}\big(C_3^1(t,x)\big)^2+O(\varepsilon)\Big)\Big(\sigma_\tx^2+\varepsilon\frac{\sigma_\ty^4}{2\gamma^3}\Big)^{-1},
\end{align*}
 so that we have 
\begin{align}\label{1-x}
\Big(\sigma_\tx^2+\varepsilon\frac{\sigma_\ty^4}{2\gamma^3}\Big)^{-1}= \sigma_\tx^{-2}\Big(1+ \varepsilon\frac{\sigma_\ty^4}{2\sigma^2_\tx\gamma^3}\Big)^{-1} = \sigma_\tx^{-2}-\varepsilon\sigma_\tx^{-4}\frac{\sigma_\ty^4}{2\gamma^3}+ O(\varepsilon^2).
\end{align}
Thus,
\begin{align}\label{InfN0}
\Big\vert \Big(\sigma_\tx, \sqrt{\varepsilon}\frac{\sigma_\ty^2}{\sqrt{2\gamma^3}}\Big)^{-1}\varTheta^\varepsilon_{\mu\nu^\fF}(t,x)\Big\vert^2 = &\Big( \sigma_\tx^{-2}-\varepsilon\sigma_\tx^{-4}\frac{\sigma_\ty^4}{2\gamma^3}+O(\varepsilon^2)\Big)\Big(\frac{1}{\rho^{\mu}_t(x)^2}C_3^1(t,x)^2+O(\varepsilon)\Big)\notag \\[.1cm]
&= \sigma_\tx^{-2}\frac{1}{\rho^{\mu}_t(x)^2}C_3^1(t,x)^2\Big(1-\varepsilon\frac{\sigma_\ty^4}{2\sigma_\tx^2\gamma^3}\Big)+O(\varepsilon)\notag\\[.2cm]
&=\bigg\vert \sigma_\tx^{-1}C_3^1(t,x)\frac{1}{\rho^{\mu}_t(x)}\sqrt{1-\varepsilon\frac{\sigma_\ty^4}{2\sigma_\tx^2\gamma^3}}\bigg\vert^2+O(\varepsilon).
\end{align}
In this case, the KL-divergence bound (\ref{InfN}), can be written as 
\begin{align}\label{InfN1}
\DK(\mu_t\|\nu_t^\fF)\leqslant\int_0^{t}\int_{\R}\bigg\vert \sigma_\tx^{-1}C_3^1(s,x)\frac{1}{\rho^{\mu}_s(x)}\Big(1-\varepsilon\frac{\sigma_\ty^4}{2\sigma_\tx^2\gamma^3}\Big)^{1/2}\bigg\vert^2\rho^{\mu}_s(x)dxds+ O(\varepsilon).
\end{align}
We conclude that whenever  $\varepsilon\frac{\sigma_\ty^4}{2\sigma_\tx^2\gamma^3}\leqslant 1$, the probability measure $\nu_t^\fF$ represents a  better approximation\footnote{\,Note that, despite the fact that $\DK(\mu_t||\nu_t^\fF)-\DK(\mu_t||\nu_t^\fI)\sim \mathcal{O}(\varepsilon)$, the presence of an $\mathcal{O}(\varepsilon)$ term, $\varepsilon \,\sigma_\tx^{-1}C_3^1(s,x)/\rho_s(x)$, in the bound (\ref{InfN1}) on $\DK(\mu_t||\nu_t^\fF)$  directly affects the leading-order term in that bound. It can be easily checked that there is no equivalent  $\mathcal{O}(\varepsilon)$ term in the bound (\ref{InfAV2}) on $\DK(\mu_t||\nu_t^\fI)$. Given that $(1+x)^{-1}$ is analytic on $(0,\,\infty)$, there is no restriction imposed by the expansion in (\ref{1-x}) on $\delta_\varepsilon := \varepsilon \sigma_\ty^4/ (2\sigma_\tx\gamma^3)$ and, in principle, the leading-order term in (\ref{InfN1}) may be significantly smaller than that in (\ref{InfAV2}). While such effects appear in simple numerical simulations, a rigorous justification would require asserting that the leading-order terms in the $\varepsilon$-asymptotic expansions underlying the bounds in remain dominant for $\delta_\epsilon\sim\mathcal{O}(1)$. This is, however, beyond the scope of this illustrative example, and it is not necessary for outlining the  general scheme. 
  } of  $\mu_t$ in the sense that  the leading-order term of  (\ref{InfN1}) is smaller than that of (\ref{InfAV2}), which leads to a tighter bound on the error in estimating  path-based observables via the information inequality (\ref{obs_bnd}).  While the above result is not particularly surprising, it serves as a simple yet nontrivial illustration of the developed framework. Detailed analysis of path-based predictions from reduced-order approximations in multi-scale systems and the Lagrangian uncertainty quantification within our framework  are postponed to a subsequent~work.

\newpage
\section{Conclusions}\label{conclusions}

\noindent We developed a new framework for Lagrangian Uncertainty Quantification (LUQ) which is aimed at estimating and mitigating uncertainty in estimates of path-based observables evaluated on trajectories of a dynamical system representing an approximation of the original dynamics.  

\medskip
Specifically, given the paths $t\mapsto \phi^\mm_{t_0,t}(\xxt,\om)$ generated by the original dynamical system on a smooth finite-dimensional manifold $\XXt = \Xt\times \mathcal{Y}$, and paths $t\mapsto \phi^\nu_{t_0,t}(x,\om)$ generated by the approximating dynamics on $\Xt\subseteq\XXt$,  we obtained a hierarchy of bounds which relate the error in estimates of path-based observables to divergences between the underlying probability measures. For time-marginal probability measures $\mu_t, \nu_t\in \PP(\Xt)$,  these bounds have the form  
\begin{align*}
\hat{\mathcal{K}}_{\varphi,\sff}^{\nu}\big(-\Df(\mu_t\|\nu_t)\big)&\leqslant  \E^{\mu_t}[\ff] - \E^{\nu_t}[\ff] \leqslant \mathcal{K}_{\varphi,\sff}^{\nu}\big(\Df(\mu_t\|\nu_t)\big), \qquad  \,t\in \Ic = \big[t_0,\,t_0+T\big], \;\mu_{t_0} = \nu_{t_0}.
\end{align*}
These bounds  are tight and general; i.e., they are not restricted to Markovian processes or solutions of SDE's/ODE's. When $\{\phi^\mm_{t_0,t}, \;t\geqslant t_0\}$ and $\{\phi^\nu_{t_0,t}, \;t\geqslant t_0\}$ are induced by sufficiently regular SDE dynamics, generated by $(\brr^\mm,\sigma^\mm)$ and $(\brr^\nu,\sigma^\nu)$, we further have  (\S\ref{Info_ineq}, and \S\ref{Inf_IP})
\begin{align*}
\hat{\mathcal{K}}_{\varphi,\sff}^{\nu}\big(-\Df(\mu_t\|\nu_t)\big)&\leqslant  \E\big[f\big(\pi^\nu_\mm\circ\Solm^\mm\big)\big] - \E\big[f\big(\Solm^\nu\big)\big] \leqslant \mathcal{K}_{\varphi,\sff}^{\nu}\big(\Df(\mu_t\|\nu_t)\big), \\[.2cm]
&\Df(\mu_t\|\nu_t)\leqslant \Upsilon_{\varphi,t}^{\mm,\nu}\big(\,\brr^\mm,\brr^\nu, \sigma^\mm,\sigma^\nu\,\big),
\end{align*}
where  $\pi^\nu_\mm: \XXt\rightarrow \Xt$ is the natural  projection.  The bound on  $\Df(\mu_t\|\nu_t)$ in terms of the functional $\Upsilon_\varphi^{\mu\nu}$  (\ref{InformaBB})  provides an analytically tractable link between  Eulerian (field-based) error  and  uncertainties in Lagrangian~(path-based)~predictions. There bounds were extended to path space probability measures~$\Ps_{t_0}^\mm\in \PP(\mathcal{W}_\ell)$, $\Ps_{t_0}^\nu\in \PP(\mathcal{W}_d)$  induced by the laws of  
$\phi^\mm_{t_0,\cccdot}$, $\phi^\nu_{t_0,\cccdot}$ in \S\ref{Path-space}.

 \smallskip
 Moreover, we derived another general bound on the $\varphi$-divergence in the form   (\S\ref{ftdr_sec})
 $$\hspace{1.3cm}\Df(\mu_t\|\nu_t)\leqslant \big\vert \Df( \mu_t\|\mu_{t_0})-\Df(\nu_t\|\nu_{t_0})\big\vert, \hspace{1.3cm}  t\in [t_0,\,t_0+T], \;\;\mu_{t_0} = \nu_{t_0}.$$
  This bound applies to both stochastic and deterministic dynamics, and it can be cast in terms of  differences between {\it finite-time divergence rate} ($\varphi$-FTDR) fields  which utilise local expansion rates in stochastic flows \cite{MBUda18}.  Importantly, such bounds  can be exploited within a computational framework to mitigate the error in Lagrangian predictions by tuning the fields of expansion rates. We also showed that the above bounds generalise to probability measures on path spaces;  while such  bounds are more difficult to deal with in practice, they take into account much more information about spatio-temporal correlations than the bounds based on families of one-point marginals, $(\mu_t)_{t\in \Ic}$, $(\nu_t)_{t\in \Ic}$; we will exploit this approach in future~work.

\smallskip
Another  strand of a follow-up research involves uncertainty quantification and optimal  path space tuning of classes of dynamical models which are of practical interest in studies focused on path-based  evolution -- specifically, in the context of transport and mixing in dynamical systems and in oceanographic applications utilising  either analytically simplified or data-driven models.

\smallskip
An explicit use of  the abstract geometry imposed by $\varphi$-divergences on the space of probability measures to analyse the approximating capability of specific classes of dynamical systems for optimal path-based predictions is a subject of an ongoing work. This approach allows  for information-geometric study of   statistical estimation and inference, large deviations, and it is also useful  in the analysis of learning efficiency and robustness of  neural network estimates, or the path-space analysis of data assimilation~techniques.  These issues will be addressed in a separate work. 
 
\newpage
\begin{appendices}
\section*{Appendix. Further proofs}\label{f_prf}
\addtocontents{toc}{\protect\setcounter{tocdepth}{1}}
\renewcommand{\thesubsection}{\Alph{subsection}}
\renewcommand{\theequation}{\Alph{subsection}.\arabic{equation}}
\setcounter{section}{1}

\subsection{Proof of Theorem \ref{iNF_IN}}\label{app_thm_iNF_IN}
 We present the proof of Theorem \ref{iNF_IN} in two propositions, Proposition {\ref{BU1}} and \ref{BU2},  corresponding to (\ref{Goal}) and (\ref{Goal2}) respectively.

\begin{prop}\label{BU1} \rm Let $\mu$ and $\nu$ be probability measures on a complete separable metric  (Polish) space $(\XX,\Bb(\XX))$ such that  $\Df(\mu\|\nu)<\infty$, $\mu,\nu\in \PP(\XX)$. Then, for any $f\in  L_{\varphi*}(\XX,\nu)$ (\ref{Orlicz}), we have   
\begin{align*}
\mathfrak{B}_{\varphi,-}(\mu\|\nu;f)\leqslant \E^{\mu}[f] - \E^{\nu}[f]\leqslant \mathfrak{B}_{\varphi,+}(\mu\|\nu;f).
\end{align*}

\end{prop}
\noindent {\it Proof.} We note that $\Df(\mu\|\nu)<\infty$ implies that there exists $0<\eta = (d\mu/d\gamma)/(d\nu/d\gamma)\in L^1(\XX;\nu)$  such that $d\mu = \eta d\nu$ (e.g.,~\citep{Gmeas}). Then, for any $f\in  L_{\varphi*}(\XX,\nu)$, the Fenchel-Young inequality $  \eta f \leqslant \varphi^*(f)+\varphi\big(\eta\big) $ implies that 
\begin{align}\label{466}
\int_{\Xt} f d\mu - \Df(\mu\|\nu)\leqslant \int_{\Xt} \varphi^*(f)d\nu<\infty.
\end{align}

 Clearly,   $\lambda\big(f-\E^{\nu}[f]\big)\in  L_{\varphi*}(\XX,\nu)$ for  $\lambda \in \Rp$ and  for any $f\in  L_{\varphi*}(\XX,\nu)$; thus we have 
\begin{align*}
\lambda\big( \E^{\mu}[f] - \E^{\nu}[f]\big)-\Df(\mu\|\nu)\leqslant \int_{\Xt} \varphi^*\big(\lambda \big(f-\E^{\nu}[f]\big)\big)d\nu
\end{align*}
 which leads to 
\begin{align}\label{4.2}
 \E^{\mu}[f] - \E^{\nu}[f]\leqslant \frac{1}{\lambda}\bigg (\int_{\Xt} \varphi^*\big(\lambda \big(f-\E^{\nu}[f]\big)\big)d\nu + \Df(\mu\|\nu)\bigg).
\end{align}
For $\lambda>0,$ the function $\lambda (\E^{\nu}[f] -f)\in  L_{\varphi*}(\XX,\nu)$ for $f\in  L_{\varphi*}(\XX,\nu)$ with $\Df(\mu\|\nu)<\infty,$ and we also obtain
 \begin{align*}
-\lambda \big ( \E^{\mu}[f] - \E^{\nu}[f]\big) \leqslant \int_{\Xt} \varphi^*\big(\lambda (\E^{\nu}[f]-f)\big)d\nu+\Df(\mu\|\nu).
\end{align*}
Finally, for any $f\in L_{\varphi*}(\XX,\nu),$ we have 
\begin{align}\label{4.3}
\E^{\mu}[f]-\E^{\nu}[f]\geqslant -\frac{1}{\lambda}\bigg(\int_{\Xt} \varphi^*\big(\lambda (\E^{\nu}[f]-f)\big)d\nu+\Df(\mu\|\nu)\bigg).
\end{align}
Combining (\ref{4.2}) and (\ref{4.3}), and taking the infimum over $\lambda>0$ leads to the desired formula
\begin{align*}
\hspace{3cm} \mathfrak{B}_{\varphi,-}(\mu\|\nu;f)\leqslant \E^{\mu}[f] - \E^{\nu}[f]\leqslant \mathfrak{B}_{\varphi,+}(\mu\|\nu;f). \hspace{4.5cm}\qed
\end{align*}

\begin{prop}\label{BU2} \rm\;
Let $\mu$ and $\nu$ be probability measures on the Polish space $(\XX,\Bb(\XX))$ such that $\Df(\mu\|\nu)<\infty$. Assume that the strictly convex function $\varphi \,{:} \;\Rp{\rightarrow}\,\R$ satisfies (\ref{Normality}) and is twice continuously differentiable. Then, the functionals $\mathfrak{B}_{\varphi, \pm}(\mu\|\nu; f)$ have the following properties:
\btm
\item[(i)] $\mathfrak{B}_{\varphi, +}(\mu\|\nu; f)\geqslant 0$ and $\mathfrak{B}_{\varphi, +}(\mu\|\nu; f) = 0,$ iff $\mu=\nu$ or $f$ is constant $\nu\,\textrm{-a.e.,}$
\item[(ii)] $\mathfrak{B}_{\varphi, -}(\mu\|\nu; f)\leqslant 0$ and $\mathfrak{B}_{\varphi, -}(\mu\|\nu; f) =0,$ iff $\mu=\nu$ or $f$ is constant $\nu\,\textrm{-a.e.}$
\etm
\end{prop}
\noindent{\it Proof.}
We only prove part (i) of the proposition since part (ii) can be derived  by changing the sign of $\lambda$. First, consider the strictly convex function $\lambda\mapsto \mathcal{G}_{\varphi,\nu}(\lambda; f)$ defined by 
\begin{align}\label{GG}
\mathcal{G}_{\varphi,\nu}(\lambda; f) = \int_{\Xt}\varphi^*\Big(\lambda(\E^{\nu}[f] -f)\Big)d\nu.
\end{align}
 Then, the functional  $\mathfrak{B}_{\varphi, +}(\mu\|\nu; f)$ becomes 
\begin{align*}
 \mathfrak{B}_{\varphi, +}(\mu\|\nu; f) = \inf_{\lambda>0}\bigg\{ \frac{1}{\lambda}\mathcal{G}_{\varphi,\nu}(\lambda; f)+ \frac{1}{\lambda}\Df(\mu\|\nu)\bigg\}.
\end{align*}
Observe that $\Df(\mu\|\nu)\geqslant 0$ by definition; thus, it only remains to show  that $\mathcal{G}_{\varphi,\nu}(\lambda; f) \geqslant 0$. By Jensen's inequality together with the normality conditions (\ref{Normality}), we have 
\begin{align*}
\mathcal{G}_{\varphi,\nu}(\lambda; f) =  \int_{\Xt}\varphi^*\Big(\lambda(\E^{\nu}[f] -f)\Big)d\nu\geqslant \varphi^*\bigg(\int_{\Xt}\lambda(\E^{\nu}[f] -f)d\nu\bigg) = \varphi^*(0) = 0,
\end{align*}  
which implies that $ \mathfrak{B}_{\varphi, +}(\mu\|\nu; f)\geqslant 0$.

\smallskip
Next, if $\mu = \nu$, then $\Df(\mu\|\nu) = 0$  and we have 
 \begin{align}\label{reff}
\mathfrak{B}_{\varphi, +}(\mu\|\nu; f) = \inf_{\lambda>0}\Big\{\frac{1}{\lambda} \mathcal{G}_{\varphi,\nu}(\lambda; f) \Big\} = \lim_{\lambda \rightarrow 0}\frac{1}{\lambda}\mathcal{G}_{\varphi,\nu}(\lambda; f) = \frac{d}{d\lambda}\mathcal{G}_{\varphi,\nu}(\lambda; f) \big |_{\lambda =0} \notag\\= \int_{\Xt}\nabla\varphi^{*}(0)(\E^{\nu}[f]-f)d\nu = 0,
\end{align}
since  $\varphi^*$ is strictly convex and $\varphi^*(0) =0.$

Conversely, assume $  \mathfrak{B}_{\varphi, +}(\mu\|\nu; f)=0$ and $f\ne \E^\nu[f]$.  
First, recall that 
\begin{equation}\label{infG_infD}
\inf_{\lambda>0}\Big\{\frac{1}{\lambda} \mathcal{G}_{\varphi,\nu}(\lambda; f) \Big\}=0, \qquad \inf_{\lambda>0}\Big\{\frac{1}{\lambda} \Df(\mu\|\nu) \Big\}=0.
\end{equation}
Thus, based on the properties of $\mathcal{G}_{\varphi,\nu}(\lambda; f)$, and $\Df(\mu\|\nu)$, and the properties of the infimum, the constraint  
\begin{align}\label{inf0}
0 =  \mathfrak{B}_{\varphi, +}(\mu\|\nu; f)& =  \inf_{\lambda>0}\bigg\{ \frac{1}{\lambda}\mathcal{G}_{\varphi,\nu}(\lambda; f)+ \frac{1}{\lambda}\Df(\mu\|\nu)\bigg\}
\end{align}
implies that  $\mu=\nu$  when $f\ne \E^\nu[f]$.  If, on the other hand,  $f = \E^{\nu}[f] \;\nu\textrm{-a.s.}$, then $\mathcal{G}_{\varphi,\nu}(\lambda;f) = 0$, since $\varphi^*(0) \,{=}\,0$ by (\ref{Phistar}) and (\ref{Normality}), and one arrives at (since $\Df(\mu\|\nu)\geqslant 0$)
\begin{align*}
\mathfrak{B}_{\varphi, +}(\mu\|\nu; f) = \inf_{\lambda>0}\Big\{\frac{1}{\lambda}\Df(\mu\|\nu)\Big\} = 0.
\end{align*}
Finally, if  $\mathfrak{B}_{\varphi, +}(\mu\|\nu; f) =0$ and $\mu\ne \nu$, then (\ref{infG_infD}) and (\ref{inf0})  imply that $\E^\nu[f]=f$.

\subsection{Proof of Proposition \ref{InformaRe}}\label{InformaRe_ap}

The representation formula for the bound $\mathfrak{B}_{\varphi,\pm}(\mu\|\nu;f)$ with  $\mu, \nu \,{\in}\, \PP(\XX) $ and $f\,{\in} \,L_{\varphi*}(\XX,\nu)$, requires the notion of a {\it pseudo-inverse} of a real-valued function. 

\vspace{-0.1cm}\bn[Pseudo-inverse]\rm\label{Pseudof}
For a nondecreasing function, $\eta: \R\rightarrow [0, \infty],$ set
\begin{align*}
\mathcal{J} =\begin{cases} \big[\inf_{x}\eta(x),\;\,\sup_{x}\eta(x)\big], \;\text{ if }\; \eta \;\text{is bounded},\\[.1cm]
 \big[\inf_{x}\eta(x), \;\,\infty\big),\hspace{1.8cm}\text{otherwise}.
 \end{cases}
 \end{align*}
The function $\tilde{\eta}^{-1}: \mathcal{J}\rightarrow [0, \infty]$ is the {\it pseudo-inverse} of $\eta$ and is given by  
\begin{align*}
 \tilde{\eta}^{-1}(y):= \inf\{x: \eta(x)\geqslant y\}, \quad\forall \;y\in \mathcal{J}.
\end{align*}
The pseudo-inverse $\tilde{\eta}^{-1}(y)$ is uniquely determined almost everywhere on $\mathcal{J}.$
\en
\noindent \underline{\it Part 1 of proof of Proposition \ref{InformaRe}:}
Define $$\Theta_{+}(\lambda,R) := \frac{1}{\lambda}\mathcal{G}_{\varphi,\nu}(\lambda;f)+\frac{1}{\lambda}R^2, \qquad  \Theta_{-}(\lambda,R) := -\frac{1}{\lambda}\mathcal{G}_{\varphi,\nu}(-\lambda;f)-\frac{1}{\lambda}R^2, $$  for $\lambda>0$,  where $R^2 = \Df(\mu\|\nu)$.  Then 
\begin{equation}\label{Minimize}
\mathfrak{B}_{\varphi,+}(\mu\|\nu;f) = \inf_{\lambda> 0}\Theta_{+}(\lambda,R), \quad \mathfrak{B}_{\varphi,-}(\mu\|\nu;f) = \sup_{\lambda > 0}\Theta_{-}(\lambda,R).
\end{equation}
Below, we prove the representation formula for $\mathfrak{B}_{\varphi,+}(\mu\|\nu;f)$; the formula for $\mathfrak{B}_{\varphi,-}(\mu\|\nu;f)$ is obtained in an analogous fashion by replacing the sign of $\lambda$ in the steps below.
First, notice that $\lambda\mapsto \mathcal{G}_{\varphi,\nu}(\lambda;f)$ is convex and $\mathcal{G}_{\varphi,\nu}(0;f) =0.$ It then follows that $\delta\mapsto\mathcal{G}_{\varphi,\nu}^*(\delta;f)$ is nonnegative and nondecreasing convex function on $[0, \infty]$ with $\mathcal{G}_{\varphi,\nu}^*(0;f) =0$.
We make the following observations: 
\btm 
\item[(i)]  From the proof of Proposition \ref{BU2} the strict convexity of $\mathcal{G}_\nu(\ccdot\,;f)$ for $\E^\nu[f]\neq f, \; \nu\textrm{\,-\,a.s.},$ follows from the strict convexity of $\varphi^*$.

\item[(ii)] Part of the normality conditions (\ref{Normality}) read, $\varphi(1) =0$ and $\varphi(u)>-\infty$ for all $u\in \Rp$; this implies that $\text{dom}\;\varphi^*\neq\emptyset$ and $\varphi^*(\xi)>-\infty$ for all $\xi\in\Rp.$ 
\etm
The observations (i) and (ii) imply that  $\mathcal{G}_{\varphi,\nu}(\ccdot\,;f)$ is a proper convex function. Next, if $R> 0,$ we have seen from the proof of Proposition \ref{BU2}, that $\Theta_{+}(\lambda;R)\rightarrow \infty$ as $\lambda\downarrow 0$ and as $\lambda\uparrow \infty.$ This implies that the infimum is achieved at some point $\lambda^\dagger\in (0, \infty).$
Suppose that the infimum is not unique; say, there exists an infimum $L>0$  such that for $0<\lambda_1^\dagger<\lambda_2^\dagger<\infty,$ we have 
\begin{align*}
\mathcal{G}_{\varphi,\nu}(\lambda_1^\dagger;f)+R^2 = \lambda_1^\dagger L \;\;\text{and}\;\; \mathcal{G}_{\varphi,\nu}(\lambda_2^\dagger;f)+R^2 = \lambda_2^\dagger L.
\end{align*}
 Given that  the function $\mathcal{G}_{\varphi,\nu}(\ccdot;f)$ is a proper convex function, we have that for $\tilde{\lambda} = \frac{1}{2}(\lambda_1^\dagger+\lambda_2^\dagger),$
\begin{align*}
\mathcal{G}_{\varphi,\nu}(\tilde{\lambda};f)+R^2<\tilde{\lambda}L,
\end{align*}
which contradicts the minimality of $\lambda_1^\dagger$ and $\lambda_2^\dagger.$ Thus, the minimiser of $\Theta_{+}(\lambda,R)$ is unique and finite for $R> 0$.  
For $R=0$,  we get a continuous extension\footnote{\,This extension also follows from the normality conditions (\ref{Normality}), which ensures the convex function $\varphi$ is non-increasing on the interval $(0,1]$ and nondecreasing on the extended interval $[1,\infty].$ Thus, its Legendre-Fenchel conjugate $\varphi^*$ is nondecreasing on $[0, \infty]$ with $\varphi^*(0) =0.$  } of $\Theta_{+}(\lambda,0)$ such that $$\Theta_{+}(\lambda,0)|_{\lambda=0} = \Theta_{+}(0,0) =0,$$ since $\mathcal{G}_{\varphi,\nu}(0;f) = \nabla_\lambda\mathcal{G}_{\varphi,\nu}(0;f) =0$, 
and we obtain a unique minimiser $\lambda^\dagger(0) =0$. By the lower semicontinuity of $(\lambda,R)\mapsto\Theta_{+}(\lambda,R)$, we extend the minimisation problem to (\ref{Minimize}) for all $R\geqslant 0$. 
Next, we observe that the Fenchel-Young inequality yields
\begin{align*}
\mathcal{G}_{\varphi,\nu}^*(\delta;f)\geqslant \delta\lambda - \mathcal{G}_{\varphi,\nu}(\lambda;f), \;\; \forall \lambda\in \mathbb{R},
\end{align*}
 which implies that $\mathcal{G}_{\varphi,\nu}^*(\delta;f)$ is strictly convex, nonnegative and unbounded. Thus, the inverse\footnote{\,In this case, the pseudo-inverse coincides with the inverse.} $\tilde{\mathcal{G}}^{*^{-1}}_{\varphi,\nu}(t;f)\,{=:}\,\mathcal{K}_{\varphi,f}^\nu(t)$ exists for all $t\geqslant 0.$ Recalling the definition of $\mathfrak{B}_{\varphi,+},$ we have for $R> 0$
\begin{align}\label{thdagg}
\mathfrak{B}_{\varphi,+}(\mu\|\nu;f) = \inf_{\lambda>0}\left\{\frac{1}{\lambda}\mathcal{G}_{\varphi,\nu}(\lambda;f)+\frac{1}{\lambda}R^2\right\}=\frac{1}{\lambda^\dagger(R)}\mathcal{G}_{\varphi,\nu}(\lambda^\dagger(R);f)+\frac{1}{\lambda^\dagger(R)}R^2=:\theta^\dagger(R),
\end{align}
where $0<\lambda^\dagger(R)<\infty$ is unique for $R\ne 0$, and  $\mathfrak{B}_{\varphi,+}(\mu\|\nu;f)=0$ for $R=0$ given the properties discussed above.  It remains to show that 
\begin{align*}
\theta^\dagger(R) = \mathcal{K}_{\varphi,f}^\nu(R) = \inf\{\delta\geqslant 0: \mathcal{G}^*_{\varphi,\nu}(\delta;f)>R^2\},
\end{align*}
which is equivalent to  showing that $\mathcal{G}_{\varphi,\nu}^*(\delta;f)>R^2$ iff $\delta>\theta^\dagger$. To see this, consider first $R> 0$ so that $0<\lambda^\dagger(R)<\infty$. 
If $\delta>\theta^\dagger$, the Fenchel-Young inequality and (\ref{thdagg}) lead to  
\begin{align*}
\mathcal{G}_{\varphi,\nu}^*(\delta;f) +\mathcal{G}_{\varphi,\nu}(\lambda^\dagger;f)&\geqslant \delta\lambda^\dagger > \theta^\dagger \lambda^\dagger= \left(\frac{1}{\lambda^\dagger}\mathcal{G}_{\varphi,\nu}(\lambda^\dagger;f)+\frac{1}{\lambda^\dagger}R^2\right)\lambda^\dagger = \mathcal{G}_{\varphi,\nu}(\lambda^\dagger;f)+R^2,
\end{align*}
which yields the result. For $R=0$, we have $\lambda^\dagger=0$ and the Fenchel-Young inequality leads to 
$$\mathcal{G}_{\varphi,\nu}^*(\delta;f) +\mathcal{G}_{\varphi,\nu}(\lambda^\dagger;f)\geqslant \delta\lambda^\dagger\quad \Longrightarrow \quad \mathcal{G}_{\varphi,\nu}^*(\delta;f) \geqslant 0$$
since $\mathcal{G}_{\varphi,\nu}(0;f)=0$. 
If $\mathcal{G}^*_{\varphi,\nu}(\delta;f)>R^2$, let $\alpha(\lambda) := \lambda \delta-\lambda\theta^\dagger$, $\lambda>0$, and note that 
\begin{align*}
\alpha(\lambda) = \delta\lambda -\lambda\inf_{\lambda>0}\Big\{\frac{1}{\lambda}\mathcal{G}_\nu(\lambda;f)+\frac{1}{\lambda}R^2\Big\}\geqslant \delta\lambda  -  \mathcal{G}_{\varphi,\nu}(\lambda;g)-R^2
>\delta\lambda  -\mathcal{G}_{\varphi,\nu}(\lambda;f)-\mathcal{G}^*_{\varphi,\nu}(\delta;f),
\end{align*}
which implies 
\begin{align*}
\lambda \delta-\mathcal{G}_{\varphi,\nu}(\lambda;f)< \mathcal{G}^*_{\varphi,\nu}(\delta;f)+\alpha(\lambda).
\end{align*}
Taking the supremum of both sides over $\lambda>0$ and recalling the definition of Fenchel's convex conjugate $\mathcal{G}^*_{\varphi,\nu}(\delta;f),$ we have $\sup\{\alpha(\lambda): \lambda>0\}>0.$ Since $\alpha(\lambda)$ is a linear function with $\alpha(0) =0$ and  $\sup\{\alpha(\lambda): \lambda>0\}>0,$ we have that 
$\delta\lambda -\theta^\dagger\lambda = \alpha(\lambda)>0$ for all $\lambda>0$ which implies that $\delta>\theta^\dagger$.  Consequently, we arrive at $\inf\{\delta\geqslant 0: \mathcal{G}_{\varphi,\nu}^*(\delta;f)>R^2\} = \theta^\dagger.$
\qed

\medskip
\noindent \underline{\it Part 2 of proof of Proposition \ref{InformaRe}:} 
Similar to part 1 of the proof we focus on $\mathfrak{B}_{\varphi,+}(\mu\|\nu;f)$; the representation formula for $\mathfrak{B}_{\varphi,-}(\mu\|\nu;f)$ is obtained in an analogous fashion by replacing the sign of $\lambda$ in the steps below.
Given part 1 of the proof,   $\mathcal{G}_{\varphi,\nu}(\ccdot\,;f)$ is a proper convex function and twice continuously differentiable on its effective domain $|\lambda|<\infty$. For $0\,{<}\, R\,{<}\, \infty$ the unique minimiser $0\,{<}\,  \lambda^\dagger=\lambda^\dagger(R)\,{<}\, \infty$ of $\Theta_+(\ccdot\,,R)$ satisfies
 \begin{align}\label{Optimality}
-\frac{1}{(\lambda^\dagger)^2}\mathcal{G}_{\varphi,\nu}(\lambda^\dagger;f)+\frac{1}{\lambda^\dagger}\nabla_\lambda\mathcal{G}_{\varphi,\nu}(\lambda^\dagger;f)-\frac{1}{(\lambda^\dagger)^2}R^2 =0, 
\end{align}
or equivalently
\begin{equation}\label{Uniq_minimum}
-\mathcal{G}_{\varphi,\nu}(\lambda^\dagger;f)+\lambda^\dagger\nabla_\lambda\mathcal{G}_{\varphi,\nu}(\lambda^\dagger;f) = R^2.
\end{equation}
Now, let $\mathcal{G}^*_{\varphi,\nu}(\delta;f)$ be the Legendre-Fenchel conjugate of $\mathcal{G}_{\varphi,\nu}(\ccdot\,;f)$ for $\lambda> 0$; i.e.,
\begin{align}\label{GG*}
\mathcal{G}_{\varphi,\nu}^*(\delta;f) = \sup_{\lambda> 0}\big\{\lambda\delta-\mathcal{G}_{\varphi,\nu}(\lambda;f)\big\}.
\end{align}
By the lower semicontinuity and proper convexity of $\mathcal{G}_{\varphi,\nu}(\lambda;f)$ with $\mathcal{G}_{\varphi,\nu}(\lambda;f)\,{<}\, \infty$\footnote{\,Given that $f\,{\in} \,L_{\varphi*}(\XX,\nu)$ we immediately have $\mathcal{G}_{\varphi,\nu}(\lambda;f)\,{<}\, \infty$; see (\ref{Orlicz}).} in an open neighbourhood of $\lambda\,{=}\,0$  together with $\nabla_\lambda\mathcal{G}_{\varphi,\nu}(0;f) =0,$ the sub-differential $\partial \mathcal{G}_{\varphi,\nu}(s;f)\neq \emptyset$ \footnote{\,The set $\partial\hspace{.02cm}\mathcal{G}_{\varphi,\nu}(s;f)$ is defined by $\partial\hspace{.02cm}\mathcal{G}_{\varphi,\nu}(s;f) :=\{t\in \R: \mathcal{G}_{\varphi,\nu}(\lambda;f) - \mathcal{G}_{\varphi,\nu}(s;f)\geqslant t(\lambda -s)\rangle,\;\;\ \forall \lambda\in \R\}$.} for $s$ in an open neighbourhood of $\lambda\,{=}\,0$.  Importantly, $0\in \partial\mathcal{G}_{\varphi,\nu}(s;f)$ which  implies that $\mathcal{G}^*_{\varphi,\nu}(\delta;f)$ has a unique minimum at $\delta=0$, and $\mathcal{G}_{\varphi,\nu}^*(\delta;f)\rightarrow \infty$ as $\delta\rightarrow \infty$ (e.g., \cite{Rockafeller74}; see also~\cite{dupuis16} for the case of KL-divergence). For $\lambda>0$ let $\hat\delta(\lambda) := \nabla_\lambda\mathcal{G}_{\varphi,\nu}(\lambda;f) > 0$ be the unique solution of $\mathcal{G}^*_{\varphi,\nu}(\hat\delta;f) =\lambda\hat\delta-\mathcal{G}_{\varphi,\nu}(\lambda;g)$; then, by convex duality, we have that 
\begin{align}\label{HHg}
\mathcal{H}_g(\lambda) = -\mathcal{G}_{\varphi,\nu}(\lambda;f)+\lambda\nabla_\lambda\mathcal{G}_{\varphi,\nu}(\lambda;f) = \mathcal{G}^*_{\varphi,\nu}(\hat\delta(\lambda);f) 
\end{align}
is non-negative and strictly increasing for $\lambda> 0$ (this is due to the fact that $\mathcal{G}^*_{\varphi,\nu}(\delta;f)$ in (\ref{GG*}) is strictly increasing for $\delta> 0$, and $\hat\delta(\lambda)$ is strictly increasing since $\nabla_\lambda^2\mathcal{G}_{\varphi,\nu}(\lambda;f) >0$).
Thus, the pseudo-inverse $\tilde{\mathcal{H}}^{-1}_{+,f}$ of $\mathcal{H}_f$ defined on $(0,\infty)$ exists uniquely almost everywhere for $\lambda> 0$. Then, from (\ref{Uniq_minimum}), we have for $0<R<\infty$ that 
\begin{align}\label{Optval}
\lambda^\dagger = \lambda^\dagger(R) = \tilde{\mathcal{H}}^{-1}_{+,f}(R^2).
\end{align}
Substituting (\ref{Optval}) into the minimisation problem (\ref{Minimize}) and then using (\ref{Uniq_minimum}), we arrive at
\begin{align}\label{Rep_formula}
\mathfrak{B}_{\varphi,+}(\mu\|\nu;f) = \Theta_{+}(\lambda^\dagger(R),R) = \nabla_\lambda\mathcal{G}_{\varphi,\nu}(\lambda^\dagger(R);f) = \nabla_\lambda\mathcal{G}_{\varphi,\nu}(\tilde{\mathcal{H}}^{-1}_{+,f}(R^2);f).
\end{align}
Finally, note that the above representation can be extended to include $R\,{=}\,0$ given that $\lambda^\dagger(0)\,{=}\,0$, $\tilde{\mathcal{H}}^{-1}_{f}(0)=0$, and $\nabla_\lambda\mathcal{G}_{\varphi,\nu}(0;f)\,{=}\,0$.
To complete the proof, we note that since $\mathcal{H}_{f}(\lambda)$ in (\ref{HHg}) is strictly increasing for $\lambda>0$, we have $\lim_{R\rightarrow \infty}\mathcal{H}_f(R^2) = \infty$  and the representation (\ref{Rep_formula}) continues to hold in terms of a pseudo-inverse $\tilde{\mathcal{H}}^{-1}_{+,f}$ of ${\mathcal{H}}_{f}$.

\qed

\subsection{Proof of Corollary \ref{UQ_lineariz}}\label{UQ_lineariz_ap} 
Similar to the proof of Proposition \ref{InformaRe} we define 
$$\Theta_{+}(\lambda,R) := \frac{1}{\lambda}\mathcal{G}_{\varphi,\nu}(\lambda;f)+\frac{1}{\lambda}R^2, \qquad  \Theta_{-}(\lambda,R) := -\frac{1}{\lambda}\mathcal{G}_{\varphi,\nu}(-\lambda;f)-\frac{1}{\lambda}R^2, $$  for $\lambda>0$,  where $R^2 = \Df(\mu\|\nu)$  so that 
\begin{equation}\label{Minimize22}
\mathfrak{B}_{\varphi,+}(\mu\|\nu;f) = \inf_{\lambda> 0}\Theta_{+}(\lambda,R), \qquad \mathfrak{B}_{\varphi,-}(\mu\|\nu;f) = \sup_{\lambda > 0}\Theta_{-}(\lambda,R).
\end{equation}
We describe  the proof for  $\mathfrak{B}_{\varphi,+}(\mu\|\nu;f)$, since the proof concerning  $\mathfrak{B}_{\varphi,-}(\mu\|\nu;f)$ is obtained by changing the sign of $\lambda$. First, we note that the existence of a unique minimiser $\lambda^\dagger$ of 
\begin{align}\label{PosP}
 \tag{$P^+$}  \inf_{\lambda>0}\bigg\{\frac{1}{\lambda}\mathcal{G}_{\varphi,\nu}(\lambda;f)+\frac{1}{\lambda}R^2\bigg\}
\end{align}
follows from the proof of Proposition \ref{InformaRe}; in particular,  for $0\leqslant R<\infty$ the unique minimiser $0\leqslant \lambda^\dagger(R)<\infty$ of $\Theta_+(\ccdot\,,R)$ satisfies
 \begin{align}\label{minim}
-\frac{1}{(\lambda^\dagger)^2}\mathcal{G}_{\varphi,\nu}(\lambda^\dagger;f)+\frac{1}{\lambda^\dagger}\nabla_\lambda\mathcal{G}_{\varphi,\nu}(\lambda^\dagger;f)-\frac{1}{(\lambda^\dagger)^2}R^2 =0,
\end{align}
with $\lambda^\dagger(0)=0$. Next, consider a function 
\begin{equation*}
H(\lambda,R):=-\frac{1}{\lambda^2}\Big(\mathcal{G}_{\varphi,\nu}(\lambda;f)-\lambda\nabla_\lambda\mathcal{G}_{\varphi,\nu}(\lambda;f)+R^2\Big),
\end{equation*}
and note that, given the properties of  $\mathcal{G}_{\varphi,\nu}(\lambda;f)$, we have
\begin{equation}\label{HHL}
H(\lambda,R):=-\frac{1}{2\lambda^2}\Big(\lambda^2\nabla^2_\lambda\mathcal{G}_{\varphi,\nu}(0;f)-2R^2+\mathcal{O}(\lambda^3)\Big),
\end{equation}
where $\nabla^2_\lambda\mathcal{G}_{\varphi,\nu}(0;g)= \nabla^2\varphi^*(0)\text{Var}_{\nu}(g)>0$ for $\E[g]\ne g$.
In particular, given the strictly increasing map $[0,\infty)\ni R\mapsto\lambda^\dagger(R)$ representing the unique solution of (\ref{minim}),
we have for $0\,{<}R\,{<}\,\infty$ that 
\begin{equation}\label{opt1}
(\lambda^\dagger)^2\nabla^2_\lambda\mathcal{G}_{\varphi,\nu}(0;f)-2R^2+\mathcal{O}((\lambda^\dagger)^3)=0,
\end{equation}
and, consequently,
\begin{equation}\label{lam_loc}
\lambda^\dagger(R) =C_{\varphi^*}R+\mathcal{O}(R^2),  \qquad C_{\varphi^*}:= \sqrt{2} \big(\nabla^2\varphi^*(0)\text{Var}_\nu(f)\big)^{-1/2}.
\end{equation}
It remains to note that  (\ref{lam_loc}) uniquely solves (\ref{minim}) for all $0\leqslant R<\infty$.  In a similar fashion, by a change of variable $\lambda\mapsto -\lambda,$ we obtain $-\lambda^\dagger(R)$ as the unique solution to the optimisation problem 
 \begin{align}\label{NegP}
 \tag{$P^{-}$} \sup_{\lambda>0}\bigg\{-\frac{1}{\lambda}\mathcal{G}_{\varphi,\nu}(-\lambda;f)-\frac{1}{\lambda}\Df(\mu\|\nu)\bigg\}.
 \end{align}
 Now, it only remains to expand the representation formula $\mathfrak{B}_{\varphi,\pm}(R;f) = \nabla_\lambda\mathcal{G}_{\varphi,\nu}(\tilde{\mathcal{H}}^{-1}_{f}(R^2);f)$ derived in part (2) of Proposition \ref{InformaRe} around $R\,{=}\,0$, with $R^2 = \Df(\mu\|\nu).$  Combining the expansion (\ref{lam_loc}) and $\lambda^\dagger(R) = \tilde{\mathcal{H}}_{f}^{-1}(R^2)$, $\tilde{\mathcal{H}}^{-1}_f(0) =0$ (cf.~(\ref{Optval})), with the fact that $\mathcal{G}_{\varphi,\nu}(0;f) =\nabla_\lambda\mathcal{G}_{\varphi,\nu}(0;f) =0$, we have 
 \begin{align*}
 \hspace{.5cm}\mathfrak{B}_{\varphi,+}(\mu\|\nu;f) = \nabla\mathcal{G}_{\varphi,\nu}(\tilde{\mathcal{H}}^{-1}_f(R^2)) &= \nabla_\lambda\mathcal{G}_\nu(0;f)+\nabla^2_\lambda\mathcal{G}_{\varphi,\nu}(0;f)\tilde{\mathcal{H}}^{-1}_f(R^2)+\mathcal{O}(|\tilde{\mathcal{H}}_f^{-1}(R^2)|^2)\\
 &= \nabla^2\varphi^*(0)\text{Var}_\nu(f)\lambda^\dagger(R) +\mathcal{O}(\lambda^\dagger(R)^2)\\
 &= \sqrt{2 \nabla^2\varphi^*(0)\text{Var}_\nu(f)}R+ \mathcal{O}(R^2)\\
 &= \sqrt{2 \nabla^2\varphi^*(0)\text{Var}_\nu(f)}\sqrt{\Df(\mu\|\nu)}+\mathcal{O}(\Df(\mu\|\nu)).  \hspace{1.2cm}\qed
 \end{align*}

\subsection{Proof of Proposition \ref{Lim_procedure}}\label{Uniell_prof2}
Given that this proposition first appears in the context when $\Xt=\XXt$, we follow the notation used in the corresponding section, where $\brr^\mm \equiv \brr^\mu$, $\sigma^\mm\equiv\sigma^\nu$. The proof of this proposition is same  when $\XXt\ne \Xt$, since the statement involves only the original dynamics in (\ref{SDE1}); the only difference would be replacing $\mu$ with $\mm$ and $\Xt$ with $\XXt$.   

The regularity assumptions on $(\brr^\mu, \sigma^\mu)$  in Theorem~\ref{Info_ineq2}  ensure that $x\mapsto \phi^\mu_{t_0,t}(x,\om)$ is a $\mathcal{C}^3$-diffeomorphism on $\Xt$, and that the derivative flow $h_t=D_x\phi^\mu_{t_0,t}(x,\om)h$ at $x\in \Xt$ in the direction of $h\in \Xt$, $|h|<\infty$, satisfying 
\begin{equation}
d h_t =D \brr^\mu(t,\phi^\mu_{t_0,t})h_tdt +D\sigma^\mu(t,\phi^\mu_{t_0,t})h_t\circ dW_t^\mu
\end{equation}
exists for almost all $\om \in \Om$ and $t\in\Ic := [t_0,\,t_0+T)$; e.g.,~\cite{Kunitabook} or Theorem \ref{SDE_flow}. In the above, $D\brr^\mu(t,\phi^\mu_{t_0,t}):= \nabla_\xi \brr^\mu(t,\xi)\big|_{\xi = \phi^\mu_{t_0,t}(x,\om)}$, and $D\sigma^\mu(t,\phi^\mu_{t_0,t})h_t:= \sum_k\big(\nabla_\xi \sigma_k^\mu(t,\xi)\big|_{\xi= \phi^\mu_{t_0,t}(x,\om)}\big)h_t$, where $\sigma^\mu_k$ are columns of $\sigma^\mu$. Note that, given the assumed growth and regularity conditions on  $(\brr^\mu, \; \sigma^\mu)$, the moments of $h_t$ are bounded on $\Ic$; this follows from the fact under such assumptions  the  derivative flow is a global $\mathcal{C}^2$-diffeomorphisms on $\Ic$ (see \cite{Kunitabook},  or Theorem~\ref{SDE_flow} and Remark~\ref{cf_remk}).

Consider the map $(r,x)\mapsto \mathcal{P}_{t_0,t-r}f(x)$ for $t_0\leqslant r\leqslant t$ with $f\in \mathcal{C}^2_\infty(\Xt)$. Then, application of It\^o's formula to  $\mathcal{P}_{t_0,t-r}f(\phi^\mu_{r,t_0}(x,\om))$, and taking the limit $r\rightarrow t$ leads to 
\begin{align}\label{SPoin}
f(\phi_{t_0,t}^\mu(x,\om)) = \mathcal{P}_{t_0,t}f(x)+\int_{t_0}^tD_x\big(\mathcal{P}_{s,t}f\big)(\phi_{t_0,s}^\mu(x,\om))\sigma^{\mu}\big(s,\phi^\mu_{t_0,s}(x,\om)\big)dW_s(\om).
\end{align}
Multiplying (\ref{SPoin}) by $ \int_{t_0}^t\big\langle \tilde{\sigma}^{\mu,-1}(s,\phi^\mu_{t_0,s}(x,\om)) D_x\phi^{\mu}_{t_0,s}(x,\om)h, dW_s(\om)\big\rangle$,  taking the expectation, and applying It\^o's isometry and Fubini's theorem leads to 
\begin{align*}
&\E\left[ f\big(\phi_{t_0,t}^\mu(x)\big)\int_{t_0}^t\big\langle \tilde{\sigma}^{\mu,-1}\big(s,\phi^\mu_{t_0,s}(x)\big) D_x\phi^{\mu}_{t_0,s}(x) h, dW_s\big\rangle\right]\\[.1cm]
&\hspace{6cm} = \E\left[ \int_{t_0}^t \big\langle D_x\big(\mathcal{P}_{s,t}f\big)\big(\phi_{t_0,s}^\mu(x)\big) D_x\phi^\mu_{t_0,s}(x), h\big\rangle ds\right]\\[.2cm]
&\hspace{6cm} = \int_{t_0}^t \big\langle D_x\E\big[\mathcal{P}_{s,t}f\big(\phi^\mu_{t_0,s}(x)\big)\big], h\big\rangle ds.
\end{align*}
By the flow property of the solutions $\{\phi_{t_0,t},\; t\geqslant t_0\},$ we have 
\begin{align*}
\E\big[\mathcal{P}_{s,t}f(\phi^\mu_{t_0,s}(x))\big] = \E\big[\E\big[f(\phi^\mu_{s,t}\circ\phi^\mu_{t_0,s}(x))\big]\big] =\E\big[ f(\phi^\mu_{t_0,t}(x))] = \mathcal{P}_{t_0,t}f(x),
\end{align*}
so that, for any $f\in \mathcal{C}_{\infty}^2(\Xt),$ we arrive at
\begin{align}\label{BsEll}
 (t-t_0) \big\langle D_x(\mathcal{P}_{t_0,t}f)(x), h\big\rangle=  \E\left[ f\big(\phi_{t_0,t}^\mu(x)\big)\int_{t_0}^t\big\langle\tilde{\sigma}^{\mu,-1}\big(s,\phi^\mu_{t_0,s}(x)\big) D_x\phi^{\mu}_{t_0,s}(x)h, dW_s\big\rangle\right].
\end{align}
Since $\mathcal{C}_{\infty}^2(\Xt)$ is dense in $\mathcal{C}_{\infty}^1(\Xt),$ we obtain a version of the {\it Bismut--Elworthy--Li formula}~(e.g.,~\cite{ EL93,Daprato08, Nulart, EL90}).
Moreover, since $\mathcal{C}_{\infty}^{1}(\Xt)$ is dense in $\mathcal{C}_{\infty}(\Xt),$ we have $(f_n)_{n\in \N}$, $f_n\in  \mathcal{C}_{\infty}^1(\Xt)$, such that $f_n\nearrow f\in \mathcal{C}_{\infty}(\Xt)$ and for 
\begin{align*}
\lim_{n\rightarrow \infty}\mathcal{P}^\mu_{t_0,t}f_n(x) &= \mathcal{P}^\mu_{t_0,t}f(x),\\
(t-t_0)\lim_{n\rightarrow \infty} \big\langle D_x(\mathcal{P}^\mu_{t_0,t}f_n)(x), h\big\rangle &= \E\Big[ f\big(\phi^\mu_{t_0,t}(x)\big)\int_{t_0}^t\big\langle\tilde{\sigma}^{\mu,-1}\big(s,\phi^\mu_{t_0,s}(x)\big)  D_x\phi^\mu_{t_0,s}(x) h, dW_s\big\rangle\Big].
\end{align*}
Given the regularity assumptions on the coefficients ${b}^\mu$, $\sigma^\mu$, of (\ref{SDE1}) there exists a fundamental solution of the Kolmogorov equations,  $0< p^\mu_{t_0,t}\in \mathcal{C}_{\infty}^\infty(\Xt)\times \mathcal{C}_{\infty}^\infty(\Xt)$, such that  (e.g., \cite{Roc-Kry})
\begin{alignat}{2}
  \mathcal{P}^\mu_{t_0,t}f(x) &= \int_\Xt p^\mu_{t_0,t}(x,y) f(y)dy,&\qquad &f\in \mathcal{C}_\infty(\Xt),\label{Ppdy}\\[.1cm]
 \mu_t(B)  = \int_B\rho_t^\mu(x)dx &= \int_B \int_\Xt p^\mu_{t_0,t}(x,y)\rho^\mu_{t_0}(y) dxdy, &\qquad &\rho^\mu_{t_0}\in L^1(\Xt,dx)\cap L^\infty(\Xt,dx),\label{Rpdy}
 \end{alignat} 
 which implies that 
 \begin{align*}
 \lim_{n\rightarrow \infty} \big\langle D_x(\mathcal{P}_{t_0,t}f_n)(x), h\big\rangle &= \lim_{n\rightarrow \infty}\int_{\Xt} f_n(y) \big\langle D_xp^\mu_{t_0,t}(x,y), h\big\rangle dy\\
 &=\int_{\Xt} f(y) \big\langle D_x p^\mu_{t_0,t}(x,y), h\big\rangle dy = \big\langle D_x(\mathcal{P}_{t_0,t}f)(x),h\big\rangle,
 \end{align*}
 so that (\ref{BsEll}) holds for all $f\in \mathcal{C}_{\infty}(\Xt),$ i.e., 
 \begin{align}\label{BsEll2}
(t-t_0) \big\langle D_x(\mathcal{P}_{t_0,t}f)(x), h\big\rangle=  \E\left[ f\big(\phi_{t_0,t}^\mu(x)\big)\int_{t_0}^t\big\langle\tilde{\sigma}^{\mu,-1}\big(s,\phi^\mu_{t_0,s}(x)\big) D_x\phi^{\mu}_{t_0,s}(x)h, dW_s\big\rangle\right].
 \end{align}
Next,  utilising (\ref{Ppdy}), (\ref{Rpdy}) 
and the law of total expectation in (\ref{BsEll2}) we have 
\begin{align}\label{BsEll3}
\notag &(t-t_0)\Big\langle D_x \int_{\Xt} f(y)p^{\mu}_{t_0,t}(x,y)dy, h\Big\rangle \\
&\hspace{.5cm} = \int_{\Xt}\int_{\Xt} f(y)p^\mu_{t_0,t}(\xi,y)\rho^\mu_{t_0}(\xi)\E\bigg[\int_{t_0}^t \langle\tilde{\sigma}^{\mu,-1}(s,\phi^{\mu}_{t_0,s}(x))  D_x\phi^\mu_{t_0,s}(x)h, dW_s\rangle\Big| \phi^\mu_{t_0,t}(x)=y\bigg] d\xi dy\notag\\
&\hspace{.5cm} \leqslant  \mathfrak{C}\int_{\Xt} f(y)p^\mu_{t_0,t}(x,y)\E\bigg[\int_{t_0}^t \langle\tilde{\sigma}^{\mu,-1}(s,\phi^{\mu}_{t_0,s}(x)) D_x\phi^\mu_{t_0,s}(x)h, dW_s\rangle\Big| \phi^\mu_{t_0,t}(x)=y\bigg] dy.
\end{align}

\noindent In particular, (\ref{BsEll3}) holds for any $f\in \mathcal{C}_c^{+}(\Xt),$ so that for $t>t_0$ we have 
\begin{align*}
&\langle D_x p^\mu_{t_0,t}(x,y),h\rangle\\
&\hspace{2.1cm} \leqslant p^\mu_{t_0,t}(x,y)\frac{\mathfrak{C}}{t-t_0}\E\bigg[\int_{t_0}^t \langle\tilde\sigma^{\mu,-1}(s,\phi^{\mu}_{t_0,s}(x))  D_x\phi^\mu_{t_0,s}(x) h, dW_s\rangle\Big| \phi^\mu_{t_0,t}(x)=y\bigg],\hspace{.3cm}
\end{align*}
 which can be written as 
\begin{align*}
\langle D_x\log p^\mu_{t_0,t}(x,y), h\rangle \leqslant  \frac{\mathfrak{C}}{t-t_0}\E\bigg[\int_{t_0}^t \langle\tilde\sigma^{\mu,-1}(s,\phi^{\mu}_{t_0,s}(x)) D_x\phi^\mu_{t_0,s}(x) h, dW_s\rangle\Big| \phi^\mu_{t_0,t}(x)=y\bigg],
\end{align*}
with the convention that $D_x \log p^\mu_{t_0,t}(x,y) =0$ if $p^\mu_{t_0,t}(x,y) =0$. Next, apply Jensen's inequality for the conditional expectation and It\^o isometry to obtain
\begin{align*}
\vert \langle D_x\log p^\mu_{t_0,t}(x,y), h\rangle \vert^2 &\leqslant \frac{\mathfrak{C}}{(t-t_0)^2}\bigg(\E\left[\int_{t_0}^t\langle \tilde\sigma^{\mu,-1}(s,\phi^\mu_{t_0,s}(x)) D_x\phi^\mu_{t_0,s}(x)h, dW_s\rangle\Big | \phi^\mu_{t_0,t}(x)=y \right] \bigg)^2\notag\\
&\leqslant \frac{\mathfrak{C}}{(t-t_0)^2}\E\bigg[\Big(\int_{t_0}^t \langle \tilde\sigma^{\mu,-1}(s,\phi^{\mu}_{t_0,s}(x)) D_x\phi^\mu_{t_0,s}(x) h, dW_s\rangle\Big)^2\Big| \phi^\mu_{t_0,t}(x)=y\Big]\notag\\
&\leqslant \frac{\mathfrak{C}}{(t-t_0)^2}\E\bigg[\int_{t_0}^t\big\vert\tilde\sigma^{\mu,-1}(s,\phi^{\mu}_{t_0,s}(x))D_x\phi^\mu_{t_0,s}(x) h\big\vert^2 ds\,\Big| \phi^\mu_{t_0,t}(x)=y\Big].
\end{align*}
Next, note that 
\begin{align*}
\langle D_x\rho_t^\mu(x), h\rangle &= \int_{\Xt}\langle D_x p^\mu_{t_0,t}(x,y), h\rangle \rho^\mu_{t_0}(y)dy
=\int_{\Xt}\langle D_x\log p^\mu_{t_0,t}(x,y), h\rangle p^\mu_{t_0,t}(x,y) \rho^\mu_{t_0}(y)dy\notag\\
&\leqslant \sup_{y\in \Xt}\big\vert \langle D_x\log p^\mu_{t_0,t}(x,y), h\rangle\big\vert \int_{\Xt} p^\mu_{t_0,t}(x,y)\rho^\mu_{t_0}(y)dy\notag\\
&= \sup_{y\in \Xt}\big\vert \langle D_x\log p^\mu_{t_0,t}(x,y), h\rangle\big\vert \rho^\mu_t(x).
\end{align*}
Finally, combining the two bounds above leads to  
\begin{align*}
\big\vert \langle D_x\log \rho^\mu_t(x), h\rangle\big\vert^2\leqslant \frac{\mathfrak{C}}{(t-t_0)^2}\sup_{y\in \Xt}\E\bigg[\int_{t_0}^t\big\vert \tilde\sigma^{\mu,-1}(s,\phi^\mu_{t_0,s}(x))\cdot D_x\phi^\mu_{t_0,s}(x)h\big\vert^2 ds\Big| \phi^\mu_{t_0,t}(x)=y\bigg]. 
\end{align*}
In particular, take $|h|=1$  and use the fact that $0\,{<}\mathfrak{C}\,\|\tilde\sigma^{\mu,-1}\|_\textsc{hs}\leqslant  \mathfrak{C}_{\sigma^\mu}\,{<}\,\infty$ to obtain 
\begin{align*}
\notag \hspace{1.cm}&\big\vert \nabla_x\log\rho^{\mu}_{t}(x)\big\vert^2\leqslant \frac{ \mathfrak{C}_{\sigma^\mu}}{(t-t_0)^2}\sup_{y\in \Xt}\E\left[ \int_{t_0}^t \big\Vert D_x\phi^{\mu}_{t_0,s}(x)\big\Vert^2_{\textsc{hs}}ds\Big| \Solm^{\mu}(x) = y\right]<\infty, \quad t>t_0\hspace{.5cm}  \qed
\end{align*}
which follows from the existence  of the  second moment of the derivative flow for  $t\in [t_0,\,t_0+T]$. The bound on the first two derivatives of $\nabla_x\log\rho_t^\mu$ is asserted by following similar steps, utilizing the fact that the moments of $D_{xx} \phi_{t_0,t}^\mu(x)$, and $D_{xxx} \phi_{t_0,t}^\mu(x)$ exist (given the conditions on the coefficients; see \cite[Corollary 4.6.7]{Kunitabook}) and the regularity of $p^\mu_{t_0,t}$).

\subsection{Proof of Lemma \ref{Ide_Lem}}\label{Ide_Lem_ap}
Given the assumed regularity of the coefficients $(b,\sigma)$, we proceed by chain and product rules of differentiation, to obtain (recall that $a = \sigma\sigma^*$, $a_{ij} = \sigma_{ik}\sigma_{jk}$)
\begin{align*}
\LG_t^*\varphi(f) &= -\partial_{x_i}(b_i\varphi(f))+\frat\partial^2_{x_ix_j}(\sigma_{ik}\sigma_{jk}\varphi(f))\\
&= -\varphi(f)\partial_{x_i}b_i-\varphi^{\prime}(f)b_i\partial_{x_i}f+\frat\Big( \varphi^{\prime}(f)\sigma_{ik}\sigma_{jk}\partial^{2}_{x_ix_j}f+ \varphi^{\prime\prime}(f)\sigma_{ik}\sigma_{jk}\partial_{x_i}f\partial_{x_j}f\\
&\hspace{2cm} + \varphi^{\prime}(f)\partial_{x_j}(\sigma_{ik}\sigma_{jk})\partial_{x_i}f + \varphi^{\prime}(f)\partial_{x_j}(\sigma_{ik}\sigma_{jk})\partial_{x_j}f+\varphi(f)\partial^2_{x_ix_j}(\sigma_{ik}\sigma_{jk})\Big)\\
&= -\varphi(f)\partial_{x_i}b_i+\varphi^{\prime}(f)f\partial_{x_i}b_i-\varphi^{\prime}(f)\partial_{x_i}(b_if)+ \frat\Big(\varphi^{\prime}(f)\partial^2_{x_ix_j}(\sigma_{ik}\sigma_{jk}f)\\
&\hspace{2cm}-\varphi^{\prime}(f)f\partial^2_{x_ix_j}(\sigma_{ik}\sigma_{jk})+\varphi(f)\partial^2_{x_ix_j}(\sigma_{ik}\sigma_{jk}) +  \varphi^{\prime\prime}(f)\sigma_{ik}\sigma_{jk}\partial_{x_i}f\partial_{x_j}f\Big)\\
&=\varphi^{\prime}(f)\big[ -\partial_{x_i}(b_if) +\frat\partial_{x_ix_j}(\sigma_{ik}\sigma_{jk}f)\big]+ \left(\varphi^{\prime}(f)f-\varphi(f)\right)\partial_{x_i}\left( b_i -\partial_{x_j}(\sigma_{ik}\sigma_{jk})\right)\\
& \hspace{2cm} + \frat \varphi^{\prime\prime}(f)\sigma_{ik}\sigma_{jk}\partial_{x_i}f\partial_{x_j}f.
\end{align*}
Next, we verify the second identity.  In a fashion similar to the procedure used above, we obtain
\begin{align*}
\hspace{.8cm}\LG_t^*(fg)& = -\partial_{x_i}(b_ifg)+\frat\partial^2_{x_ix_j}(\sigma_{ik}\sigma_{jk}fg)\\
&= -fg\partial_{x_i}b_i -f[\partial_{x_i}(b_ig) -g\partial_{x_i}b_i] -g[\partial_{x_i}(b_if) -f\partial_{x_i}b_i]\\
&\hspace{.5cm}+\frat\Big( -fg\partial^2_{x_ix_j}(\sigma_{ik}\sigma_{jk})-\partial_{x_i}(fg)\partial_{x_j}(\sigma_{ik}\sigma_{jk})+f\partial^2_{x_ix_j}(\sigma_{ik}\sigma_{jk}g)\\
&\hspace{1.5cm} +g\partial^2_{x_ix_j}(\sigma_{ik}\sigma_{jk}f)+2\sigma_{ik}\sigma_{jk}\partial_{x_i}f\partial_{x_j}g+[g\partial_{x_i}f+f\partial_{x_i}g]\partial_{x_j}(\sigma_{ik}\sigma_{jk})\Big)\\
&=  -f\partial_{x_i}(b_ig)+\frat f\partial^2_{x_ix_j}(\sigma_{ik}\sigma_{jk}g)  -g\partial_{x_i}(b_if)+\frat g\partial^2_{x_ix_j}(\sigma_{ik}\sigma_{jk}f)\\
&\hspace{.5cm} + \sigma_{ik}\sigma_{jk}\partial_{x_i}f\partial_{x_j}g+fg\partial_{x_i}[b_i-\frat\partial_{x_j}(\sigma_{ik}\sigma_{jk})]. \hspace{5.5cm}\qed
\end{align*}

\subsection{Proof of Lemma \ref{Pre_lem2}}\label{Pre_lem2_ap}
Proceeding as in \cite[Lemma 2.4]{Rockner16}, we recall from Lemma \ref{Ide_Lem} that for all test functions $f,g\in \mathcal{C}^{2}(\Ic\times \Xt)$ and $\varphi\in \mathcal{C}^{2}(\Rp),$ we have (for $h^\nu:=b^\nu-\frat\nabla a^\nu$)
\begin{align}\label{GW1}
\begin{cases}
\mathcal{L}_t^{\nu*}\varphi(f) &= \varphi^{\prime}(f)\LG_t^{\nu^*}f+ \frac{1}{2}\varphi^{\prime\prime}(f)\langle a^{\nu}\nabla f, \nabla f\rangle + \left(f\varphi^{\prime}(f) - \varphi(f)\right)\nabla\cdot h^\nu,\\
\LG_t^{\nu^*}(fg)& = f\LG_t^{\nu^*}g+g\LG_t^{\nu^*}f + \langle a^{\nu} \nabla f, \nabla g\rangle + fg \,\nabla\cdot h^\nu.
\end{cases}
\end{align}
Next, from the assumptions on the coefficients of the SDEs and the initial data, we have that the densities $\rho_t^\mu$ and $\rho_t^\nu$ are strictly positive,  so that $\eta_t = \rho_t^\mu/\rho_t^\nu$ is strictly positive, differentiable and finite.  This follows from the fact that for probability measures on the initial conditions with strictly positive (Lebesgue) densities $\rho^\mu_{t_0},\rho^\nu_{t_0}$, and the dynamics generated by flows of diffeomorphisms, solutions of the forward Kolmogorov equations~(\ref{F_Kol}) are absolutely continuous w.r.t.~the Lebesgue measure and have  strictly positive densities. Moreover,  densities associated with uniformly elliptic dynamics are differentiable for  $t>t_0$ regardless of the initial density (see, e.g.,~\cite{Stroock79, Roc-Kry, Figali, Gmeas}).
Now, consider the solutions of the forward Kolmogorov equations  
\begin{equation}
(a)\;\;\;\partial_t\rho_t^\nu = \LG_t^{\nu^*}\rho_t^\nu,\qquad\quad  (b)\;\;\;\partial_t\rho_t^\mu = \LG_t^{\nu^*}\rho_t^\mu -\nabla\cdot(\varTheta_{\mu\nu}\rho_t^\mu),
\end{equation}
and the equation 
\begin{align}\label{GW2}
\partial_t\rho_t^\mu -\eta_t\partial_t\rho_t^\nu = \LG_t^{\nu^*}\!\rho_t^\mu-\eta_t\LG_t^{\nu^*}\!\rho_t^\nu -\nabla\cdot(\betaLk \rho_t^\mu).
\end{align}
Observe, that $\eta_t\rho_t^\nu = \rho_t^\mu$ implies that 
$\partial_t\rho_t^\mu-\eta_t\partial_t\rho_t^\nu = \rho_t^\nu\partial_t\eta_t$ and the identities (\ref{GW1}) lead to
\begin{align*}
\LG_t^{\nu^*}\rho_t^\mu -\eta_t\LG_t^{\nu^*}\!\rho_t^\nu = \LG_t^{\nu^*}\!(\eta_t\rho_t^\nu) -\eta_t\LG_t^{\nu^*}\!\rho_t^\nu = \rho_t^\nu\LG_t^{\nu^*}\!\eta_t+\langle a^\nu\nabla\rho_t^\nu, \nabla\eta_t\rangle+\rho_t^\nu\eta_t\nabla\cdot h^\nu.
\end{align*}
The equation (\ref{GW2}) becomes 
\begin{align}\label{GW3}
\rho_t^\nu\partial_t\eta_t = \rho_t^\nu\LG_t^{\nu^*}\!\eta_t + \langle a^\nu\nabla\rho_t^\nu, \nabla\eta_t\rangle+\rho_t^\nu\eta_t\nabla\cdot h^\nu -\nabla\cdot (\betaLk \rho_t^\mu).
\end{align}
Multiplying both sides of (\ref{GW3}) by $\varphi^{\prime}(\eta_t)$ and noticing that 
\begin{align*}
\partial_t\varphi(\eta_t) = \varphi^{\prime}(\eta_t)\partial_t\eta_t, \quad \text{and}\quad \nabla\varphi(\eta_t) = \varphi^{\prime}(\eta_t)\nabla\eta_t,
\end{align*}
we have 
\begin{align}\label{GW4}
\rho_t^\nu\partial_t\varphi(\eta_t) = \varphi^{\prime}(\eta_t)\rho_t^\nu\LG_t^{\nu^*}\!\eta_t+\langle a^\nu\nabla\rho_t^\nu, \nabla\varphi(\eta_t)\rangle+ \eta_t\varphi^{\prime}(\eta_t)\rho_t^\nu\nabla\cdot h^\nu -\varphi^{\prime}(\eta_t)\nabla\cdot(\betaLk\rho_t^\mu).
\end{align}
From the identities (\ref{GW1}), we write $\varphi^{\prime}(\eta_t)\LG_t^{\nu^*}\eta_t $ as follows
\begin{align}\label{Gww}
\varphi^{\prime}(\eta_t)\LG_t^{\nu^*}\eta_t = \LG_t^{\nu^*}\!\varphi(\eta_t)-\frat\varphi^{\prime\prime}(\eta_t)\langle a^\nu\nabla\eta_t, \nabla\eta_t\rangle- (\eta_t\varphi^{\prime}(\eta_t) -\varphi(\eta_t))\nabla\cdot h^\nu, 
\end{align}
and substituting (\ref{Gww}) into (\ref{GW4}), 
we have 
\begin{align}\label{GW5}
\notag & \rho_t^\nu\partial_t\varphi(\eta_t) = \rho_t^\nu\LG_t^{\nu^*}\!\varphi(\eta_t)-\frat\rho_t^\nu\varphi^{\prime\prime}(\eta_t)\langle a^\nu\nabla\eta_t, \nabla\eta_t\rangle +\langle a^\nu\nabla\rho_t^\nu, \nabla\varphi(\eta_t)\rangle \\ &\hspace{2cm} +\varphi(\eta_t)\rho_t^\nu\nabla\cdot h^\nu  -\varphi^{\prime}(\eta_t)\nabla\cdot(\betaLk\rho_t^\mu).
\end{align}
Addition of $\varphi(\eta_t)\partial_t\rho_t^\nu = \varphi(\eta_t)\LG_t^{\nu^*}\!\rho_t^\nu$ to the (\ref{GW5}) leads to
\begin{align*}
\notag \rho_t^\nu\partial_t\varphi(\eta_t)+\varphi(\eta_t)\partial_t\rho_t^\nu &= \rho_t^\nu\LG_t^{\nu^*}\!\varphi(\eta_t)+\varphi(\eta_t)\LG_t^{\nu^*}\rho_t^\nu + \langle a^\nu\nabla\rho_t^\nu, \nabla\varphi(\eta_t)\rangle\\ &\hspace{.5cm} + \varphi(\eta_t)\rho_t^\nu\nabla\cdot h^\nu -\frat\rho_t^\nu\varphi^{\prime\prime}(\eta_t)\langle a^\nu\nabla\eta_t, \nabla\eta_t\rangle -\varphi^{\prime}(\eta_t)\nabla\cdot(\betaLk\rho_t^\mu).
\end{align*}
By the product rule together with the identities (\ref{GW1}), we arrive at 
\begin{align*}
\hspace{1cm}\partial_t(\varphi(\eta_t)\rho_t^\nu) = \LG_t^{\nu^*}(\varphi(\eta_t)\rho_t^\nu)-\frat\rho_t^\nu\varphi^{\prime\prime}(\eta_t)\langle a^\nu\nabla\eta_t, \nabla\eta_t\rangle -\varphi^{\prime}(\eta_t)\nabla\cdot(\betaLk\rho_t^\mu). \hspace{2cm}\qed
\end{align*}

\subsection{Proof of Lemma \ref{Integrand}}\label{Integrand_ap}

Recall from Lemma \ref{Pre_lem2}  that 
\begin{align}\label{Pre_eq1}
\partial_t(\varphi(\eta_t)\rho_t^\nu) = \LG_t^{\nu^*}(\varphi(\eta_t)\rho_t^\nu)-\frat\rho_t^\nu\varphi^{\prime\prime}(\eta_t)\langle a^\nu\nabla\eta_t, \nabla\eta_t\rangle -\varphi^{\prime}(\eta_t)\nabla\cdot(\betaLk\rho_t^\mu).
\end{align} 
 Multiply both sides of (\ref{Pre_eq1}) by a test function $f\in \mathcal{C}_c^{\infty}(\Xt)$ and integrate to arrive at 
\begin{align}\label{GW20}
\notag &\int_{\tau}^{t}\int_{\Xt}\partial_s(\varphi(\eta_s)\rho_s^\nu)f(x)dxds+\frac{1}{2}\int_{\tau}^t\int_{\Xt}\varphi^{\prime\prime}(\eta_s)\langle a^\nu\nabla\eta_s, \nabla\eta_s\rangle f(x)\nu_s(dx)ds\\
&\hspace{3cm} = \int_{\tau}^t\int_{\Xt}\varphi(\eta_s)\LG_t^\nu f(x)\nu_s(dx)ds -\int_{\tau}^t\varphi^{\prime}(\eta_s)\nabla\cdot(\betaLk\rho_s^\mu)f(x)dxds.
\end{align}
By the Newton-Leibniz formula, we have 
\begin{align}\label{GW21}
\int_{\tau}^{t}\int_{\Xt}\partial_s(\varphi(\eta_s)\rho_s^\nu)f(x)dx = \int_{\tau}^{t}\int_{\Xt}\varphi(\eta_t)f(x)\nu_t(dx)-\int_{\tau}^{t}\int_{\Xt}\varphi(\eta_\tau)f(x)\nu_\tau(dx).
\end{align}
Also, one has 
\begin{align}\label{GW22}
-\int_{\tau}^t\!\!\int_{\Xt}\varphi^{\prime}(\eta_s)\nabla\cdot(\betaLk\rho_s^\mu)f(x)dxds= \int_{\tau}^\tau\!\!\int_{\Xt}\!\left[ \varphi^{\prime\prime}(\eta_s)\langle \betaLk, \nabla\eta_s\rangle f(x) +\langle \betaLk, \nabla f(x)\rangle\right]\eta_s\nu_s(dx)ds.
\end{align}
Substituting (\ref{GW21}) and (\ref{GW22}) into (\ref{GW20}), we obtain the required equality (\ref{GW19}).
\qed

\subsection{Proof of Theorem \ref{difference_bound}}\label{difference_bound_app}

For a finite time interval $\Ic$ it is sufficient to consider $\Df(\mu_t\|\nu_t)$; the result can be then extended to divergence rates on unbounded intervals provided that the respective ratios are finite. Since, $\mu_{t_0}=\nu_{t_0}$ is assumed from the outset, we use denote both initial measures by $\mu_{t_0}$ throughout. We want to show that 
\begin{align}\label{statm}
\Df(\mu_t\|\nu_t)\leqslant \begin{cases} \hspace{0.7cm}\Df(\mu_t\|\mu_{t_0})-\Df(\nu_t\|\mu_{t_0}), \; \text{if}\; \Df(\mu_t\|\mu_{t_0})\geqslant \Df(\nu_t\|\mu_{t_0}),\\ -\left(\Df(\mu_t\|\mu_{t_0})-\Df(\nu_t\|\mu_{t_0})\right), \; \text{if}\; \Df(\nu_t\|\mu_{t_0})\geqslant \Df(\mu_t\|\mu_{t_0}).\end{cases}
\end{align}
Consider the convex set $ L^{1,\varphi}(\Xt; \mu_{t_0}) =\{f: \Xt\rightarrow [a, b]: f, \varphi(f)\in L^{1}(\Xt;\mu_{t_0})\}$ for some $0<a\leqslant b<\infty$, and set $d\mu_t = \eta_t^\mu d\mu_{t_0}$.  By the assumptions of the theorem, we have  $\varphi(\eta_t^\mu) = \varphi\left(d\mu_t/d\mu_{t_0}\right)\in L^1(\Xt;\mu_{t_0})$. Let $h\in L^{1,\varphi}(\Xt;\mu_{t_0})$ and  define $\beta: [0, 1]\rightarrow\Rp$ by
\begin{align}\label{ConvEn}
\beta_t(s) = \big\langle \varphi(s\eta_t^\mu+(1-s)h), \mu_{t_0}\big\rangle,
\end{align} 
so that $\beta_t(1) = \Df(\mu_t\|\mu_{t_0})$. One can check that $s\mapsto\beta_t(s)$ is convex, and as $\varphi\in \mathcal{C}^2(\Rp)$ that  the normality conditions (\ref{Normality}) ensure that $s\mapsto \beta_t(s)$ is twice continuously differentiable with bounded derivatives. In particular, we have 
\begin{align*}
\beta_t^{\prime}(s) = \big\langle \nabla\varphi\left(s\eta^\mu_t+(1-s)h\right)(\eta^\mu_t-h), \mu_{t_0}\big\rangle.
\end{align*}
Recalling that a smooth convex function is the envelope of its tangents (e.g.,~\cite{Rockafeller74}), we arrive at 
\begin{align*}
\Df(\mu_t\|\mu_{t_0})=\beta_t(1) = \sup_{s\in [0,1]}\big\{\beta_t(s)+ \beta_t^{\prime}(s)(1-s)\big\} \geqslant &\beta_t(0) +\beta_t^{\prime}(0)\\
 &= \big\langle \varphi(h), \mu_{t_0}\big\rangle + \big\langle \nabla\varphi(h)(\eta_t^\mu-h), \mu_{t_0}\big\rangle.
\end{align*}
This implies that 
\begin{align}\label{Duality2}
\Df(\mu_t\|\mu_{t_0})\geqslant\big\langle \nabla\varphi(h)(\eta_t^\mu-h), \mu_{t_0}\big\rangle + \big\langle \varphi(h), \mu_{t_0}\big\rangle \quad \forall \,h\in L^{1,\varphi}(\Xt; \mu_{t_0}).
\end{align}
In particular, take $h_1\in L^{1,\varphi}(\Xt;\mu_{t_0})$ s.t. $\varphi(h_1) = 2 \varphi(\eta_t^\nu)$ with $\eta_t^\nu= d\nu_t/d\mu_{t_0}$ and observe that $$\big\langle \varphi(h_1), \mu_{t_0}\big\rangle = 2\int_{\Xt}\varphi(\eta_t^\nu)d\mu_{t_0} = 2\int_{\Xt}\varphi\left(\frac{d\nu_t}{d\mu_{t_0}}\right)d\mu_{t_0} = 2\hspace{0.02cm}\Df(\nu_t\|\mu_{t_0}).$$
Then, from (\ref{Duality2}), we obtain
\begin{align}\label{Duality3}
 \frac{1}{2}\Df(\mu_t\|\mu_{t_0})\geqslant  \frac{1}{2}\big\langle\nabla\varphi(h_1)(\eta_t^\mu-h_1), \mu_{t_0}\big\rangle +  \Df(\nu_t\|\mu_{t_0}).
\end{align}
On the other hand, take $h_2\,{\in}\, L^{1,\varphi}(\Xt;\mu_{t_0})$ s.t. \!$\varphi(h_2)\,{=}\, 2 \varphi(\eta_t)\eta^\nu_t $ with $\eta_t \,{=}\, d\mu_t/d\nu_t, \; \eta_t^\nu \,{=}\, d\nu_t/d\mu_{t_0}$ and notice that 
\begin{align*}
\big\langle\varphi(h_2), \mu_{t_0}\big\rangle= 2\int_{\Xt}\varphi\Big(\frac{d\mu_t}{d\nu_t}\Big)\frac{d\nu_t}{d\mu_{t_0}}d\mu_{t_0} = 2\int_{\Xt}\varphi\Big(\frac{d\mu_t}{d\nu_t}\Big)d\nu_t = 2\,\Df(\mu_t\|\nu_t).
\end{align*}
From the inequality (\ref{Duality2}), we have 
\begin{align}\label{Duality4}
 \frac{1}{2}\Df(\mu_t\|\mu_{t_0})\geqslant   \frac{1}{2}\langle\nabla\varphi(h_2)(\eta_t^\mu-h_2), \mu_{t_0}\rangle+ \Df(\mu_t\|\nu_t).
\end{align}
Combining (\ref{Duality3}) and (\ref{Duality4}) yields \begin{align}\label{Duality6}
\Df(\mu_t\|\nu_t)\leqslant \Df(\mu_t\|\mu_{t_0}) - \Df(\nu_t\|\mu_{t_0})+ C(\mu_{t_0}, \mu_t, \nu_t, \nabla\varphi),
\end{align}
where 
\begin{align*}
2C(\mu_{t_0}, \mu_t, \nu_t, \nabla\varphi) := -\langle\nabla\varphi(h_2)(\eta_t^\mu-h_2), \mu_{t_0}\rangle- \langle\nabla\varphi(h_1)(\eta_t^\mu-h_1), \mu_{t_0}\rangle.
\end{align*}
Next, we derive a bound on $C(\mu_{t_0}, \mu_t, \nu_t, \nabla\varphi)$  as follows:
\begin{align}
\notag & 2C(\mu_{t_0}, \mu_t, \nu_t,\nabla\varphi) -2\Df(\nu_t\|\mu_{t_0})+2\Df(\nu_t\|\mu_{t_0})+2\Df(\mu_t\|\nu_t)\\
\nonumber &\hspace{1cm}\leqslant 2C(\mu_{t_0}, \mu_t, \nu_t, \nabla\varphi )+ \langle\varphi(h_1), \mu_{t_0}\rangle+ \langle \varphi(h_2), 
\mu_{t_0}\rangle -2\Df(\nu_t\|\mu_{t_0})\\ 
\notag &\hspace{1.cm}= -\langle\nabla\varphi(h_2)(\eta_t^\mu-h_2), \mu_{t_0}\rangle + \langle \varphi(h_2), \mu_{t_0}\rangle\\
&\hspace{2.2cm}- \langle\nabla\varphi(h_1)(\eta_t^\mu-h_1), \mu_{t_0}\rangle + \langle \varphi(h_1), \mu_{t_0}\rangle-2\Df(\nu_t\|\mu_{t_0}) \notag\\
 &\hspace{1cm}\leqslant 2\sup_{h\in L^{1,\varphi}(\mu_{t_0})}\Big\{ \langle \nabla\varphi(h)(\eta_t^\mu-h), \mu_{t_0}\rangle + \langle\varphi(h), \mu_{t_0}\rangle\Big\} -2\Df(\nu_t\|\mu_{t_0})  \notag\\
 &\hspace{1cm}= 2\Df(\mu_t\|\mu_{t_0}) -2\Df(\nu_t\|\mu_{t_0}), \label{Cderiv}
\end{align}
where the last equality is based on Proposition 1.2 in \cite{Chafi04}. The above bound implies that 
\begin{align}\label{Duality7}
C(\mu_{t_0}, \mu_t, \nu_t, \nabla\varphi)\leqslant  - \Df(\mu_t\|\nu_t)+\Df(\mu_t\|\mu_{t_0})-\Df(\nu_t\|\mu_{t_0})
\end{align}
and, consequently, the inequality (\ref{Duality6}) implies that for $\Df(\mu_t\|\mu_{t_0})\geqslant \Df(\nu_t\|\mu_{t_0})$ we have 
\begin{align}\label{Spil1}
\Df(\mu_t\|\nu_t)\leqslant \Df(\mu_t\|\mu_{t_0})-\Df(\nu_t\|\mu_{t_0}).
\end{align}

To complete the proof, and prove the second part of (\ref{statm}) we exchange the role of $\eta^\mu_t$ with  that of $d\nu_t=\eta_t^\nu d\mu_{t_0}$ in  (\ref{ConvEn}).  Following steps analogous to those above leads to 
\begin{align*}
\hspace{4cm}\Df(\mu_t\|\nu_t)\leqslant -\Big(\Df(\mu_t\|\mu_{t_0})-\Df(\nu_t\|\mu_{t_0})\Big)\hspace{4.2cm} 
\end{align*} 
when $\Df(\nu_t\|\mu_{t_0})\geqslant \Df(\mu_t\|\mu_{t_0})$.
\qed

\subsection{Proof of Lemma \ref{Absestm}}\label{Absestm_Ap}
Parts (a) and part (b) were proved long ago in \cite{Pinsky04}. Since $\Ms, \Ns$ are probability measures on the Polish space $\mathcal{W}_d$, there exists a unique Lebesgue decomposition of $\Ns$ given by $\Ns = \Ns^{ac}+\Ns^s$, such that $\Ns^{ac}\ll \Ms$ and $\Ns^s\perp \Ms.$  Given that  $\Ns|_{\F_n}\ll \Ms|_{\F_n}$ for all $n\in \N,$  we have (see, e.g.,~\cite{Pinsky04, Yor}),
\begin{align}\label{Abs1}
\Ns^{ac}(A) &= \int_{A}\lim_{n\rightarrow\infty}\frac{d\Ns|_{\F_n}}{d\Ms|_{\F_n}}d\Ms, \qquad A\in \F_t;\\
\Ns^{s}(A) &=\left(\left\{\om: \limsup_{n\rightarrow\infty}\frac{d\Ns|_{\F_n}}{d\Ms|_{\F_n}}=\infty\right\}\cap A\right), \qquad A\in \F_t.\label{Abs11}
\end{align}
The equalities (\ref{Abs1}) and (\ref{Abs11}) yield part (a) and part (b). For part (c), set 
\begin{align*}
D_n = \E^\Ms\left[\frac{d\Ns}{d\Ms}\Big|_{\F_n}\right] = \frac{d\Ns|_{F_n}}{d\Ms|_{\F_n}}, \qquad \Ms\; \text{-a.s.};
\end{align*}
it is straightforward from part (a) that $(D_n)_{n\in\N}$ is a $\Ms$-uniform integrable martingale. As $\F_t$ is compactly generated by $(\F_n)_{n\in\N},$ we obtain that 
\begin{align*}
\lim_{n\rightarrow\infty}D_n = D:=\frac{d\Ns}{d\Ms}, \qquad \Ms\; \text{-a.s.}
\end{align*} 
Next, given that  the strictly convex function $\varphi$ satisfies the normality condition (\ref{Normality}) and the assumptions in part (c), by Fatou's lemma, we have 
\begin{align}\label{JEN1}
\liminf_{n\rightarrow\infty}\Df\big(\Ns|_{\F_n}\|\Ms|_{\F_n}\big)=\liminf_{n\rightarrow\infty}\E^\Ms\big[\varphi\left(D_n\right)\big]\geqslant \E^\Ms\big[\varphi(D)\big] = \Df\big(\Ns\|\Ms\big).
\end{align} 
Conversely,  Jensen's inequality for conditional expectation yields
\begin{align*}
\E^\Ms\big[\varphi(D)|_{\F_n}\big]\geqslant \varphi\left(D_n\right), \qquad \Ms\;\text{-\,a.s.,}
\end{align*}
which implies  that $\Df\big(\Ns\|\Ms\big)\geqslant \Df\big(\Ns|_{\F_n}\|\Ms|_{\F_n}\big)$.
\qed

\subsection{Proof of Lemma \ref{Martinest}}\label{Martinest_Ap}
We proceed as in \cite{Pinsky04}. The exponential martingale $\mathcal{E}(\Mt_t)$ is clearly non-negative local martingale, and  it follows from Fatou's lemma for conditional expectation that $\mathcal{E}(\Mt_t)$ is a supermartingale, so, by Doob's theorem (e.g.,~\cite{Yor, Kunitabook}), it converges almost surely. Since, $
\mathcal{E}(\Mt_t) = \mathcal{E}\left(\frac{1}{2}\Mt_t\right)\exp\left(-\frac{1}{4}\langle \Mt\rangle_t\right)$, 
we have 
\begin{align*}
\left\{\om: \langle \Mt\rangle_{\infty} = \infty\right\}\subseteq\left\{\om: \lim_{t\rightarrow\infty}\mathcal{E}(\Mt_t) =0\right\}, \quad \p\text{\,-\,a.s.}
\end{align*}
Conversely, the identity 
$\mathcal{E}(-\Mt_t) = \mathcal{E}(\Mt_t)^{-1}\exp\left(\langle \Mt\rangle_t\right)$, yields
\begin{align*}
\hspace{3.5cm} \left\{\om: \lim_{t\rightarrow\infty}\mathcal{E}(\Mt_t) = \infty\right\}\subseteq \left\{\om: \langle \Mt\rangle_{\infty} = \infty\right\}.\hspace{4cm}\qed
\end{align*}

\subsection{Proof of Proposition \ref{DifFO}}\label{DifFO_app}
Fix $x\in \Xt$  and let $\tau_n = \inf\{t\geqslant t_0: \vert \phi^\nu_{t_0,t}\vert\geqslant n\}\wedge(t_0+T)$, with $ X^\nu_t(\om) = \phi^\nu_{t_0,t}(x,\om)$, $\p$\,-\,a.a $\om$, solving (\ref{SDE2}) for $t\in \Ic = [t_0,t_0+T]$ such that $X^\nu_{t_0}=x$.  Then $\tau_n$ is a localising sequence for the local martingale  
\begin{align*}
\Mt_t(x,\om) = \int_{t_0}^t\Big\langle \big({a}^{\nu,-1}(\tilde{b}-b^\nu)\big)\big(s,  \phi^\nu_{t_0,s}(x,\om)\big), d\St_s\Big\rangle-\int_{t_0}^t\left\langle {a}^{\nu,-1}(\tilde{b}-b^\nu), \tilde{b} \right\rangle\big(s, \phi^\nu_{t_0,t}(x,\om)\big)ds,
\end{align*}
where $\tilde{b} = \brr^\nu+\betaLk$, with $\betaLk$ as defined in (\ref{Per_vF}), $a^{\nu,-1}$ is the inverse of $a^\nu = \sigma^\nu(\sigma^\nu)^*$, which is coercive by assumption,  and $\St_t$ denotes an $\Ns_{t_0,x}^\nu$ martingale, $\Ns_{t_0,x}^\nu\in \PP(\mathcal{W}_d)$, given by 
\begin{align}\label{Mbb}
\begin{cases} \hspace{.4cm}\St_t(x,\om) =  \phi^\nu_{t_0,t}(x,\om)-x -\int_{t_0}^t b^\nu\big(s,\phi^\nu_{t_0,s}(x,\om)\big)ds, \quad t\in\Ic,\\[.2cm]
\big\langle \St\big\rangle_t(x,\om) = \int_{t_0}^ta^\nu\big(s, \phi_{t_0,s}^\nu(x,\om)\big)ds, \hspace{3.2cm} t\in\Ic.
\end{cases}
\end{align}
Next, for fixed $n\in\N,$ $\Mt_{\tau_n\wedge t}$ is an $\Ns_{t_0,x}^\nu$-\,adapted martingale (e.g.~\cite{Pinsky04, Yor}) and the corresponding quadratic variation process $\langle \Mt\rangle_t$ is given by 
\begin{align}\label{Recon}
\langle \Mt\rangle_t(x,\om) &= \int_{t_0}^t\left\langle \tilde{b}-b^\nu, {a}^{\nu,-1}(\tilde{b}-b^\nu)\right\rangle\big(s, \phi_{t_0,s}^\nu(x,\om)\big)ds\notag \\
&= \int_{t_0}^t\langle \betaL, a^\nu\betaL\rangle\big(s,\phi_{t_0,s}^\nu(x,\om)\big)ds,
\end{align}
where $\betaL = {a}^{\nu,-1}\betaLk$. Now, consider the following two probability measures  on $\mathcal{W}_d$ given by $\Ps_{t_0,x}^{\mu,\tau_n} :=\Ps_{t_0, x}^\mu|_{\F_{\tau_n}}$, and $\Ns_{t_0,x}^{\nu,\tau_n}:=\Ns_{t_0,x}^\nu|_{\F_{\tau_n}}$. Then, by the Girsanov theorem (e.g., \cite{Pinsky04, Yor}), $\Ps_{t_0,x}^{\mu,\tau_n}\ll \Ns_{t_0,x}^{\nu,\tau_n}$ and $\Ns_{t_0,x}^{\nu,\tau_n}\ll \Ps_{t_0,x}^{\mu,\tau_n}$ with 
\begin{align*}
\frac{\vphantom{\scaleobj{1}{|}}d\Ms^{\mu,\tau_n}_{t_0,x}}{\vphantom{\scaleobj{1.3}{|}}d\Ns^{\nu,\tau_n}_{t_0,x}} = \frac{1}{\mathcal{E}(-\Mt_{\tau_n})}, \qquad \Ns_{t_0,x}^\nu \; \text{-\,a.s.}
\end{align*}
This implies that
\begin{align*}
\left\{\om: \limsup_{n\rightarrow\infty}\frac{d\Ms_{t_0,x}^{\mu,\tau_n}}{\vphantom{\scaleobj{1.3}{|}}d\Ns^{\nu,\tau_n}_{t_0,x}}=\infty\right\}= \bigg\{\om: \lim_{n\rightarrow\infty}\mathcal{E}(-\Mt_{\tau_n}) = 0 \bigg\}, \qquad \Ns_{t_0,x}^\nu \, \text{-a.s.},
\end{align*}
and, according to Lemma \ref{Martinest}, this implies  
\begin{align*}
\bigg\{\om: \limsup_{n\rightarrow\infty}\frac{d\Ms_{t_0,x}^{\mu,\tau_n}}{\vphantom{\scaleobj{1.3}{|}}d\Ns^{\nu,\tau_n}_{t_0,x}}=\infty\bigg\} = \bigg\{\om: \langle \Mt\rangle_{\infty}=\infty\bigg\}.
\end{align*}
Then, by parts (a) and (b) of Lemma \ref{Absestm}, parts (i) and  (ii) of Proposition \ref{DifFO} hold.

\smallskip
For part (iii), define
\begin{align*}
D_s = \frac{d\Ms^{\mu,s}_{t_0,x}}{\vphantom{\scaleobj{1.3}{|}}d\Ns^{\nu,s}_{t_0,x}}, \qquad \Ns_{t_0,x}^\nu\; \text{-\,a.s.,}\quad s\in [t_0,\; t], \quad t\in\Ic,
\end{align*}
and   notice that the regularity of the convex function $\varphi$  implies that $t\mapsto\int_{t_0}^t\varphi^{\prime}(D_s)d\langle D\rangle_s$ is continuous and finite on the interval $\Ic$.

Next, consider a localising sequence $(T_n)_{n\in \N}$ defined by 
\begin{align*}
T_n = \inf\left\{t\in [t_0, t_0+T]: D_t\leqslant 1/n \; \text{or}\; \int_{t_0}^t\varphi^{\prime}(D_s)d\langle D\rangle_s\geqslant n\right\},
\end{align*}
so that $T_n\uparrow (t_0+T)$ as $n\rightarrow\infty.$ Next, by It\^o's formula, we have 
\begin{align*}
\varphi\left(D_{t\wedge\tau_n\wedge T_n}\right) = \varphi(D_{t_0})+\int_{t_0}^{t\wedge \tau_n\wedge T_n}\varphi^{\prime}\left(D_s\right)dD_s+\frac{1}{2}\int_{t_0}^{t\wedge \tau_n\wedge T_n}\varphi^{\prime\prime}(D_s)d\langle D\rangle_s.
\end{align*}
Since $\varphi(D_{t_0})<\infty$ (recall that $\varphi(D_{t_0}) = \varphi(1) =0$),  by  Jensen's inequality we have that $(\varphi(D_t))_{t\in\Ic}$ is a uniformly integrable submartingale w.r.t. $\Ns_{t_0,x}^\nu$. Taking the expectation of both sides w.r.t. $\Ns_{t_0,x}^\nu$, we have 
\begin{align*}
\E^{\Ns_{t_0,x}^\nu}\left[\varphi(D_{t\wedge\tau_n\wedge T_n})\right] =\frac{1}{2}\E^{\Ns_{t_0,x}^\nu}\left[\int_{t_0}^{t\wedge \tau_n\wedge T_n}\varphi^{\prime\prime}(D_s)d\langle D\rangle_s\right].
\end{align*}
As $\varphi$ is continuous and locally bounded on $\Rp$ due to its convexity and the normality condition (\ref{Normality}), we have that $\varphi(D_{t\wedge \tau_n\wedge T_n})\rightarrow \varphi(D_{t\wedge({t_0+T})}) $, as $n\rightarrow\infty$. By uniform integrability of $(\varphi(D_t))_{t\in\Ic}$ w.r.t. \!\!$\Ns_{t_0,x}^\nu$, we have that 
$\E^{\Ns_{t_0,x}^\nu}[\varphi(D_{t\wedge\tau_n\wedge T_n})]\rightarrow \E^{\Ns_{t_0,x}^\nu}[\varphi(D_t)]$ as $n\rightarrow\infty.$ Consequently, by the monotone convergence theorem, we have 
\begin{align*}
\E^{\Ns_{t_0,x}^\nu}[\varphi(D_t)] =  \frac{1}{2}\E^{\Ns_{t_0,x}^\nu}\left[\int_{t_0}^t\varphi^{\prime\prime}(D_s)d\langle D\rangle_s\right]\!, \quad t\in \Ic.
\end{align*} 
Finally, observe that $D_t$ is nothing but the exponential martingale of $ \Mt_t$ appearing in part (i) and (ii) (see~(\ref{Recon})), so that $d\langle D\rangle_s = D_s^2d\langle \Mt\rangle_s$. This fact together with the assumption  $\varphi\left(\frac{\scaleobj{1.1}{d\Ms_{t_0,x}^\mu}}{\scaleobj{1.1}{\vphantom{\scaleobj{1.3}{|}}d\Ns_{t_0,x}^\nu}}\right)\in L^1(\mathcal{W}_d; \Ns^\nu_{t_0,x})$ and part (c) of Lemma \ref{Absestm} yields for $t\in\Ic$
\begin{align*}
\;\;\;\Df\left(\Ms_{t_0,x}^\mu\|\Ns_{t_0,x}^\nu\right)&=
\lim_{n\rightarrow\infty}\E^{\Ns^\nu_{t_0,x}}\big[\varphi(D_{t\wedge\tau_n\wedge T_n})\big]\\ & = \E^{\Ns^\nu_{t_0,x}}\big[\varphi(D_t)\big]  = \frac{1}{2}\E^{\Ns^\nu_{t_0,x}}\bigg( \int_{t_0}^t\varphi^{\prime\prime}(D_s)D_s^2\big\langle \betaL, a^\nu\betaL\big\rangle(s,\phi^\nu_{t_0,s}(x)) ds\bigg). \;\;\qed
\end{align*}

\subsection{Proof of Theorem \ref{path_proj}}\label{path_proj_app}
First, we note that for a strictly convex function $\varphi\in \mathcal{C}^2(\Rp),$ $\varphi$-projections $\Ms^{\mu\nu}_{t_{0}}$ are uniquely determined since $\mathbb{C}_{\varphi,\Ic}^{{\mu\nu}}$ is a closed convex subset of $\PP(\Wd)$; this is a consequence of the Hahn--Banach theorem. Thus, there exists a unique $\Ps^{\mu\nu}_{t_{0}}\in \mathbb{C}_{\varphi,\Ic}^{{\mu\nu}}$ such that 
\begin{align}\label{phiproj}
\Df\big(\Ms^{\mu\nu}_{t_{0}}\|\Ns^\nu_{t_{0}}\big) = \inf\Big\{\Df(\Ms_{t_0}\|\Ns^\nu_{t_{0}}): \;\;\Ms_{t_0}\in \mathbb{C}_{\varphi,\Ic}^{{\mu\nu}}, \;\;\Ns_{t_{0}}^\nu\in\PP(\Wd)\Big\}.
\end{align}
It remains to show that there exists a Markovian version of $\Ms^{\mu\nu}_{t_{0}}\in \PP(\mathcal{W}_d)$  in the sense that for any $t_1,t_2\in \Ic$ with $t_0\leqslant t_1<t_2\leqslant t_0+T$, the conditional law of $X^\nu_{t_2}(\om) = \phi^\nu_{t_0,t_1}(x,\om)$ w.r.t.~$\Ms^{\mu\nu}_{t_{0}}$ given the $\mathfrak{S}$-algebra $\F_{t_1}$, is only a function of $t_2$ and $(t_1,\phi^\nu_{t_0,t_1}(x,\om))$.  

\smallskip
In order to simplify notation in what follows,   we set $\Ms_{t_0}(d\om) :=\int_\Xt \Ms_{t_0,x}(d\om)\mu_{t_0}(dx)$, for any $\Ps_{t_{0}}\in \mathbb{C}_{\varphi,\Ic}^{{\mu\nu}}$, and $X^\nu_t(\om)=\phi^\nu_{t_0,t}(x,\om)$ solving the SDE (\ref{SDE2}), we adopt the following notation: $$\E^{{\Ms}_{t_{0}}}\left[ f(t,\phi^\nu_{t_0,t})\right]\equiv \int_\Om f(t,\phi^\nu_{t_0,t}(\om))\Ms_{t_0}(d\om)\equiv\int_\Om\int_\Xt f(t,\phi^\nu_{t_0,t}(x,\om))\Ms_{t_0,x}(d\om)\mu_{t_0}(dx).$$ 
First, we note that similar to the steps following Lemma \ref{THT_def}, and  the assumption of this theorem that  $\Xt=\XXt$, and that the  regularity of the coefficients of the SDE's (\ref{SDE1}) and (\ref{SDE2}) is 
\begin{itemize}[leftmargin=0.7cm]
\item[--] $\brr^\mu,\brr^\nu\in \mathcal{C}\big(\Ic; \tilde{\mathcal{C}}^{3,\delta}(\Xt;\Xt)\big)$, 
\item[--] $\sigma^\mu,\sigma^\nu$ have uniformly bounded right inverses, 
\item[--] the columns of $\sigma^\mu$ and $\sigma^\nu$ are s.t.~$\sigma_k^\mu, \sigma_{k'}^\nu\in \bar{\mathcal{C}}^{4,\delta}(\Xt;\Xt)$, $1\leqslant k\leqslant m$, $1\leqslant k'\leqslant m'$,
\end{itemize}
there exists a predictable process $\hat{\beta}_{\mu\nu}\,{:}\; \Ic\times\Wd\rightarrow\Xt$ s.t.\! 
 $\Ms_{t_{0}}\in \mathbb{C}_{\varphi,\Ic}^{{\mu\nu}}\subset\PP(\Wd)$ is a martingale solution of an SDE associated with the generator $\LG:=\frac{\partial}{\partial t}+\LG_t^\nu+ a^\nu\hat{\beta}_{\mu\nu}\nabla$, $a^\nu = \sigma^\nu(\sigma^\nu)^*$, on the domain $\mathcal{C}_c^\infty(\Ic\times\Xt)$; this follows from the Girsanov theorem (e.g., \cite{Yor}, or \cite[Proposition~3.5]{Catti95}). To show that $\Ms_{t_{0}}{\in} \,\mathbb{C}_{\varphi,\Ic}^{{\mu\nu}}$  has a Markovian version, we only need to show the uniqueness of the conditional expectation of $\hat{\beta}_{\mu\nu}$ given $\F_t.$ To this end, we consider the pre-Hilbert spaces
 \begin{align*}
\mathbb{H}_{\mu} &:= \bigg\{ B\in \mathbb{M}_\infty(\Ic\times\Xt; \Xt): \left(\int_{\Ic}\int_\Xt \langle B,a^\nu B\rangle(s,x)\mu_s(dx)ds\right)^{1/2}<\infty\bigg\}
\end{align*}
and 
\begin{align*}
\mathbb{H}_{\Ms_{t_0}} :=\bigg\{\beta\in \mathbb{M}_\infty(\Ic\times\Om; \Xt): \left(\int_{\Ic}\int_\Om \langle\beta(s,\om), a^\nu\big(s, \phi^\nu_{t_0,s}(\om)\big)\beta(s,\om)\rangle \Ms_{t_0}(d\om) ds\right)^{1/2}<\infty\bigg\}.
\end{align*}
It follows from Proposition \ref{DifFO} that $\hat{\beta}_{\mu\nu}(s, \phi^\nu_{t_0,s}(x,\om))$ is in $ \mathbb{H}_{\Ms_{t_0}}$. 

Now, let $\Ms_{t_{0}}\in \mathbb{C}_{\varphi,\Ic}^{{\mu\nu}}$ and let the predictable process $\hat{\beta}_{\mu\nu}\in \mathbb{H}_{\Ms_{t_0}}$ be given. Then, by the Riesz representation theorem (e.g.,~\cite{Rudin}), there exists a unique $B\in \mathbb{H}^{-1}_{\mu}$ such that for all $f\in \mathbb{H}^{-1}_{\mu}$,  
\begin{align}\label{Riesz}
\int_{\Ic\times\Om} \big\langle\hat{\beta}_{\mu\nu}(s, \phi^\nu_{t_0,s}(\om)), (a^\nu f)(s, \phi^\nu_{t_0,s}(\om))\big\rangle \Ms_{t_{0}}(d\om)ds = \int_{\Ic\times\Xt}\langle B,a^\nu f\rangle(s,x)\mu_s(dx)ds,
\end{align}
where $\mathbb{H}_\mu^{-1} := \text{closure of } \Big\{ B\in \mathbb{M}_{\infty}(\Ic\times\Xt; \Xt): B = \nabla_x f, \; f\in \mathcal{C}_c^\infty(\Ic\times\Xt)\Big\}$ in $\mathbb{H}_\mu.$
This implies that $B$ is the unique Markovian version of $\hat{\beta}_{\mu\nu}$, i.e., $\E^{\Ms_{t_{0}}}\big[\hat{\beta}_{\mu\nu}| \F_s\big] = B\big(s, \phi^\nu_{t_0,s}\big)$; we denote the Markovian version of $\Ms_{t_{0}}$ by $\widehat{\Ms}_{t_{0}}$.  The first part of the proof will be complete if we justify the following claim:\\[.3cm]
{\bf Claim I:}  Let $\widehat{\Ms}_{t_{0}}\in \PP(\mathcal{W}_d)$ be the Markovian version of $\Ms_{t_{0}}\in \mathbb{C}_{\varphi,\Ic}^{{}\mu\nu}$. Then $\widehat{\Ms}_{t_{0}}\in\mathbb{C}_{\varphi,\Ic}^{\mu\nu}$ and it is such that  $\Df\big(\widehat{\Ms}_{t_{0}}\| \Ns_{t_{0}}^\nu\big)\leqslant\Df(\Ms_{t_{0}}\| \Ns_{t_{0}}^\nu)$ for all $\Ms_{t_{0}}\in \mathbb{C}_{\varphi,\Ic}^{{}\mu\nu}$.\\

\noindent
The difficult part of the claim is to show that $\widehat{\Ms}_{t_{0}}\circ \phi^{\nu,-1}_{t_0,t} = \mu_t, \; t\in \Ic$. This is achieved by first establishing the  following domination property
\begin{align}\label{Domination}
\E^{\widehat{\Ms}_{t_{0}}}\left[ f(t,\phi^\nu_{t_0,t})\right]\leqslant \int_{\Xt}f(t,x)\mu_t(dx), \quad t\in\Ic, \; f\in \mathbb{M}_\infty(\Ic\times\Xt).
\end{align}
In order to derive the above inequality, we follow the same localisation procedure as the one used in \cite{Catti94, Catti95}  in the case of KL-divergence. First,  it is clear from (\ref{Riesz}) and the $L^2$-contractivity of conditional expectation  that 
\begin{align}\label{DConex}
&\int_\Ic\int_{\Xt}\!\!\left\langle B(s,x), (a^{\nu}B)(s, x)\right\rangle\mu_s(dx)ds\notag\\
&\hspace{3cm}\leqslant \!\int_\Ic\int_{\Om^\mm}\!\!\left\langle \hat{\beta}_{\mu\nu}(s,\phi^\nu_{t_0,s}(\om)), (a^{\nu}\hat{\beta}_{\mu\nu})(s, \phi^\nu_{t_0,s}(\om))\right\rangle \Ms_{t_{0}}(d\om)ds<\infty.
\end{align}
We define localising sequences $(T_k)_{k\in \N}$ and $(S_k)_{k\in \N}$ as follows:
\begin{align*}
\begin{cases}
T_k = \inf\left\{t\in [t_0, t_0+T]: \quad \int_{t_0}^t\langle \hat{\beta}_{\mu\nu}, a^\nu\hat{\beta}_{\mu\nu}\rangle(s,  \phi^\nu_{t_0,s}(x,\om))ds>k  \right\}\wedge (t_0+T), \\[.2cm]
S_k = \inf\left\{t\in [t_0, t_0+T]: \quad \int_{t_0}^t\langle B, a^\nu B\rangle(s,  \phi^\nu_{t_0,s}(x,\om)) ds>k\right\}\wedge (t_0+T).
\end{cases}
\end{align*}
 Next, consider the Radon-Nikodym derivative
\begin{align}\label{DtTk}
\frac{d\Ms_{t_{0}}^{t\wedge T_k}}{\vphantom{\scaleobj{1.5}{|}}\;\;\;d\Ns_{t_{0}}^{\nu,t\wedge T_k}}=: D_{t\wedge T_k}, 
\end{align}
where $\Ms_{t_{0}}^{t\wedge T_k} :=\Ms_{t_{0}}|_{\F_{t\wedge T_k}}$, $\Ns_{t_{0}}^{\nu,t\wedge T_k} :=\Ns_{t_{0}}^\nu|_{\F_{t\wedge T_k}}$. Then, by Novikov's criterion (or~Lemma~\ref{Absestm}), we have that $\Ms^{t\wedge T_k}_{t_{0}}\in \PP(\Wd)$,  which leads to the following  for $f\in\mathbb{M}_\infty(\Ic\times\Xt)$ 
\begin{align}\label{Step0}
&\E^{\Ms^{t\wedge T_k}_{t_{0}}}[f(t,\phi^\nu_{t_0,t})] = \E^{\Ns^{\nu,t\wedge T_k}_{t_{0}}}\left[f(t,\phi^\nu_{t_0,t})D_{t\wedge T_k}\right]\notag\\[.2cm]
&\hspace{5.5cm}\underset{k\rightarrow \infty}\longrightarrow \E^{\Ns^{\nu}_{t_{0}}}\left[f(t,\phi^\nu_{t_0,t})D_t\right] = \E^{\Ms_{t_{0}}}[f(t,\phi^\nu_{t_0,t})].
\end{align}
Next, let $B_{(n,m)}=(B^i_{(n,m)})_{i\geqslant 1}$ be a sequence Borel functions defined by 
\begin{align*}
\begin{cases}
B_{(n,m)}^i(s,x)= B^i(s,x)\I_{\{\vert B^i(s,x)\vert\leqslant n\}}\I_{\{\max_{1\leqslant j\leqslant n}\vert a^\nu_{ij}|\leqslant m\}}, \quad &\text{if}\; i\leqslant n,\\
B^i_{(n,m)}(s,x) = 0,  &  \text{if}\; i\geqslant n,
\end{cases}
\end{align*}
where $B^i\in\mathbb{M}_\infty(\Ic{\times}\Xt; \Xt)$. For each $(n,m)\in \N\times\N$, we see that $B_{(n,m)}$ is bounded Borel measurable and $\langle B_{(n,m)}, a^\nu B_{(n,m)}\rangle$ is bounded. 
This implies that the following probability measure $\tilde \Ms^{t\wedge T_k}_{t_0(n,m)}=\int_\Xt \tilde \Ms^{t\wedge T_k}_{(t_0,x)(n,m)}\mu_{t_0}(dx)$ such that  
\begin{align*}
&\hspace{.5cm}\tilde \Ms^{t\wedge T_k}_{(t_0,x)(n,m)}(d\om)\mu_{t_0}(dx)= \mathcal{E}(t\wedge T_k,x,\om)\,\Ns^{\nu}_{t_{0},x}(d\om)\nu_{t_0}(dx),\\[.2cm] 
&\hspace{.5cm}\mathcal{E}(t\wedge T_k,x,\om) = \exp\bigg( \sum_{k\geqslant 1}\int_{t_0}^{t\wedge T_k}B_{(n,m)}^kd\tilde{\Mt}^k_s-\frac{1}{2}\int_{t_0}^{t\wedge T_k}\langle B_{(n,m)}, a^\nu B_{(n,m)}\rangle(s, \phi^{\nu}_{t_0,s}(x,\om))ds\bigg),
\end{align*}
 where $\tilde{\Mt}^k_t, \; k\geqslant 1$, are square-integrable $\Ns^{\nu}_{t_{0}}$ local martingales defined by 
\begin{align*}
\begin{cases}
\tilde{\Mt}_t^k(x,\om) := \phi^\nu_{t\wedge T_k,t_0}(x,\om)-x-\int_{t_0}^{t\wedge T_k}b^\nu(s, \phi^\nu_{t_0,s}(x,\om))ds,\\[.1cm]
\langle \tilde{\Mt}\rangle_{t\wedge T_k}(x,\om) = \int_{t_0}^{t\wedge T_k}a^{\nu}(s, \phi^\nu_{t_0,s}(x,\om))ds,
\end{cases}
\end{align*}
is well-defined.
 Next,  consider the difference
\begin{align}\label{Nstep}
&\int_{\Xt}f(t, x)\mu_t(dx)-\int_{\Xt}\E^{\tilde \Ms^{t\wedge T_k}_{t_0(n,m)}}[f(t,\phi^\nu_{t_0,t})]\mu_{t_0}(dx) \notag\\&\hspace{4cm}= \int_{\Ic\times\Xt}\Big\langle a^{\nu}( B-B_{(n,m)}), \nabla_x f\Big\rangle(s,x)\mu_s(dx)ds,
\end{align}
for  $f\in \mathcal{C}^{\infty}_c(\Ic\times\Xt)$.
Given that $\nabla_x f$ is bounded with a compact support,  by the Cauchy--Schwartz inequality there exists $C>0$ (cf.~\cite{Catti94}) such that (reverting to the shorter notation $X_t^\nu(\om) = \phi^\nu_{t_0,t}(x,\om)$) we have 
\begin{align*}
\Big\vert \int_{t_0}^t\langle a^\nu(B-B_{(n,m)}), \nabla_x f\rangle(s,X^\nu_s)ds\Big\vert ^{2}\leqslant C\Vert \nabla_x f\Vert_{\infty}^2\Big\vert \int_{t_0}^t\langle B-B_{(n,m)}, a^{\nu}(B-B_{(n,m)})\rangle(s,X^\nu_s)ds\Big\vert^2.
\end{align*}
Consequently, for   $f\in \mathcal{C}_c^\infty(\Ic\times\Xt)$  we have (see \cite[Section 4]{Catti94})
\begin{align*}
\E^{\tilde \Ms^{t\wedge T_k}_{t_0(n,m)}}\left[f(t,X^\nu_t)\right]\leqslant \int_{\Xt}f(t,x)\mu_t(dx)+C\Vert f\Vert_{\infty}\Vert B-B_{(n,m)}\Vert^2_{\mathbb{H}_\mu}\left(1+ \Vert B-B_{(n,m)}\Vert^2_{\mathbb{H}_\mu}\right).
\end{align*}
Then, for any $f\in\mathbb{M}_\infty(\Ic\times\Xt)$,  we use the density argument to arrive at 
\begin{align}\label{Nstep2}
\E^{\tilde{\Ms}^{t\wedge T_k}_{t_0(n,m)}}\left[f(t,X^\nu_t)\right]\leqslant \int_{\Xt}f(t,x)\mu_t(dx)+C\Vert f\Vert_{\infty}\Vert B-B_{(n,m)}\Vert^2_{\mathbb{H}_\mu}\left(1+ \Vert B-B_{(n,m)}\Vert^2_{\mathbb{H}_\mu}\right).
\end{align}
 Next, it is relatively straightforward to verify via the variational representation of $\Df$ that $\lim_{n,m\rightarrow\infty}\Df(\tilde \Ms^{t\wedge T_k}_{t_0(n,m)}\| \Ms^{t\wedge T_k}_{t_{0}})\leqslant 2 \Df(\Ms_{t_{0}}\|\Ns^{\nu}_{t_{0}})<\infty$.  Then,  by the $\varphi$-divergence formula in part~(iii) of Proposition \ref{DifFO} and (\ref{Step0}), we obtain 
\begin{align*}
\notag &\hspace{-0.cm}\lim_{n\rightarrow \infty}\lim_{m\rightarrow\infty}\Df(\tilde \Ms_{t_0(n,m)}^{t\wedge T_k}\|\Ms^{t\wedge T_k}_{t_{0}}) \\ &\hspace{.5cm}= \lim_{n\rightarrow\infty}\lim_{m\rightarrow \infty}\E^{\Ns^{\nu}_{t_{0}}}\!\!\left[\int_{t_0}^{t\wedge S_k}\!\!\!\varphi^{\prime\prime}(G^{n,m}_{s})(G^{n,m}_{s})^2D_{s}\big\langle B\,{-}\,B_{(n,m)}, a^{\nu}(B\,{-}\,B_{(n,m)})\big\rangle(s,X^\nu_s)ds\right]{=}\,0,
\end{align*} 
where $D_s$ is as defined in (\ref{DtTk}), and $G_{s}^{n,m} = d\Ms^{s}_{t_0(n,m)}/d\Ms^{s}_{t_{0}}$.  It follows that $(D_{t\wedge S_k}, B_{(n,m)}, \mu_{t_0})\rightarrow (D_{t\wedge S_k}, B_{(n,\infty)}, \mu_{t_0})$ as $m\rightarrow \infty$ $\Ns^{\nu}_{t_{0}}\,\text{-\,a.s.}$ and there exists a subsequence $(D_{t\wedge S_k}, B_{(n_{j}, \infty)}, \mu_{t_0})$ s.t.~$(D_{t\wedge S_k}, B_{(n_{j}, \infty)}, \mu_{t_0})\rightarrow (D_{t\wedge S_k}, B, \mu_{t_0})$ as $j\,{\rightarrow}\, \infty, \; \Ns^{\nu}_{t_{0}}\,\text{-\,a.s.}$
Then, applying Fatou's lemma twice, we have for any $f\in\mathbb{M}^{+}(\Ic\times\Xt)$,
\begin{align}\label{Fatou_Step}
\notag \E^{\widehat{\Ms}_{t_{0}}}\left[ f(t,X^\nu_t)\I_{\{t<S_k\}}\right]&= \E^{\Ns^\nu_{t_{0}}}\left[f(t,X^\mu_t)D_{t\wedge S_k}\right]\leqslant \liminf_{n\rightarrow\infty}\liminf_{m\rightarrow\infty}\;\E^{\tilde \Ms_{t_0(n,m)}^{t\wedge S_k}}\left[ f(t,X^\nu_t)\I_{\{t<S_k\}}\right]\\[.2cm]
&\hspace{0cm}\leqslant  \liminf_{n\rightarrow\infty}\liminf_{m\rightarrow\infty}\;\E^{\tilde \Ms_{t_0(n,m)}^{t\wedge S_k}}\left[ f(t,X^\nu_t)\right]\leqslant \int_{\Xt}f(t,x)\mu_t(dx), 
\end{align}
where the last step follows from (\ref{Nstep2}) and the fact that
\begin{align*}
\lim_{n\rightarrow\infty}\lim_{m\rightarrow\infty}\Vert B-B_{(n,m)}\Vert_{\mathbb{H}_\mu}= 0
\end{align*}
by the bounded convergence theorem.

 Combining (\ref{Nstep}), (\ref{Nstep2}) and (\ref{Fatou_Step}) together with application of the monotone convergence theorem, we arrive at the domination property 
\begin{align}\label{Dominat}
\E^{\widehat{\Ms}_{t_{0}}}\left[ f(t,X^{\nu}_{t})\right]\equiv\E^{\widehat{\Ms}_{t_{0}}}\left[ f(t,\phi^{\nu}_{t_0,t})\right]\leqslant \int_{\Xt}f(t,x)\mu_t(dx).
\end{align}
Since $\widehat{\Ms}_{t_{0}}\circ \phi^{\nu,-1}_{t_0,t}$ and $\mu_t$ are probability measures on $\Xt$,  the domination property (\ref{Dominat}) implies that $\widehat{\Ms}_{t_{0}}\circ \phi^{\nu,-1}_{t_0,t} = \mu_t$ and hence, $\widehat{\Ms}_{t_{0}}\in \mathbb{C}_{\varphi,\Ic}^{\mu\nu}$. This conclusion follows from the property of $\Df$, i.e., given the domination property (\ref{Dominat}), we have from the variational representation of $\Df$,
 \begin{align*}
\Df(\widehat{\Ms}_{t_{0}}\circ \phi^{\nu,-1}_{t_0,t}\|\mu_t) &= \sup_{f(t,\ccdot)\in \mathcal{C}_\infty(\Xt)}\bigg\{\int_{\Xt}f(t,x)(\widehat{\Ms}_{t_{0}}\circ \phi^{\nu,-1}_{t_0,t})(dx)- \int_{\Xt}\varphi^*(f(t,x))\mu_t(dx)\bigg\}\\
&=  \sup_{f(t,\ccdot)\in \mathcal{C}_\infty(\Xt)}\bigg\{\E^{\widehat{\Ms}_{t_{0}}}[f(t,\phi^{\nu}_{t_0,t})]- \int_{\Xt}\varphi^*(f(t,x))\mu_t(dx)\bigg\}\\[.2cm]
& \leqslant \Df(\mu_t\|\mu_t) =0.
\end{align*}
The second assertion, $\Df(\widehat{\Ms}_{t_{0}}\|\Ns^{\nu}_{t_{0}})\leqslant \Df({\Ms}^{\mu\nu}_{t_{0}}\| \Ns^{\nu}_{t_{0}})$,  of Claim I  is obtained by noticing that 
\begin{align}\label{Uniqn}
\hspace{-.1cm} \notag \Df(\widehat{\Ms}_{t_{0}}\|\Ns^{\nu}_{t_{0}}) &= \hspace{-.1cm}\sup_{f\in\mathbb{M}_\infty(\Wd)}\left\{\int_{\Wd} fd\widehat{\Ms}_{t_{0}}-\!\int_{\Wd}\varphi^*(f)d\Ns^{\nu}_{t_{0}}\right\}\\ 
&=\hspace{-.1cm} \sup_{f\in \mathcal{C}_\infty(\Wd)}\left\{\int_{\Wd} fd\widehat{\Ms}_{t_{0}}-\!\int_{\Wd}\varphi^*(f)d\Ns^{\nu}_{t_{0}}\right\} \leqslant \Df({\Ms}_{t_0}\|\Ns^{\nu}_{t_{0}}), \; \;\;\forall \;{\Ms}_{t_0}\in \mathbb{C}_{\varphi,\Ic}^{\mu\nu},
\end{align} 
where we have used the fact that $\Ms_{t_0}\circ \phi^{\nu,-1}_{t_0,t} = \mu_t = \widehat{\Ms}_{t_0}\circ \phi^{\nu,-1}_{t_0,t}$. 
Note also that (\ref{phiproj}) and the inequality~(\ref{Uniqn}) imply that  $\widehat{\Ms}_{t_0}=\Ms^{\mu\nu}_{t_0}$ so that the unique projection  of $\Ns^{\nu}_{t_{0}}\in \PP(\mathcal{W}_d)$ onto the closed and convex non-empty set $\mathbb{C}_{\varphi,\Ic}^{\mu\nu}$  has a Markovian version.

\smallskip
Now, we derive the representation of $\Df(\Ms^{\mu\nu}_{t_{0}}\| \Ns^\nu_{t_{0}})$ where $\Ms^{\mu\nu}_{t_{0}}, \Ns^\nu_{t_{0}}\in \mathbb{C}_{\varphi,\Ic}^{\mu\nu} \subset \PP(\mathcal{W}_{d})$.  To this end, consider the Polish spaces $\big(\Xt,\mathcal{B}({\Xt})\big)$ and $\big(\Ic,\mathcal{B}({\Ic}),\rho\big)$, $\Ic:=[t_{0}, \;t_{0}+T]$, where $\rho$ is a nonnegative measure on $\mathcal{B}(\Ic)$. The set of paths in $\mathcal{W}_{d}=C(\Ic,\Xt)$ is endowed with the relative $\mathfrak{S}$-field associated with $\mathcal{B}(\Xt)^{{\otimes \Ic}}$.  We utilise the primal-dual representation of constrained convex optimisation problem (e.g.,~\cite{Rockafeller74, Leonard01, Leonard01b}). The primal problem is given by
 \begin{align}\label{Primal2}
\Df\big(\Ms^{\mu\nu}_{t_{0}}\| \Ns^\nu_{t_{0}}\big) = \inf\left\{ \Df\big(\Ms_{t_{0}}\|\Ns^\nu_{t_{0}}\big): \Ms_{t_{0}}\ll \Ns^{\nu}_{t_{0}},\;\;\Ms_{t_{0}}\circ \phi^{\nu,-1}_{t_0,t} = \mu_t, \; \forall\,t\in \Ic\right\}.
\end{align}
The dual problem associated with the primal problem (\ref{Primal2}) is given by (e.g.,~\cite{Rockafeller74, Leonard01b} and, in particular, \cite[Proposition 6.2]{Leonard01} with appropriate notational adjustments)
\begin{align}\label{Dual2}
\notag \sup\bigg\{ \int_{\Ic\times\Xt}\sum_{i=1}^n\alpha_i(t)f_i(x)\mu_{t}(dx)\rho(dt) -\int_{\Wd}\varphi^*\Big(\int_\Ic\sum_{i=1}^n\alpha_i(t)f_i(\phi^{\nu}_{t_0,t}(\om))\rho(dt)\Big) \Ns_{t_{0}}^{\nu}(d\om):\\ \forall \;1\leqslant n <\infty ,f_{i}\in \mathcal{C}_\infty(\Xt), \;   \alpha_i\in L^{1}(\Ic,\rho), 1\leqslant i\leqslant n\bigg\}.
\end{align}
By Fenchel's primal-dual equality (e.g.,~\cite{Rockafeller74, Leonard01, Leonard01b} ), we have 
\begin{align*}
&\inf\left\{ \Df\big(\Ms_{t_{0}}\|\Ns^\nu_{t_{0}}\big): \Ms_{t_{0}}\ll \Ns^{\nu}_{t_{0}},\;\;\Ms_{t_{0}}\circ \phi^{\nu,-1}_{t_0,t} = \mu_t, \; \forall\,t\in \Ic\right\}\\
& \hspace{1cm}= \sup\bigg\{ \int_{\Ic\times\Xt}\sum_{i=1}^n\alpha_i(t)f_i(x)\mu_{t}(dx)\rho(dt) -\int_{\Wd}\varphi^*\Big(\int_\Ic\sum_{i=1}^n\alpha_i(t)f_i(\phi^{\nu}_{t_0,t}(\om))\rho(dt)\Big) \Ns_{t_{0}}^{\nu}(d\om);\\& \hspace{6cm}\forall \;1\leqslant n <\infty ,f_{i}\in \mathcal{C}_\infty(\Xt), \;   \alpha_i\in L^{1}(\Ic,\rho), 1\leqslant i\leqslant n\bigg\}.
\end{align*}
Since  $\Ns_{t_{0}}^{\nu}$ is an extremal solution to the martingale problem with generator $\mathcal{L}_t^\nu$, we have 
\begin{align}\label{PD3}
&\inf\left\{ \Df\big(\Ms_{t_{0}}\|\Ns^\nu_{t_{0}}\big): \Ms_{t_{0}}\ll \Ns^{\nu}_{t_{0}},\;\;\Ms_{t_{0}}\circ \phi^{\nu,-1}_{t_0,t} = \mu_t, \; \forall\,t\in \Ic\right\}\notag\\
& \hspace{2cm}= \sup\bigg\{ \int_{\Ic\times\Xt}\sum_{i=1}^n\alpha_i(t)f_i(x)\mu_{t}(dx)\rho(dt) -\int_{\Xt}\varphi^*\Big(\int_\Ic\sum_{i=1}^n\alpha_i(t)f_i(x)\rho(dt)\Big) \nu_{t}(dx);\notag\\ &\hspace{5cm} \forall \;1\leqslant n <\infty ,f_{i}\in \mathcal{C}_\infty(\Xt), \;   \alpha_i\in L^{1}(\Ic,\rho), 1\leqslant i\leqslant n\bigg\}.
\end{align}
Moreover, if we consider the discrete topology on $\Ic\subset\R$ with $\rho$ the counting measure on $\Ic$, then (\ref{PD3}) can be written in a more general form \cite[Proposition 6.1]{Leonard01}
\begin{align}\label{PD4}
&\inf\left\{ \Df\big(\Ms_{t_{0}}\|\Ns^\nu_{t_{0}}\big): \Ms_{t_{0}}\ll \Ns^{\nu}_{t_{0}},\;\;\Ms_{t_{0}}\circ \phi^{\nu,-1}_{t_0,t} = \mu_t, \; \forall\,t\in \Ic\right\}\notag\\
& \hspace{2cm}= \sup\bigg\{ \int_{\Xt}\sum_{i=1}^n f_i(x)\mu_{s_i}(dx) -\int_{\Xt}\varphi^*\Big(\sum_{i=1}^n f_i(x)\Big) \nu_{s_i}(dx);\notag\\ &\hspace{5cm} \forall \;1\leqslant n <\infty ,\;f_{1}, \dots, f_{n}\in \mathbb{M}_\infty(\Xt), \;   s_1,\dots,s_n\in \Ic\bigg\},
\end{align}
where the supremum is taken over all $n$-tuple partitions of $\Ic$, and $n$ functions in $\mathbb{M}_\infty(\Xt)$ for all $1\leqslant n <\infty$.
The above procedure  simplifies further if one considers the $\varphi$-projection for finite-dimensional marginal probability measures  $\Ms_{\Ic_{n}}, \Ns^\nu_{\Ic_{n}}$ on $\otimes_{i=1}^nA_i$, $A_i\in \Bb(\Xt)$ in which case one obtains 
\begin{align}\label{PD5}
&\inf\left\{ \Df(\Ms_{\Ic_n}\|\Ns^\nu_{\Ic_n}): \Ms_{\Ic_n}\ll \Ns^{\nu}_{\Ic_n},\;\;\Ms_{t_0}\circ \phi^{\nu,-1}_{t_i,t_0} = \mu_{t_i}, \;\;\; t_0, t_1, \dots, t_n\in \Ic\right\}\notag\\
& \hspace{2cm}= \sup\bigg\{ \int_{\Xt}\sum_{i=1}^n f_i(x)\mu_{t_i}(dx) -\int_{\Xt}\varphi^*\Big(\sum_{i=1}^n f_i(x)\Big) \nu_{t_i}(dx);\notag\\ &\hspace{8cm}\;f_{1}, \dots, f_{n}\in \mathbb{M}_\infty(\Xt), \;   t_1,\dots,t_n\in \Ic\bigg\}.
\end{align}
The second part of the theorem follows from the observation that the normality conditions (\ref{Normality}) imply that $\varphi^*(0) =0,$ which in particular implies that $\varphi^*$ is a super-additive convex function, so the right-hand side of (\ref{PD4}) can bounded by the sum of the difference of FTDR fields, as in Theorem \S\ref{difference_bound} of \S\ref{ftdr_sec}. To see this, the supper-additivity of $\varphi^*$ implies that $$-\varphi^*\left(\sum_{i=1}^n f_i(x)\right)\leqslant -\sum_{i=1}^n\varphi^*(f_i(x)),$$ so that from (\ref{PD4}), we have 
\begin{align*}
&\hspace{0cm} \inf\left\{ \Df\big(\Ms_{t_{0}}\|\Ns^\nu_{t_{0}}\big): \Ms_{t_{0}}\ll \Ns^{\nu}_{t_{0}},\;\;\Ms_{t_{0}}\circ \phi^{\nu,-1}_{t_0,t} = \mu_t, \; \forall\,t\in \Ic\right\}\notag\\
& \hspace{3.5cm}\leqslant \sup\bigg\{\sum_{i=1}^n\int_{\Xt}f_i(x)\mu_{s_i}(dx)-\sum_{i=1}^n\int_{\Xt}\varphi^*(f_i(x))\nu_{s_i}(dx);\notag\\ 
&\hspace{5.4cm} \forall \;1\leqslant n <\infty ,\;f_{1}, \dots, f_{n}\in \mathbb{M}_\infty(\Xt), \;   s_1,\dots,s_n\in \Ic\bigg\}\\
&\hspace{3.5cm}\leqslant \sup\bigg\{\sum_{i=1}^n\Big\vert \Df(\mu_{s_i}\|\mu_{t_0})-\Df(\nu_{s_i}\|\mu_{t_0})\Big\vert, \;\forall \;1\leqslant n <\infty , \;   s_1,\dots,s_n\in \Ic\bigg\},
\end{align*}
where the supremum is taken over all $n$-tuple partitions of $\Ic$ for all $1\leqslant n <\infty$. Analogously, for finite-dimensional marginal probability measures  $\Ms_{\Ic_{n}}, \Ns^\nu_{\Ic_{n}}$ on $\otimes_{i=1}^nA_i$, $A_i\in \Bb(\Xt)$ we obtain 
 \begin{align*}
&\hspace{0cm} \inf\left\{ \Df\big(\Ms_{\Ic_n}\|\Ns^\nu_{\Ic_n}\big): \Ms_{\Ic_n}\ll \Ns^{\nu}_{\Ic_n},\;\;\Ms_{t_{0}}\circ \phi^{\nu,-1}_{t_0,t} = \mu_{t_i}, \; \,t_0,t_1,\dots,t_n\in \Ic\right\}\notag\\
& \hspace{5.5cm}\leqslant \sup\bigg\{\sum_{i=1}^n\int_{\Xt}f_i(x)\mu_{t_i}(dx)-\sum_{i=1}^n\int_{\Xt}\varphi^*(f_i(x))\nu_{t_i}(dx);\notag\\ 
&\hspace{8.7cm}  \;f_{1}, \dots, f_{n}\in \mathbb{M}_\infty(\Xt), \;   t_1,\dots,t_n\in \Ic\bigg\}\\
&\hspace{5.5cm}\leqslant \sum_{i=1}^n\Big\vert \Df(\mu_{t_i}\|\mu_{t_0})-\Df(\nu_{t_i}\|\mu_{t_0})\Big\vert,  \;\;\;   t_1,\dots,t_n\in \Ic.
\end{align*}
Combining the above with (\ref{Primal2}) leads to 
\begin{align*}
\hspace{2.1cm}\Df\big(\Ms_{t_0}^{\mu\nu}\|\Ns_{t_0}^{\nu}\big) &\leqslant  \sum_{i=1}^n \Big\vert \Df(\mu_{t_i}\|\mu_{t_0})- \Df(\nu_{t_i}\|\mu_{t_0})\Big\vert, \qquad   t_1,\dots,t_n\in \Ic. \hspace{2cm}\qed
\end{align*}

\end{appendices}

\newpage

\vspace*{-1.2cm}
\section*{Glossary}\label{glossary}

\vspace{-.2cm}\noindent Here, we list further definitions and notation which recurs throughout the paper. 

\noindent{\small \bf (1) Domains. }{\small Throughout,  $\XXt$ is a finite-dimensional smooth manifold taken to be either $\XXt = \R^\ell$ or a \\ \hspace*{.5cm} flat torus $\XXt = \bar{\mathbb{T}}^\ell$, $1\leqslant \ell <\infty$. $\Xt\subseteq \XXt$ is a linear subspace, $\textrm{dim}(\Xt)=d$, $1\leqslant d \leqslant \ell$. }

\smallskip
\noindent{\small \bf (2) Probability spaces and function spaces}
{\small
\vspace*{-.0cm}
\begin{itemize}[leftmargin=0.4cm]
\item[\tiny$\bullet$]{\it Wiener space.} We  fix the probability space $\PS$ as the Wiener space, i.e.,
$\Om=\mathcal{C}_{0}(\Ic;\R^m)$, $m\in \N$, $\Ic:= [t_0, t_0+T]\subset\R$, $T>0$, is a subspace of continuous functions $\mathcal{C}(\R;\R^m)$ which are zero at $t_0\in \Ic$. $\F$ is the Borel  $\mathfrak{S}$-algebra generated by open subsets  in the compact-open topology on $\Om$ defined~via
\begin{align*}
\varrho(\om,\hat{\om})= \sum_{\ell=0}^\infty\frac{1}{2^\ell}\frac{\Vert \om -\hat{\om}\Vert_\ell}{1+\Vert \om -\hat{\om}\Vert_\ell}, \qquad \Vert \om -\hat{\om}\Vert_{\ell} := \sup_{t\in [-\ell, \,\ell\,]}\vert \om(t) -\hat{\om}(t)\vert, \quad \om,\hat\om\in \Omega, 
\end{align*}
with $|{\ccdot}|$ the Euclidean norm on $\R^m$. Finally, $\p$ is the Wiener measure on $\mathcal{F}$.

\vspace{.05cm}
\item[\tiny$\bullet$]  $\mathcal{W}_\ell:=\mathcal{C}(\Ic,\XXt)$ and $\mathcal{W}_d:=\mathcal{C}(\Ic,\Xt)$ denote path spaces defined over $\XXt$ and $\Xt$, respectively. Borel $\mathfrak{S}$-algebras, $\mathcal{B}(\mathcal{W}_\ell)$ and $\mathcal{B}(\mathcal{W}_d)$, on $\mathcal{W}_\ell$ and $\mathcal{W}_d$ are defined analogously to those in the Wiener space.

\vspace{.05cm}
\item[\tiny$\bullet$] $(\Om_\mm, \mathcal{H}_\mm, \Ms^\mm_{t_0})$ is a probability space, $\Om_\mm\simeq \mathcal{W}_\ell$, $\mathcal{H}_\mm\simeq \mathcal{B}(\mathcal{W}_\ell)$; \!and  $(\Om_\nu, \mathcal{H}_\nu, \Ns^\nu_{t_0})$,  $\Om_\nu\simeq \mathcal{W}_d$, $\mathcal{H}_\nu\simeq \mathcal{B}(\mathcal{W}_d)$.

\vspace{.05cm}
\item[\tiny$\bullet$] For $f\,{:}\; \XX\rightarrow\R$, where $(\XX, \mathcal{B}(\XX))$ is a Polish  space equipped with a Borel $\mathfrak{S}$\,-\,algebra, 
the following function spaces are relevant  (in particular,  $\XX = \XXt$ or $\XX = \Xt$, or $\XX = \mathcal{W}_d$):
\begin{itemize}[leftmargin=0.4cm]
 \item $\mathbb{M}_\infty(\XX)$ space of bounded Borel measurable functions  $\mathbb{M}(\XX)$ on $\XX$. 
 \item $\mathbb{M}^{+}(\XX)$ space of non-negative Borel measurable functions on $\XX$.
 \item $\mathcal{C}_\infty(\XX)$ space of bounded continuous functions on $\XX$.
\item  $\mathcal{C}^l(\XX)$, $l\geqslant 1$, space of  $l$-times continuously differentiable functions on $\XX$. 
\item  $\mathcal{C}^k_\infty(\XX)$,  $l\geqslant 1$, functions in $\mathcal{C}^l(\XX)$ which are bounded with bounded derivatives up to order $l$ on $\XX$. 
 \item $\mathcal{C}_c^{+}(\XX)$ space continuous non-negative functions on $\XX$ with compact supports.
\item  $\mathcal{C}_c^{\infty}(\XX)$ space of smooth functions on $\XX$ with compact support.
 \end{itemize}
  
\vspace{.05cm}
\item[\tiny$\bullet$]  Given the Borel measure space  $\big(\XX, \Bb(\XX), m_n\big)$, $\dim \XX<\infty$, with $m_n$ denoting the Lebesgue measure on $\XX$,  the Banach space $L^p(\XX;m_n)$ is the set of  Lebesgue-integrable functions satisfying
\begin{align*}
L^p\big(\XX;m_n\big):&=\big\{f\in \mathbb{M}(\XX): \Vert f\Vert_p<\infty\big\}, \hspace{.85cm} \Vert f\Vert_p:=\left({\textstyle \int_{\XX}}\vert f\vert^p d m_n\right)^{1/p}, \; 1\leqslant p<\infty,\\
L_{+}^p\big(\XX;m_n\big) :&= \big\{f\in L^p\big(\XX;m_n\big): f>0\big\},\\
L^\infty\big(\XX;m_n\big) :&= \big\{f\in \mathbb{M}(\XX): \Vert f\Vert_\infty<\infty\big\},  \qquad \Vert f\Vert_\infty:=\inf\big\{C\geqslant 0: |f|\leqslant C \;\;m_n \,\textrm{-a.e.} \big\}.
\end{align*}

\item[\tiny$\bullet$] $\tilde{\mathcal{C}}^{l,\delta}(\XX;\XX)$,  is the space of functions $f\,{:} \;\XX\,{\rightarrow}\, \XX$ with the countable family of semi-norms
\begin{align*}
& \tilde{\Vert} f\tilde{\Vert}_{l,\delta;N}:= \tilde{\Vert} f\tilde{\Vert}_{l;N}+\sum_{\vert \alpha\vert =l}\sup_{x,y\in \textsf{B}_{N}\!,\, x\neq y}\frac{\vert D^{\alpha}f(x) -D^{\alpha}f(y)\vert}{\vert x-y\vert^{\delta}}<\infty, \quad 0<\delta\leqslant 1,\;N\in\mathbb{N}_1,\hspace{.6cm}\\[.1cm]
&\tilde{\Vert} f\tilde{\Vert}_{l;N}:= \sup_{x\in\XX}\frac{\vert \langle f(x),x\rangle\vert}{1+\vert x\vert^2}+\sum_{1\leqslant \vert \alpha\vert \leqslant l}\sup_{x\in \textsf{B}_{N}}\vert D^{\alpha}f(x)\vert,
\end{align*}

\vspace{-.3cm}\noindent where $\textsf{B}_N:=\{x\in\XX: |x|\leqslant N\}$, $\displaystyle D^{\alpha}f(x): = \frac{\partial^{\vert \alpha\vert}f}{(\partial x_1)^{\alpha_1}\cdots(\partial x_d)^{\alpha_n}}, \; \vert \alpha\vert := \sum_{i=1}^n\alpha_i$, $\alpha_i\in \mathbb{N}_0$, and \;$D^0\equiv \textrm{Id}$.

\smallskip
\item[\tiny$\bullet$] $\mathcal{C}\big(\Ic; \tilde{\mathcal{C}}^{l,\delta}(\XX;\XX)\big)$, $\Ic\subseteq \R$,  is the set of all continuous fields $f\,{:}\; \R\times\XX\rightarrow\XX$ such that $f(t,\ccdot)\in \tilde{\mathcal{C}}^{l,\delta}(\XX;\XX)$. 

\smallskip
\item[\tiny$\bullet$] $\bar{\mathcal{C}}^{l,\delta}(\XX;\XX)$,  is the space of functions $f\,{:} \;\XX\,{\rightarrow}\, \XX$ with the countable family of semi-norms
\begin{align*}
& \tilde{\tilde{\Vert}} f\tilde{\tilde{\Vert}}_{l,\delta;N}:= \tilde{\tilde{\Vert}} f\tilde{\tilde{\Vert}}_{l;N}+\sum_{\vert \alpha\vert =l}\sup_{x,y\in \textsf{B}_{N}\!,\, x\neq y}\frac{\vert D^{\alpha}f(x) -D^{\alpha}f(y)\vert}{\vert x-y\vert^{\delta}}<\infty, \quad 0<\delta\leqslant 1,\;N_1\in\mathbb{N},\hspace{.6cm}\\[.1cm]
&\tilde{\tilde{\Vert}} f\tilde{\tilde{\Vert}}_{l;N}:= \sup_{x\in\XX}\frac{\vert  f(x)\vert}{1+\vert x\vert}+\sum_{1\leqslant \vert \alpha\vert \leqslant l}\sup_{x\in \textsf{B}_{N}}\vert D^{\alpha}f(x)\vert,
\end{align*}

\smallskip
\item[\tiny$\bullet$] $\mathcal{C}\big(\Ic; \bar{\mathcal{C}}^{l,\delta}(\XX;\XX)\big)$, $\Ic\subseteq \R$,  is the set of all continuous fields $f\,{:}\; \R\times\XX\rightarrow\XX$ such that $f(t,\ccdot)\in \bar{\mathcal{C}}^{l,\delta}(\XX;\XX)$. 

\end{itemize}

\newpage
\noindent{\small \bf (3) Frequently used notation}

\vspace*{-0.1cm}
\begin{itemize}[leftmargin=0.4cm]
\item[\tiny$\bullet$] $X^\mm_t(\omega)$, $\om\in \Om_\mm\simeq \Om$, $t\in \Ic := [t_0, t_0+T]$, is a solution of the SDE (\ref{SDE1}) with coefficients $(\brr^\mm$, $\sigma^\mm)$ representing  the original/reference dynamics on a smooth manifold~$\XXt$; $\textrm{dim}(\XXt)=\ell<\infty$.

\vspace{.1cm}
\item[\tiny$\bullet$] $X^\nu_t(\omega)$, $\om\in \Om_\nu\simeq \Om$, $t\in \Ic$, is a solution of the SDE (\ref{SDE2}) with coefficients $(\brr^\nu$, $\sigma^\nu)$ representing an approximation  of (\ref{SDE1}) on a subspace $\Xt\subseteq \XXt$; $\textrm{dim}(\Xt)=d\leqslant \ell$.

\vspace{.1cm}
\item[\tiny$\bullet$] $\mathcal{L}_t^{\mm}$ is the generator (\ref{gen}) of It\^o diffusion solving (\ref{SDE1}) on $\XXt$ with the $L^2(\XXt,\mm)$ dual  denoted by $\mathcal{L}_t^{\mm*}$.

\vspace{.1cm}
\item[\tiny$\bullet$] $\mathcal{L}_t^{\nu}$ is the generator of \^Ito diffusion solving (\ref{SDE2}) on $\Xt$ with the $L^2(\Xt,\nu)$ dual  denoted by $\mathcal{L}_t^{\nu*}$.

\vspace{.1cm}
\item[\tiny$\bullet$] $ {b}_i^\mm(t,\xxt) \,{:=}\, \brr_i^\mm(t,\xxt)\,{+}\,c_i^\mm(t,\xxt)$,  $c_i^\mm(t,\xxt) := \frac{1}{2}\sum_{k,j=1}^{m,d} \sigma_{jk}^\mm(t,\xxt)\partial_{\xxt_j} \sigma_{ik}^\mm(t,\xxt)$, $a_{ij}^\mu :=  \sum_{k=1}^m \sigma_{ik}^\mm\sigma_{jk}^\mm$, $i=1,\dots,\ell$, is the Stratonovich-corrected drift in (\ref{SDE1}). Analogous notation holds for $ {b}_i^\nu(t,x)$, $i=1,\dots,d$.

\vspace{.1cm}
\item[\tiny$\bullet$] $\sigma^\mm_{k}$, $\sigma^\nu_{k}$ stand for the $k$-th column of the matrix fields $\sigma^\mm$ and $\sigma^\nu$ with coefficients $\sigma^\mm_{ik}$, $\sigma^\nu_{jk}$.

\vspace{.1cm}
\item[\tiny$\bullet$] $\Vert \sigma(t,x)\Vert^2_{\textsc{hs}} :=\sum_{i=1}^n \sum_{k=1}^m\vert \sigma_{ik}(t,x)\vert^2$ is the Hilbert-Schmidt (or Frobenius) norm of a matrix field $\sigma$.   

\vspace{.1cm}
\item[\tiny$\bullet$] $\big\{\phi^\mm_{t_0,t}(\ccdot,\om)\,{:}\; s,t\in\Ic\big\}$, $\om\in \Om_\mm$, denotes a stochastic flow on $\XXt$ associated with the original dynamics (usually, but not exclusively, generated by the  SDE (\ref{SDE1})); see~\S\ref{SFlow}.   

\vspace{.1cm}
\item[\tiny$\bullet$]  $\big\{\phi^\nu_{t_0,t}(\ccdot,\om)\,{:}\; s,t\in\Ic\big\}$, $\om\in \Om_\nu$, denotes a stochastic flow on $\Xt\subseteq \XXt$ associated with the approximate dynamics (usually, but not exclusively, generated by the  SDE (\ref{SDE2})). 

\vspace{.1cm}
\item[\tiny$\bullet$] $t\mapsto \phi^\mm_{t_0,t}(\xxt,\om)$, $t\in \Ic$, $\phi^\mm_{t_0,t_0}(\xxt,\om)=\xxt\in \XXt$ is a random path of the original dynamical system on $\XXt$. 

\vspace{.1cm}
\item[\tiny$\bullet$] $t\mapsto \phi^\nu_{t_0,t}(x,\om)$, $t\in \Ic$, $\phi^\nu_{t_0,t_0}(x,\om)=x\in \Xt$, is a random path of the approximate dynamics on $\Xt\subseteq\XXt$.

\vspace{.1cm}
\item[\tiny$\bullet$] $\PP(\XXt)$ is a set of all probability measures on $\XXt$. $\PP(\Xt)$ is a set of all probability measures on $\Xt$.

\vspace{.1cm}
\item[\tiny$\bullet$] $\mm_t\in \PP(\XXt)$ is the time-marginal probability measure associated with the dynamics on $\XXt$.  For dynamics induced by the SDE (\ref{SDE1}), $\mm_t$ solves (weakly) the forward Kolmogorov equation (\ref{F_Kol}a) with  $\mathcal{L}_t^{\mm*}$. 

\vspace{.1cm}
\item[\tiny$\bullet$] $\varrho^\mm_t$ is the density of $\mm_t$ w.r.t.~Lebesgue measure $m_\ell$ on $\XXt$ (whenever $\mm_t\ll m_\ell$).

\vspace{.1cm}
\item[\tiny$\bullet$] $\mu_t\in \PP(\Xt)$ is obtained by a projection of $\mm_t$ via marginalisation of its Lebesgue density over $\XXt\setminus\Xt$.

\vspace{.1cm}
\item[\tiny$\bullet$] $\rho^\mu_t$ denotes a density of $\mu_t\in \PP(\Xt)$ obtained via marginalisation of $\varrho^\mm_t$ over $\XXt\setminus \Xt$. 

\vspace{.1cm}
\item[\tiny$\bullet$] $\nu_t\in \PP(\Xt)$ is the time-marginal probability measure induced by  the approximate dynamics on $\Xt\subseteq \XXt$. For dynamics of  the SDE~(\ref{SDE2}) $\nu_t$ solves (weakly) the forward Kolmogorov equation (\ref{F_Kol}b) with~$\mathcal{L}_t^{\nu*}$.

\vspace{.1cm}
\item[\tiny$\bullet$] $\rho^\nu_t$ denotes a density of $\nu_t\in \PP(\Xt)$ w.r.t.~the Lebesgue measure $m_d$ on $\Xt$ (whenever $\nu_t\ll m_d$).

\vspace{.1cm}
\item[\tiny$\bullet$] $(\mathcal{P}^\mm_{t_0,t})_{t\geqslant s}$ is a family of transition evolutions  induced by $\phi^\mm_{t_0,t}$ acting on $f\in \mathbb{M}(\XXt)$; see (\ref{P}).

\vspace{.1cm}
\item[\tiny$\bullet$]
 $(\mathcal{P}^{\mm*}_{t_0,t})_{t\geqslant s}$ is a family of duals of $\mathcal{P}^\mm_{t_0,t}$ acting on probability measures in $\PP(\XXt)$; see (\ref{P*}).  

\vspace{.1cm}
\item[\tiny$\bullet$] $(\mathcal{P}^\nu_{t_0,t})_{t\geqslant s}$ and $(\mathcal{P}^{\nu*}_{t_0,t})_{t\geqslant s}$  - analogous to the above except that they are induced by $\phi^\nu_{t_0,t}$.

\vspace{.1cm}
\item[\tiny$\bullet$] $\Ms^\mm_{t_0,\xxt}$ solves  the martingale problem for the operator $\mathcal{L}^\mm_t$ starting at $(t_0,\xxt)\in \Ic\times\XXt$ and is associated with the original dynamics.  The law of $\Ms^\mm_{t_0,\xxt}$ is identified with a path space probability measure on $(\mathcal{W}_\ell, \mathcal{B}(\mathcal{W}_\ell))$; the same symbol is used for the law of $\Ms^\mm_{t_0,\xxt}$ and the corresponding probability measure.  

\vspace{.1cm}
\item[\tiny$\bullet$] $\Ms^\mm_{t_0}(d\om):=\int_\XXt \Ms^\mm_{t_0,\xxt}(d\om)\mu_{t_0}(d\xxt)$ is a $\mm_{t_0}$-\,dependent  path space probability measure on $(\mathcal{W}_\ell, \mathcal{B}(\mathcal{W}_\ell))$. 

\vspace{.1cm}
\item[\tiny$\bullet$] $\Ns^\nu_{t_0,\xxt}$ is a solution to the martingale problem for the operator $\mathcal{L}^\nu_t$ starting at $(t_0,x)\in \Ic\times\Xt$, which is  associated with the approximate dynamics,  and identified with a  probability measure on $(\mathcal{W}_d, \mathcal{B}(\mathcal{W}_d))$. 

\vspace{.1cm}
\item[\tiny$\bullet$] $\Ns^\nu_{t_0}(d\om):=\int_\XXt \Ns^\nu_{t_0,x}(d\om)\nu_{t_0}(dx)$ is a  $\nu_{t_0}$-\,dependent path space probability measure on $(\mathcal{W}_d, \mathcal{B}(\mathcal{W}_d))$.

\vspace{.1cm}
\item[\tiny$\bullet$] $\E[f(\phi^\mm_{t_0,t}(\xxt))]:=\int_{\Om_\mm} f(\phi^\mm_{t_0,t}(\xxt,\om))\Ms^\mm_{t_0,\xxt}(d\om)$  denotes an ``observable" based on $f\in \mathbb{M}(\XXt)$, defined on the paths $t\mapsto\phi^\mm_{t_0,t}(\xxt,\om)$. $\E[f(\phi^\mm_{t_0,\cccdot}(\xxt))]$ is an observable for $f\in \mathbb{M}(\mathcal{W}_\ell)$.

\vspace{.1cm}
\item[\tiny$\bullet$] $\E[f(\pi^\nu_\mm\circ\phi^\mm_{t_0,t}(\xxt))]:=\int_{\Om_\mm} f(\pi^\nu_\mm\circ\phi^\mm_{t_0,t}(\xxt,\om))\Ms^\mm_{t_0,\xxt}(d\om)$, with $\pi^\nu_\mm:\XXt\rightarrow \Xt$ a projection,  is an observable based on $f\in \mathbb{M}(\Xt)$, defined on the paths $t\mapsto\phi^\mm_{t_0,t}(\xxt,\om)$ of the original dynamical system.

\vspace{.1cm}
\item[\tiny$\bullet$] $\E[f(\phi^\nu_{t_0,t}(x))]=\int_{\Om_\nu} f(\phi^\nu_{t_0,t}(x,\om))\Ns^\nu_{t_0,x}(d\om)$;  denotes an observable based on $f\in \mathbb{M}(\Xt)$, $\Xt\subseteq\XXt$, evaluated on the paths $t\mapsto\phi^\nu_{t_0,t}(x,\om)$.

\vspace{.1cm}
\item[\tiny$\bullet$]
$\E^{{\Ms}^\mm_{t_{0}}}\!\left[ f(\phi^\mm_{t_0,\cccdot})\right]:= \int_{\Om_\mm} f(\phi^\mm_{t_0,\cccdot}(\om))\Ms_{t_0}^\mm(d\om)=\int_{\Om_\mm}\int_\XXt f(\phi^\mm_{t_0,\ccdot}(\xxt,\om))\Ms^\mm_{t_0,\xxt}(d\om)\mm_{t_0}(d\xxt)$ denotes an observable based on $f\in \mathbb{M}(\mathcal{W}_\ell)$, and evaluated on  paths $t\mapsto\phi^\mm_{t_0,t}(\xxt,\om)$ of the original dynamical system. This is abbreviated as $\E^{{\Ms}^\mm_{t_{0}}}\!\left[ f(\phi^\mm_{t_0,\cccdot})\right]\equiv \E\left[ f(\phi^\mm_{t_0,\cccdot})\right]$ whenever not ambiguous.

\vspace{.1cm}
\item[\tiny$\bullet$]
$\E^{{\Ns}^\nu_{t_{0}}}\left[ f(\phi^\nu_{t_0,\ccdot})\right]:= \int_{\Om_\nu} f(\phi^\nu_{t_0,\cccdot}(\om))\Ns_{t_0}^\nu(d\om)=\int_{\Om_\nu}\int_\Xt f(\phi^\nu_{t_0,\cccdot}(x,\om))\Ns^\nu_{t_0,x}(d\om)\nu_{t_0}(dx)$ denotes an observable defined  on $f\in \mathbb{M}(\mathcal{W}_d)$, and evaluated on the paths $t\mapsto\phi^\nu_{t_0,t}(x,\om)$ of the approximate dynamics. This is abbreviated as $\E^{{\Ms}^\nu_{t_{0}}}\!\left[ f(\phi^\nu_{t_0,\cccdot})\right]\equiv \E\left[ f(\phi^\nu_{t_0,\cccdot})\right]$ whenever not ambiguous.

\vspace{.1cm}
\item[\tiny$\bullet$] $\Df\big(\mu\|\nu\big)$ is a $\varphi$-divergence between measures in $\PP(\Xt)$; $\varphi$  is a strictly convex function (see~\S\ref{phi_def}).

\vspace{.1cm}
\item[\tiny$\bullet$] $\Df\big(\Ms^\mu_{t_0}\|\Ns^\nu_{t_0}\big)$ is a $\varphi$-divergence between path space measures in  $\PP(\mathcal{W}_d)$ (only for $\Xt=\XXt$); see \S\ref{Path-space}.

\end{itemize}
}

\bigskip
\begin{center} \textsc{Acknowledgements}\end{center}

{\small The research of M.B.~was supported by the Office of Naval Research grant ONR N00014-15-1-2351 and ONRG N62909-20-1-2037. K.U.~was supported by the first grant as a postdoctoral research fellow. We are grateful to the anonymous referees for their thorough read and many constructive remarks which helped improve the paper.}

\bibliographystyle{plain}

\end{document}